\documentclass[12pt]{article}
\usepackage{amsfonts}
\usepackage{amssymb}
\usepackage{polski}

\input{tcilatex}
\begin{document}

\begin{center}
\textbf{The Steinhaus-Weil property: its converse, subcontinuity and Solecki
amenability}

\textbf{by\\[0pt]
N. H. Bingham and A. J. Ostaszewski}

\bigskip
\end{center}

\noindent \textbf{Abstract.}

The Steinhaus-Weil theorem that concerns us here is the simple, or
classical, `interior-points' property -- that in a Polish topological group
a non-negligible set $B$ has the identity as an interior point of$\ B^{-1}B.$
There are various converses; the one that mainly concerns us is due to
Simmons and Mospan. Here the group is locally compact, so we have a \textit{%
Haar} reference measure $\eta $. The Simmons-Mospan theorem states that a
(regular Borel) measure has such a Steinhaus-Weil property if and only if it
is absolutely continuous with respect to the Haar measure. In Part I
(Propositions 1-7, Theorems 1-4) we exploit the connection between the
interior-points property and a selective form of infinitesimal invariance
afforded by a certain family of \textit{selective} reference measures $%
\sigma $, drawing on Solecki's amenability at 1 (and using Fuller's notion
of subcontinuity). In Part II (Propositions 8, 9, Theorems 5, 6) we develop
a number of relatives of the Simmons-Mospan theorem. In Part III (Theorems
7, 8) we link this with topologies of Weil type. We close in Part IV with
Propositions 12-13 and Theorem B -- concerning the `composite
interior-point' property (of $AB^{-1}$) and the Borell `relative
interior-point' property (relative to the Cameron-Martin space) -- and
complements.

\bigskip

\noindent \textbf{Keywords. }Steinhaus-Weil property, amenability at 1,
measure subcontinuity, Simmons-Mospan theorem, Weil topology,
interior-points property, Haar measure, Lebesgue decomposition, left Haar
null, Cameron-Martin space.

\noindent \textbf{Classification}: Primary 22A10, 43A05; Secondary 28C10.

\bigskip

\begin{center}
\textbf{Table of Contents}
\end{center}

\noindent \hbox{          }\hbox{          }\hbox{          }0. Introduction%
\newline
\hbox{          }I. Subcontinuity and Solecki's amenability\newline
\qquad \hbox{          }\hbox{          }\hbox{          }1. Measure under
translation -- preliminaries\newline
\qquad \hbox{          }\hbox{          }\hbox{          }2. Subcontinuity
of measures\newline
II. Around the Simmons-Mospan Theorem\newline
\qquad \hbox{          }\hbox{          }\hbox{          }3. A Lebesgue
Decomposition\newline
\qquad \hbox{          }\hbox{          }\hbox{          }4. Discontinuity:
the Simmons-Mospan Theorem\newline
III. Weil-like topologies\newline
\qquad \hbox{          }\hbox{          }\hbox{          }5. Weil-like
topologies: preliminaries\newline
\qquad \hbox{          }\hbox{          }\hbox{          }6. Weil-like
topologies: theorems\newline
IV. Other interior-point properties\newline
\qquad \hbox{          }\hbox{          }\hbox{          }7.1. The Steinhaus
property $AA^{-1}$\ versus the Steinhaus property $AB^{-1}$\newline
\qquad \hbox{          }\hbox{          }\hbox{          }7.2. Borell's
interior-point property\newline
\qquad \hbox{          }\hbox{          }\hbox{          }8.\hbox{          }%
\hbox{          }\hbox{          } Complements\newline
\noindent References

\bigskip

\section*{0\qquad Introduction}

We begin by stating the Steinhaus-Weil Theorem in its simplest form
(Steinhaus [Ste] for the line, Weil [Wei, \S 11, p. 50] for a Polish locally
compact group, Grosse-Erdmann [GroE]):

\bigskip

\noindent \textbf{Theorem SW. } \textit{In a locally compact Polish group }$%
G $\textit{\ with (left) Haar measure }$\eta _{G}$\textit{, for non-null
Borel }$B$\textit{, }$B^{-1}B$ \textit{(and likewise }$BB^{-1}$\textit{)
contains a neighbourhood of the identity.}

\bigskip

The context we work in here and throughout, unless otherwise stated, is that
groups and spaces are assumed separable. This both simplifies the exposition
and emphasizes that we need only the axiom of Dependant Choices (DC -- `what
is needed to make induction work'), rather than the Axiom of Choice (AC);
cf. [BinO7]. For comments concerning non-separable settings, see \S 8.1.

The interior-point property of the measure-theoretically `non-negligible'
set $B$ of the theorem is referred to as the \textit{Steinhaus-Weil property}%
, which encompasses the category variant due to Piccard [Pic] and Pettis
[Pet], cf. Cor. 2$^{\prime }$ and Th. 7B (by reference, when appropriate, to
the \textit{quasi-interior} of a set -- the largest open set equivalent to
it modulo a meagre set). This important result has many ramifications; for
example, it is basic to the theory of regular variation -- see e.g. [BinGT,
Th. 1.1.1].

We are also concerned here with its converse, the Simmons-Mospan theorem
([Sim, Th. 1], [Mos, Th. 7], recently rediscovered in the abelian case
[BarFF, Th. 10]), stated in its simplest form.

\bigskip

\noindent \textbf{Theorem SM}. \textit{In a locally compact Polish group, a
Borel measure has the Steinhaus-Weil property if and only if it is
absolutely continuous with respect to Haar measure}.

\bigskip

\noindent The proof of this and related (converse) results emerges in \S 4
after due preparatory work; regarding the reduction of the non-separable to
the separable case, see the closing remark of \S 4. The results below hinge
on work of Solecki [Sol2] on amenability at 1 and the concept of
subcontinuity (see \S 2 and below). These are aimed at freeing up the
classical dependency on local compactness and the corresponding standard
(Haar) reference measure. To the best of our knowledge such aims, in respect
of topological groups, were last undertaken by Xia in 1972 in Chapter 3 of
[Xia], where the emphasis is on (relative) quasi-invariance (cf. \S 7.2), a
topic we pursue in the companion paper [BinO8] with tools developed here.
(The additive subgroup context of the more recent [Bana] sits in a
topological vector space.)

For $G$ a topological group with (admissible) metric $d$ (briefly: metric
group), denote by $\mathcal{M}(G)$ the family of Borel regular $\sigma $%
-finite measures on $G,$ with $\mathcal{P}(G)\subseteq \mathcal{M}(G)$ the
probability measures ([Kec, \S 17E], [Par]), by $\mathcal{P}_{\text{fin}}(G)$
the larger family of finitely-additive regular probability measures (cf.
[Bin]), and by $\mathcal{M}_{\text{sub}}(G)$ submeasures (monotone, finitely
subadditive set functions $\mu $ with $\mu (\emptyset )=0$). Here \textit{%
regular} is taken to imply both \textit{inner} regularity (inner
approximation by compact subsets, also called the \textit{Radon} property,
as in [Bog2, II \S 7.1] and [Sch]), and \textit{outer} regularity (outer
approximation by open sets). We recall that a $\sigma $-finite Borel measure
on a metric space is necessarily outer regular ([Bog2, II. Th. 7.1.7], [Kal,
Lemma 1.34], cf. [Par, Th. II.1.2] albeit for a probability measure) and,
when the metric space is complete, inner regular ([Bog2, II. Th. 7.1.7], cf.
[Par, Ths. II.3.1 and 3.2]). When $G\ $is locally compact we denote Haar
measure by $\eta _{G}$ or just $\eta $ (H denoting capital eta in Greek).
For $X$ a metric space, we denote by $\mathcal{K}=\mathcal{K}(X)$ the family
of compact subsets of $X$ (the \textit{hyperspace }of $X$ in \S 1,\ where we
view it as a topological space under the Hausdorff metric, or the Vietoris
topology). For $\mu \in \mathcal{M}(G)$ we write $_{g}\mu (\cdot ):=\mu
(g\cdot )$ and $\mu _{g}(\cdot ):=\mu (\cdot g)$; $\mathcal{M}(\mu )$
denotes the $\mu $-measurable sets of $G$ and $\mathcal{M}_{+}(\mu )$ those
of finite positive measure, and $\mathcal{K}_{+}(\mu ):=\mathcal{K}(G)\cap
\mathcal{M}_{+}(\mu )$. For $G$ a Polish group, recall that $E\subseteq G$
is \textit{universally measurable} ($E\in \mathcal{U}(G)$) if $E$ is
measurable with respect to every measure $\mu \in \mathcal{P}(G)$ -- for
background, see e.g. [Kec, \S 21D], cf. [Fre, 434D, 432], [Sho]; these form
a $\sigma $\textit{-algebra}. Examples are analytic subsets (see e.g. [Rog,
Part 1 \S 2.9], or [Kec, Th. 21.10], [Fre, 434Dc]) and the $\sigma $-algebra
that they generate. Beyond these are the \textit{provably }$\mathbf{\Delta }%
_{2}^{1}$ sets of [FenN] -- cf. [BinO7].

Recall that $E$ is \textit{left Haar null}, $E\in \mathcal{HN}$, as in
Solecki [Sol1,2,3] (following [Chr1,2]) if there are $B\in \mathcal{U}(G)$
covering $E$ and $\mu \in \mathcal{P}(G)$ with
\[
\mu (gB)=0\qquad (g\in G).
\]%
(The terminal brackets here and below indicate \textit{universal
quantification} over the free variable.) So if $B\in \mathcal{U}(G)$ is not
left Haar null, then for each $\mu \in \mathcal{P}(G)$ there is compact $%
K=K_{\mu }\subseteq B$ and $g\in G$ with%
\[
_{g}\mu (K)>0.
\]%
The question then arises whether there is also $\delta >0$ with $_{g}\mu
(Kt)>0$ for \textit{all} $t\in B_{\delta },$ for $B_{\delta }=B_{\delta
}(1_{G})$ the open $\delta $-ball centered at $1_{G}:$ a \textit{right-sided
}property\textit{\ complementing\ }the earlier\textit{\ left-sided} property
(of nullity, or otherwise). If this is the case for some $\mu ,$ then (see
Corollary 2$^{\prime }$ in \S 2) $1_{G}\in \mathrm{int}(K^{-1}K)\subseteq
\mathrm{int}(E^{-1}E)$; indeed, one has%
\begin{equation}
K\cap Kt\in \mathcal{M}_{+}(_{g}\mu )\qquad (t\in B_{\delta }),
\tag{$\ast
M$}
\end{equation}%
and this implies (see Lemma 1, \S 2):%
\[
B_{\delta }\subseteq \mathrm{int}(K^{-1}K)\subseteq \mathrm{int}(E^{-1}E).
\]%
As this clearly forces local-compactness of $G$ (see Lemma 1 below), for the
more general context we weaken the `complementing right-sided property' to
hold only \textit{selectively}: on a subset of $B_{\delta }$ of the form%
\[
\{z\in B_{\delta }:|\mu (Kz)-\mu (K)|<\varepsilon \}
\]%
(cf. $B_{\delta }^{\Delta }(\mu )$ in \S 2). This turns attention to
refining the topology $\mathcal{T}_{d}$ of $G$ via the metrics%
\[
d_{K}(x,y):=d(x,y)+|\mu (Kx)-\mu (Ky)|
\]%
(for a family of sets $K\in \mathcal{K}_{+}(\mu )$ -- cf. the Struble
sampler of \S 5), determined in Theorem 5 below, aiming at a
relative-interior-points property in the finer topology, the theme of the
companion paper [BinO8], cf. \S 8.2. See (in addition to Theorem SW above)
Kemperman [Kem] (cf. [Kuc, Lemma 3.7.2], [BinO1, Th. K], [BinO5, Th. 1(iv)]).

We are guided by the close relation between the measure-theoretic \textit{%
Steinhaus-Weil-like} property (*M) and its category version
\begin{equation}
K\cap Kt\in \mathcal{B}_{+}(\tau ),  \tag{$\ast B$}
\end{equation}%
where the latter term refers to non-meagre \textit{Baire} sets (= with the
Baire property) of $\tau $, a refinement of the ambient topology $\mathcal{T}%
_{G}$ $=\mathcal{T}_{d}$ of $G,$ the latter conveniently taken to be
generated by a \textit{left}-invariant metric $d=d_{L}^{G}$ with associated
group-norm (\S 5) $||x||:=d(x,1_{G})=d(tx,t)$ (so that $B_{\delta
}(t)=tB_{\delta }$ -- see Prop. 1). We refer to the (left) invariance of $%
\mathcal{B}_{+}(\tau )$ (under translation) as the (left) \textit{Nikodym
property} of $\tau .$

In Part I below (\S 1,2, Props. 1-7, Th. 1-4: Subcontinuity Theorem,
Aggregation Theorem, Shift-compactness Theorem, Strong Subcontinuity
Theorem), in the context of a metric or Polish group $G,$ we study
continuity properties of the maps $m_{K}:t\mapsto \mu (Kt)$ in the light of
theorems of Solecki [Sol2] and of Theorem SM above and related results. The
key here is Fuller's notion of subcontinuity, as applied to the function $%
m_{K}(t)$ at $t=1_{G}$. This yields a fruitful interpretation of Solecki's
notion of \textit{amenability at }$1_{G}$ via \textit{selective subcontinuity%
} and linkage to \textit{shift-compactness} (see Th. 3 below; the term is
borrowed from [Par, III.2]). Since commutative Polish groups are amenable at
1 [Sol2, Th. 1(ii)], this widens the field of applicability of
shift-compactness to non-Haar-null subsets of these, as in [BinO4], and
leads to a conjecture (see remarks preceding Theorem 3) as to whether $%
\mathcal{HN}$ comprises the negligible sets of some refinement topology of $%
\mathcal{T}_{d}$. In Part II (\S 3,4, Props. 8-9, Theorems 5,6:
Disaggregation Theorem, Generalized Simmons Theorem) we study measure
discontinuity and extend results of Simmons beyond his locally-compact
context by reference to measures exhibiting selective subcontinuity. This
draws on results in the companion paper [BinO8], which is concerned with the
relative-interior-point property inspired by the \textit{Cameron-Martin}
theory of Gaussian measures (\S 7.2 and \S 8.2). In Part III (\S 6, 7,
Props. 10, 11, Theorems 7, 7M, 7B, 8) we study the converse problem of the
\textit{Weil topologies} generated by the \textit{Fr\'{e}chet-Nikodym}
pseudo-norms, defined for $E\in \mathcal{M}_{+}(G)$ by:%
\[
||t||_{\mu }^{E}:=\mu (tE\triangle E)\qquad (t\in G).
\]%
For ease of reference -- duly identified, as here -- we study $EE^{-1}$ (for
$||t||_{\mu }^{E}>0$) when following Simmons and Weil, and $E^{-1}E$ when
following Solecki. In Part IV (\S 7, Props. 12-13, Th. K, Th. B) we study
the relation of the $AA^{-1}$ interior-point property to that for $AB^{-1}$,
and discuss the \textit{Borell interior-point property} relative to the
Cameron-Martin `embedded subspace'. We close with a section of complements (%
\S 8).

\part{Subcontinuity and Solecki's amenability}

\section{Measure under translation -- preliminaries}

We begin with a form of the `telescope' or `tube' lemma (cf. [Mun, Lemma
5.8]), applied in \S 2. Our usage of upper semicontinuity in relation to
set-valued maps follows [Rog], cf. [Bor].

\bigskip

\noindent \textbf{Proposition 1 }(cf. [Hey, 1.2.8]). \textit{For a metric
group }$G$ \textit{and compact }$K\subseteq G,$ \textit{the map }$t\mapsto
Kt $ \textit{is upper semicontinuous; in particular, for }$\mu \in \mathcal{M%
}(G)$\textit{, }%
\[
m_{K}:t\mapsto \mu (Kt)
\]%
\textit{is upper semicontinuous, hence }$\mu $\textit{-measurable. In
particular, if }$m_{K}(t)=0$\textit{, then }$m$ \textit{is continuous at }$t$%
\textit{.}

\bigskip

\noindent \textbf{Proof. }For $K$ compact and $V\supseteq K$ open, pick for
each $k\in K$ a $\delta (k)>0$ with $kB_{2\delta (k)}$ $\subseteq V.$ By
compactness, there are $k_{1},...,k_{n}$ with $K\subseteq
\bigcup\nolimits_{j}k_{j}B_{\delta (k_{k})}\subseteq V;$ then for $\delta
:=\min_{j}r(k_{j})>0$%
\[
Kt\subseteq \bigcup\nolimits_{j}k_{j}B_{\delta (k_{k})}t\subseteq
\bigcup\nolimits_{j}k_{j}B_{2\delta (k_{k})}\subseteq V\qquad (||t||<\delta
).
\]

To prove upper semicontinuity of $m_{K},$ fix $t\in G.$ For $\varepsilon >0$%
, as $Kt$ is compact, choose by outer regularity an open $U\supseteq Kt$
with $\mu (U)<\mu (Kt)+\varepsilon ;$ as before, there is an open ball $%
B_{\delta }$ at $1_{G}$ with $KtB_{\delta }\subseteq U$, and then $\mu
(KtB_{\delta })\leq \mu (U)<\mu (Kt)+\varepsilon .$ The final assertion
follows from positivity of $m_{K}.$ $\square $

\bigskip

We continue with an analogue. The result is folklore, cf. [BeeV, Th.
3.2(i)]; it comes close to matters touched on in [Ost1, \S 3]. Here and
below the vertical section of a set $A$ is denoted $A_{x}:=\{y:(x,y)\in A\}.$

\bigskip

\noindent \textbf{Proposition 2 (Sectional upper semicontinuity).} \textit{%
For a metric group }$G,\,$\textit{compact }$F\subseteq G$\textit{\ and
compact }$K\subseteq G^{2}$\textit{, the map}%
\[
x\mapsto K_{x}\qquad (x\in F)
\]%
\textit{is upper semicontinuous.}

\bigskip

\noindent \textbf{Proof.} For $V\subseteq G$ open with $K_{x}\subseteq V,$
suppose for $x_{n}\in F$ with $x_{n}\rightarrow x$ that $(x_{n},y_{n})\in
K\backslash (G\times V).$ By compactness of $K,$ we may suppose w.l.o.g.
that $y_{n}\rightarrow y.$ Then $(x,y)\in K\backslash (G\times V),$ and so $%
(x,y)\in \{x\}\times K_{x}$ and $y\notin V;$ but $y\in K_{x}\subseteq V,$ a
contradiction. $\square $

\bigskip

From Prop. 2 on\textit{\ upper} semicontinuity, we obtain information about $%
m_{K}:t\mapsto \mu (Kt)$ below. This links with \textit{lower}
semicontinuity. By a theorem of Fort, the $\varepsilon $-continuity points
(defined in terms of the Hausdorff metric: see [For]) of an upper
semicontinuous compact-valued mapping of a metric space into a totally
bounded metric space form a dense open set, implying in a real-valued
context such as here \textit{continuity} on a co-meagre set. We return to
this shortly in Theorem LB below.

\bigskip

\noindent \textbf{Proposition 3 (Sectional upper semicontinuity under a
measure).} \textit{For a metric group }$G$,\textit{\ compact }$F\subseteq G$
\textit{and compact }$K\subseteq G^{2}$\textit{,\ and }$\mu \in \mathcal{M}%
(G)$\textit{, the map}%
\[
m:x\mapsto \mu (K_{x})\qquad (x\in F=\mathrm{proj}_{1}K)
\]%
\textit{is upper semicontinuous, and so Borel.}

\bigskip

\noindent \textbf{Proof.} Fix $x\in F.$ Let $\varepsilon >0.$ By outer
regularity, take $V$ open in $G\ $with $K_{x}\subseteq V$ and $\mu (V)<\mu
(K_{x})+\varepsilon .$ By Prop. 2 $x\mapsto \{x\}\times K_{x}$ is upper
semicontinuous on $F$; so for some open neighbourhood $U$ of $x$%
\[
K\cap (U\times G)\subseteq K\cap (U\times V).
\]%
So for $y\in F\cap U$%
\[
K_{y}\subseteq V,
\]%
and so $\mu (K_{y})\leq \mu (V)<\mu (K_{x})+\varepsilon ,$ proving the first
assertion. The second assertion follows since%
\[
m^{-1}(a,b)=\bigcap\nolimits_{n\in \mathbb{N}}m^{-1}[0,b)\backslash
m^{-1}[0,a+1/n).\qquad \square
\]

\bigskip

For further results on Borel-measurability of regular Borel measures see
[BeeV, Th. 2.2] (there termed `Radon measures').

We will need the following result in Parts II and III (see Lemma 2, \S 4,
and Th. 7, \S 6), preferable to the usual Fubini Theorem as using
qualitative rather than quantitative measure theory (like the
Kuratowski-Ulam Theorem [FreNR]). Interestingly, it may be proved by
mimicking the proof of Prop. 1 above, yielding a simplification to that by
Eric van Dowen [vDo], itself a simplification of that in [Oxt2, Ch. 14]: for
the proof (omitted here), see \S 8.12.

\bigskip

\noindent \textbf{Theorem FN (Fubini theorem for null sets).} \textit{For a
metric group }$G$\textit{and} $A\subseteq G^{2}$ \textit{measurable under }$%
\mu \times \nu $\textit{, with }$\mu ,\nu \in \mathcal{M}(G)$\textit{: if
the `exceptional set' of points }$x $\textit{\ for which the vertical
section }$A_{x}$\textit{\ is }$\nu $\textit{-non-null is itself }$\mu $%
\textit{-null, then }$A$ \textit{is }$\mu \times \nu $\textit{-null.}

\bigskip

We close this section with a study of the continuity properties of the map $%
m_{K}:t\mapsto \mu (Kt)$ for compact $K,$ extending Prop. 3.

\bigskip

\noindent \textbf{Corollary 1 (}Fort [For]). \textit{In Proposition 2, }$%
t\mapsto \mu (Kt)$\textit{\ is lower semi-continuous (so also continuous) on
a co-meagre set.}

\bigskip

We can improve on the preceding result by recourse to a natural
generalization, for our compact sectional context, of the classical
continuity theorems of Luzin [Hal, \S 55] and Baire [Oxt2, Th. 8.1] -- see
also [Sch, Ch. 1, \S 5]. Below for a \textit{compact} metric space $X$, we
denote by $\mathcal{K}(X)$ the \textit{hyperspace of }$X,$ the space of
compact subsets of $X$ under the Hausdorff metric, or Vietoris topology;
here this is also a compact space ([Eng, 2.7.28], [Kec, Th. 4.25], [Mic]).
Then (LB for `Luzin-Baire'):

\bigskip

\noindent \textbf{Theorem LB. }\textit{For }$G$ \textit{a metric group and
compact }$K\subseteq G^{2}$, \textit{the map} $\kappa :G\rightarrow \mathcal{%
K}(G):x\mapsto K_{x}$ \textit{is Borel-measurable, and so}

\noindent (i) $\kappa $ \textit{is continuous relative to a co-meagre set}.

\textit{For }$\mu \in \mathcal{P}(G)$\textit{:}

\noindent (ii)\textit{\ for each }$\varepsilon >0$\textit{\ there is a Borel
set }$S_{\varepsilon }$\textit{\ with }$\mu (G\backslash S_{\varepsilon
})<\varepsilon $\textit{\ such that }$x\mapsto K_{x}$\textit{\ is continuous
on }$S_{\varepsilon }$\textit{; equivalently:}

\noindent (ii)$^{\prime }$\textit{\ there is an increasing sequence of Borel
sets }$S_{n}$ \textit{with union }$\mu $\textit{-almost all of }$G$ \textit{%
such that }$x\mapsto K_{x}$\textit{\ is continuous on each }$S_{n}.$

\bigskip

\noindent \textbf{Proof. }For $U,V$ open, the set $\{x:K_{x}\subseteq U$ and
$K_{x}\cap V\neq \emptyset \}$ is $\sigma $-compact; indeed, by Prop. 2, $%
\{x:K_{x}\subseteq U\}$ is open, whereas for $F:=\mathrm{proj}_{1}(K),$ a
compact set
\[
\{x:K_{x}\cap U\neq \emptyset \}=\mathrm{proj}_{1}\left[ (F\times U)\cap K)%
\right] ,
\]%
which is the projection of a $\sigma $-compact set. Then (ii) follows from
Luzin's theorem ([Hal, \S 55]), and for (ii)$^{\prime }$ see [BinO1, p.
142]. Its extension to a regular (i.a. $\mathcal{G}$-outer regular) $\sigma $%
-finite measure may be made via Egoroff's Theorem [Hal, \S 21 Th. A] -- see
[Zak]. $\square $

\bigskip

A first corollary is the following result on the continuity of the map%
\[
x\mapsto ||x||_{E}^{\mu }=\mu (xE\triangle E),
\]%
for \textit{measurable} $E,$ by compact approximation. Below, the sets $%
C_{x} $ associated with points $x$ should be interpreted as neighbourhoods
of $x$ in the spirit of a Hashimoto ideal topology for the ideal of $\mu $%
-null sets, for which see [LukMZ], or [BinO5]. This mimicks Weil's \textit{%
proof} of the `fragmentation lemma' (\S 8, Lemma 3) in [Hal, Ch. XII \S 62
Th. A] (cf. [Wei, Ch. VII, \S 31]).

\bigskip

\noindent \textbf{Proposition 4 (Almost everywhere continuity)}. \textit{For
a metric group }$G$\textit{, }$\delta >0$\textit{, }$\mu \in \mathcal{P}(G),$%
\textit{\ }$E\in \mathcal{M}_{+}(\mu ),$\textit{\ and }$F\in \mathcal{K}%
_{+}(\mu )$\textit{:}

\textit{there is a compact }$C\subseteq F$ \textit{with }$\mu (F\backslash
C)<\delta $ \textit{such that for any }$\varepsilon >0$\textit{\ and each }$%
x\in C$ \textit{there is a }$\mu $-\textit{non-null measurable }$%
C_{x}\subseteq C$\textit{\ containing }$x$ \textit{with}%
\[
|\mu (xE\triangle E)-\mu (yE\triangle E)|<\varepsilon \qquad (y\in C_{x}).
\]%
\textit{In particular, there is an increasing family of compact sets }$C_{n}$%
\textit{\ with union }$\mu $\textit{-almost all of }$G$\textit{\ satisfying
the above with }$C_{n}$\textit{\ for }$C$\textit{. }

\bigskip

\noindent \textbf{Proof. }Fix $E\subseteq G$ measurable and for $F$ compact
with $\mu (F)>0,$ put%
\[
H:=\bigcup\nolimits_{x\in F}\{x\}\times (E\triangle
xE)=\bigcup\nolimits_{x\in F}\{x\}\times ((xE\backslash E)\cup (E\backslash
xE))
\]%
\[
=\bigcup\nolimits_{x\in F}(\{x\}\times xE)\backslash \bigcup\nolimits_{x\in
F}(\{x\}\times E)\cup \bigcup\nolimits_{x\in F}(\{x\}\times E)\backslash
\bigcup\nolimits_{x\in F}\{x\}\times (xE),
\]%
which is measurable. So $F\subseteq \mathrm{proj}_{1}(H)$ and for $x\in F$%
\[
H_{x}:=xE\triangle E.
\]

We will work inductively, taking successively smaller values of $\varepsilon
.$ Fix $\varepsilon >0$ with $\varepsilon <\mu (F)$ and choose $K$ a finite
union of compact rectangles with%
\[
(\mu \times \mu )(H\triangle K)<\varepsilon ^{2}.
\]%
So%
\[
K=\bigcup\nolimits_{j\leq n}(F^{j}\times K^{j}),
\]%
say. Let $S=S_{\varepsilon }:=\{x\in F:$ $\mu ((H\triangle K)_{x})\geq
\varepsilon \}.$ Then $\mu (S)\leq \varepsilon ;$ otherwise%
\[
\varepsilon ^{2}>(\mu \times \mu )(H\triangle K)\geq \varepsilon \mu (S)\geq
\varepsilon ^{2},
\]%
a contradiction. So%
\[
\mu (F\backslash S)>\mu (F)-\varepsilon .
\]%
So we may choose a compact set $C=C_{\varepsilon }\subseteq F\backslash
S_{\varepsilon }$ with $\mu (F\backslash C)<\varepsilon $ and%
\[
\mu (H_{x}\triangle K_{x})\leq \varepsilon \qquad (x\in C_{\varepsilon }).
\]%
But $K_{x}\in \{K^{1}...,K^{n}\},$ so%
\[
C_{\varepsilon }=\bigcup\nolimits_{j\leq n}C_{\varepsilon }^{j},\qquad \text{%
where }C_{\varepsilon }^{j}:=\{x\in C_{\varepsilon }:\mu (H_{x}\triangle
K^{j})\leq \varepsilon \}.
\]%
Now for $x,y\in C_{\varepsilon }^{j}$%
\[
\mu (H_{x}\triangle H_{y})\leq \mu (H_{x}\triangle K^{j})+\mu
(K^{j}\triangle H_{y})\leq 2\varepsilon ,
\]%
and w.l.o.g. the sets $C_{\varepsilon }^{j}$ may be assumed compact.

Now fix $\delta >0$ and repeat the construction inductively, taking in turn $%
\varepsilon =\varepsilon _{n}=2^{-n}\mu (F)\delta ,$ to obtain a sequence of
compact sets $F=C^{0}\supseteq C=C^{1}\supseteq C^{2}\supseteq C^{3}...$
with $C^{n}=C_{\varepsilon _{n}}$ as above, and $\mu (C^{n}\backslash
C^{n+1})<\varepsilon _{n}.$ These have non-null intersection $C_{0}:$
\[
\mu (F\backslash C_{0})<\sum\nolimits_{n\geq 1}2^{-n}\mu (F)\delta =\mu
(F)\delta .
\]%
Now for each $\varepsilon >0$ there is $n$ with $\varepsilon
_{n}<\varepsilon .$ So for $x\in C_{0}\subseteq C^{n}=C_{\varepsilon
_{n}}=\bigcup\nolimits_{j\leq n}C_{\varepsilon _{n}}^{j}$%
\[
\mu (H_{x}\triangle H_{y})<\mu (H_{x}\triangle H_{y})+2\varepsilon _{n}<\mu
(H_{x})+2\varepsilon \qquad (x,y\in C_{\varepsilon _{n}}^{j}\cap C_{0}).
\]%
So%
\[
|\mu (H_{x})-\mu (H_{y})|<4\varepsilon \qquad (x,y\in C_{\varepsilon
_{n}}^{j}\cap C_{0}).\qquad \square
\]

\bigskip

A proof similar to but simpler than that above (omitted here -- see \S 8.13)
improves Prop. 1:

\bigskip

\noindent \textbf{Proposition 5 (Almost everywhere upper semicontinuity)}.
\textit{For a metric group }$G$\textit{, }$\delta >0$\textit{, }$\mu \in
\mathcal{P}(G),$\textit{\ }$E\in \mathcal{M}_{+}(\mu ),$\textit{\ and }$F\in
\mathcal{K}_{+}(\mu )$\textit{:}

\textit{there is a compact }$C\subseteq F$ \textit{with }$\mu (F\backslash
C)<\delta $ \textit{such that for any }$\varepsilon >0$\textit{\ each }$x\in
C$ \textit{has a neighbourhood }$U_{x}$\textit{\ with}%
\[
\mu (yE)<\mu (xE)+\varepsilon \qquad (y\in C\cap U_{x}).
\]%
\textit{In particular, there are disjoint compact sets }$C$\textit{\ with
union }$\mu $\textit{-almost all of }$G$\textit{\ for which this holds}.

\section{Subcontinuity of measures}

Proposition 1 above, on upper semicontinuity, motivates the following
definitions, the key one being an adaptation of \textit{subcontinuity} (of
functions) due to Fuller [Ful] (for which see Remark 4 below) to the context
of measures. We focus on the \textit{right-sided} version of the concept.
Subcontinuity is a natural auxiliary in the quest for fuller forms of
continuity: as one instance, see [Bou] for the step from separate to joint
continuity; as another, classic instance, note that a subcontinuous
set-valued map with closed graph (yet another relative of upper
semicontinuity) is continuous -- see [HolN] for an extensive bibliography.
Here its relevance to the Steinhaus-Weil Theorem (which seems to be new
here) yields Theorems 1 and 3, linking \textit{amenability at 1} with
\textit{shift-compactness,} for which see Theorem 3 below (the latter term
is borrowed from [Par, III.2]).

\bigskip

\noindent \textbf{Definition }([BinO5])\textbf{. }For $\mu \in \mathcal{P}_{%
\text{fin}}(G)$, and (compact) $K\in \mathcal{K}(G),$ noting that $\mu
_{\delta }(K):=$ $\inf \{\mu (Kt):t\in B_{\delta }\}$ is weakly decreasing
in $\delta $, put%
\[
\mu _{-}(K):=\sup_{\delta >0}\inf \{\mu (Kt):t\in B_{\delta }\},
\]%
and, for $\mathbf{t}=\{t_{n}\}$ a \textit{null sequence}, i.e. with $%
t_{n}\rightarrow 1_{G},$
\[
\mu _{-}^{\mathbf{t}}(K):=\lim \inf\nolimits_{n\rightarrow \infty }\mu
(Kt_{n}).
\]%
Then%
\[
0\leq \mu _{-}(K)\leq \mu (K)=\inf_{\delta >0}\sup \{\mu (Kt):t\in B_{\delta
}\},
\]%
by Proposition 1. We say that a null sequence $\mathbf{t}$ is \textit{%
non-trivial} if $t_{n}\neq 1_{G}$ infinitely often. Define (with the
quantification convention of the Introduction) as follows:

\noindent (i) $\mu $ is \textit{translation-continuous} (`\textit{continuous}%
' or `\textit{mobile}') if $\mu (K)=\mu _{-}(K)$ ($K\in \mathcal{K}(G)$)$;$

\noindent (ii) $\mu $ is \textit{maximally discontinuous} at $K\in \mathcal{K%
}(G)$ if $0=\mu _{-}(K)<\mu (K);$

\noindent (iii) $\mu $ is \textit{subcontinuous} if $0<\mu _{-}(K)\leq \mu
(K)$ $(K\in \mathcal{K}_{+}(\mu));$

\noindent (iv) $\mu $ is \textit{(selectively) subcontinuous at }$K\in%
\mathcal{K}_{+}(\mu)$ \textit{along }$\mathbf{t}$ if $\mu _{-}^{\mathbf{t}%
}(K)>0$.

\bigskip

\noindent \textbf{Remarks. }1. $m_{K}(.)$ is \textit{continuous} if $\mu $
is continuous, since $m_{K}(st)=m_{Ks}(t)$ and $Ks$ is compact whenever $K$
is compact; for directional continuity of measures in linear spaces see
[Bog3, \S 3.1]. In [LiuR] (cf. [LiuRW], [Gow1,2]) a Radon measure $\mu $ on
a space $X,$ on which a group $G\ $acts homeomorphically, is called \textit{%
mobile} if $t\mapsto \mu (Kt)$ is continuous for all $K\in \mathcal{K}(X).$

\noindent 2. For $G$ locally compact (i) holds for $\mu $ the left Haar
measure $\eta _{G}$, and also for $\mu \ll \eta _{G}$ (absolutely continuous
w.r.t. to $\eta _{G}$).

\noindent 3. A measure $\mu $ singular w.r.t. Haar measure is maximally
discontinuous for its support: this is at the heart of the analysis offered
by Simmons (and independently, much later by Mospan) -- see Corollary 2$%
^{\prime }$ below.

\noindent 4. \textit{Subcontinuity}, in the sense of [Ful], of a map $%
f:G\rightarrow (0,\infty )$ requires that, for every $t_{n}\rightarrow t\in
G,$ there is a subsequence $t_{m(n)}$ with $f(t_{m(n)})$ convergent in the
range (i.e. to a positive value). The distinguished role of null sequences
emerges below in the \textit{Subcontinuity Theorem} (Theorem 1). Null
sequences should be viewed here as selecting stepwise (or even pathwise,
under local connectedness, as suggested by Tomasz Natkaniec) `asymptotic
directions' justifying the phrase `\textit{along} $\mathbf{t}$' in (iv)
above, and allowing (iv) to be interpreted as a \textit{selective
subcontinuity} in `direction' $\mathbf{t}$. The analogous selective concept
in a linear space is `along a vector' as in [Bog3, \S 3.1].

\noindent 5. \textit{Selective versus uniform subcontinuity.} Definition
(iii) is equivalent to demanding for $K\in \mathcal{K}_{+}(\mu )$ that
\textit{any} null sequence $\mathbf{t}=\{t_{n}\}$ have a subsequence $\mu
(Kt_{m(n)})$ bounded away from $0;$ then (iii) may be viewed as demanding
`uniform subcontinuity': selective subcontinuity along \textit{each} $%
\mathbf{t}$ for all $K\in \mathcal{K}_{+}(\mu )$.

\noindent 6. \textit{Left- versus right-sided versions. }Writing $\tilde{\mu}%
(E):=\mu (E^{-1})$ for the \textit{inverse measure} captures versions
associated with right-sided translation such as $\tilde{\mu}_{\_}$ and
\[
\tilde{\mu}_{-}^{\mathbf{t}}(K):=\lim \inf\nolimits_{n\rightarrow \infty
}\mu (t_{n}K).
\]

\noindent \textbf{Definition. }We will say that $\mu $ is \textit{symmetric}
if $\mu =\tilde{\mu};$ then $B$ is null iff $B^{-1}$ is null.

\bigskip

In Lemma 1 below it suffices for $\mu $ to be a bounded, regular submeasure
which is \textit{supermodular:}
\[
\mu (E\cup F)\geq \mu (E)+\mu (F)-\mu (E\cap F);
\]%
recall, however, from [Bog2, 1.12.37] the opportunity to replace, for any $%
K\in \mathcal{K}(G),$ a supermodular submeasure $\mu $ by a dominating $\mu
^{\prime }\in \mathcal{M}_{\text{fin}}(G),$ i.e. with $\mu ^{\prime }(K)\geq
\mu (K)$.

For $K\in \mathcal{K}_{+}(\mu )$ and $\delta ,\Delta >0,$ put%
\[
B_{\delta }^{\Delta }=B_{\delta }^{K,\Delta }(\mu ):=\{z\in B_{\delta }:\mu
(Kz)>\Delta \},
\]%
which is monotonic in $\Delta :$ $B_{\delta }^{\Delta }\subseteq B_{\delta
}^{\Delta ^{\prime }}$ for $0<\Delta ^{\prime }\leq \Delta .$ Note that $%
1_{G}\in B_{\delta }^{\Delta }$ for $0<\Delta <\mu (K).$

\bigskip

The specialization below to a mobile measure (see above) may be found in
[Gow1,2].

\bigskip

\noindent \textbf{Lemma 1 }(cf. [BinO5, Th. 2.5])\textbf{.} \textit{Let }$%
\mu \in \mathcal{P}_{\text{fin}}(G)$\textit{\ for }$G$ \textit{\ a metric
group. For }$K\in \mathcal{K}_{+}(\mu ),$ \textit{if }$\mu _{-}^{\mathbf{t}%
}(K)>0$ \textit{for some non-trivial null sequence }$\mathbf{t}$\textit{,\
then there are }$\delta >0$ \textit{and }$0<\Delta <\mu _{-}^{\mathbf{t}}(K)$%
\textit{\ with }$t_{n}\in B_{\delta }^{\Delta }$\textit{\ for all large
enough }$n$\textit{\ and}%
\[
\Delta \leq \mu (K\cap Kt)\qquad (t\in B_{\delta }^{\Delta }),
\]%
\textit{so that}%
\begin{equation}
K\cap Kt\in \mathcal{M}_{+}(\mu )\qquad (t\in B_{\delta }^{\Delta }).
\tag{$\ast $}
\end{equation}%
\textit{In particular,}%
\[
K\cap Kt\neq \emptyset \qquad (t\in B_{\delta }^{\Delta }),
\]%
\textit{or, equivalently,}%
\begin{equation}
B_{\delta }^{\Delta }\subseteq K^{-1}K,  \tag{$\ast \ast $}
\end{equation}%
\textit{so that }$B_{\delta }^{\Delta }$ \textit{has compact closure.}

\textit{A fortiori, if }$\mu _{-}(K)>0$, \textit{then }$\delta ,\Delta >0$
\textit{may be chosen with }$\Delta <\mu _{-}(K)$ \textit{and} $B_{\delta
}\subseteq B_{\delta }^{\Delta }$ \textit{so that (*) and (**) hold with }$%
B_{\delta }$ \textit{replacing }$B_{\delta }^{\Delta },$\textit{\ and in
particular }$G$ \textit{is locally compact.}

\bigskip

\noindent \textbf{Proof.} For the first part take $\Delta :=\mu _{-}^{%
\mathbf{t}}(K)/4.$ Then, for $t\in B_{\delta }^{\Delta }$ with $\delta >0$
arbitrary, $\mu (Kt)>2\Delta$ $,$ and so, since $t_{n}\in B_{\delta }$ for
all large enough $n,$ also $t_{n}\in B_{\delta }^{\Delta }$ for all large
enough $n,$ so $B_{\delta }^{\Delta }(K)\backslash \{1_{G}\}$ is non-empty
for $\mathbf{t}$ non-trivial.

Put $H_{t}:=K\cap Kt\subseteq K$. By outer regularity of $\mu $, choose $U$
open with $K\subseteq U$ and $\mu (U)<\mu (K)+\Delta .$ By upper
semicontinuity of $t\mapsto Kt$, w.l.o.g. $KB_{\delta }\subseteq U$ for some
$\delta >0.$ For $t\in B_{\delta }^{\Delta },$ by finite additivity of $\mu
, $ since $2\Delta <\mu (Kt)$
\begin{eqnarray*}
2\Delta +\mu (K)-\mu (H_{t}) &\leq &\mu (Kt)+\mu (K)-\mu (H_{t})=\mu (Kt\cup
K) \\
&\leq &\mu (U)\leq \mu (K)+\Delta .
\end{eqnarray*}%
Comparing extreme ends of this chain of inequalities gives%
\[
0<\Delta \leq \mu (H_{t})\qquad (t\in B_{\delta }^{\Delta }).
\]

For $t\in B_{\delta }^{\Delta },$ as $K\cap Kt\in \mathcal{M}_{+}(\mu )$,
take $s\in K\cap Kt\neq \emptyset ;$ then $s=kt$ for some $k\in K,$ so $%
t=k^{-1}s\in K^{-1}K.$ Conversely, $t\in B_{\delta }^{\Delta }\subseteq
K^{-1}K$ yields $t=k^{-1}k^{\prime }$ for some $k,k^{\prime }\in K;$ then $%
k^{\prime }=kt\in K\cap Kt$.

By the compactness of $K^{-1}K,$ $B_{\delta }^{\Delta }$ has compact closure.

As for the final assertions, if $\mu _{-}(K)>0,$ take $\Delta :=\mu
_{-}(K)/2.$ Then $\inf \{\mu (Kt):t\in B_{\delta }\}>\Delta$ for all small
enough $\delta >0,$ and so in particular $\mu (Kt)>\Delta$ for $t\in
B_{\delta },$ i.e. $B_{\delta }\subseteq B_{\delta }^{\Delta }.$ So the
argument above applies for such $\delta >0$ with $B_{\delta }$ in lieu of $%
B_{\delta }^{\Delta }$, just as before. Here the compactness of $K^{-1}K$
now implies local compactness of $G$ itself. $\square $

\bigskip

As an immediate and useful corollary, we have

\bigskip

\noindent \textbf{Lemma 1}$^{\prime }$\textbf{.} \textit{For }$\mu \in
\mathcal{P}_{\text{fin}}(G),$\textit{\ with }$G$ \textit{\ a metric group, }%
\textit{\ any null sequence }$\mathbf{t}$ \textit{and any }$K\in \mathcal{K}%
(G):$ \textit{if }$\mu _{-}^{\mathbf{t}}(K)>0$,\textit{\ then there is }$m$
\textit{\ }$\in \mathbb{N}$ \textit{with }%
\begin{equation}
0<\mu _{-}^{\mathbf{t}}(K)/4<\mu (K\cap Kt_{n})\qquad (n>m).
\tag{$\ast
^{\prime }$}
\end{equation}%
\textit{In particular,}%
\begin{equation}
t_{n}\in K^{-1}K\qquad (n>m).  \tag{$\ast \ast ^{\prime }$}
\end{equation}

\noindent \textbf{Proof.} Apply Lemma 1 to obtain $\Delta ,\delta >0;$ for $%
t\in B_{\delta }^{\Delta }$, $\mu (Kt)>\Delta,$ so as above $t_{n}\in
B_{\delta }^{\Delta }$ for all large enough $n.$ $\square $

\bigskip

This permits a connection with left Haar null sets; recall that a group $G$
is \textit{amenable at} $1$ [Sol2] (see below for the origin of this term)
if, given $\mathbf{\mu }:=:\{\mu _{n}\}_{n\in \mathbb{N}}\subseteq \mathcal{P%
}(G)$ with $1_{G}\in $ \textrm{supp}$(\mu _{n}),$ for $n\in \mathbb{N}$
there are $\sigma $ and $\sigma _{n}$ in $\mathcal{P}(G)$ with $\sigma
_{n}\ll \mu _{n}$ satisfying:
\[
\sigma _{n}\ast \sigma (K)\rightarrow \sigma (K)\qquad (K\in \mathcal{K}%
(G)).
\]%
We refer to $\sigma $ (or $\sigma (\mathbf{\mu })$ if context requires) as a
\textit{Solecki measure} and to the measures $\sigma _{n}$, if needed, as
\textit{associated measures }(corresponding to the sequence $\{\mu
_{n}\}_{n\in \mathbb{N}}).$

Solecki explains ([Sol2, end of \S 2]) the use of the term `amenability at
1' as a localization (via the restriction that supports contain $1_{G})$ of
a \textit{Reiter-like condition }[Pat, Prop. 0.4] which characterizes
amenability: for $\mu \in \mathcal{P}(G)$ and $\varepsilon >0,$ there is $%
\nu \in \mathcal{P}(G)$ with%
\[
|\nu \ast \mu (K)-\nu (K)|<\varepsilon \qquad (K\in \mathcal{K}(G)).
\]

\bigskip

Lemma 1 and the next several results disaggregate Solecki's Interior-point
Theorem [Sol2, Th 1(ii)] (Corollary 2 below), shedding more light on it and
in particular connecting it to shift-compactness (Theorem 3 below). Indeed,
we see that interior-point theorem itself as an `aggregation' phenomenon.
Theorem 5 of Part II below identifies subgroups with a `disaggregation'
topology, refining $\mathcal{T}_{G}$ by using sets of the form $B_{\delta
}^{K\Delta }(\sigma )$, the measures $\sigma $ being provided in our first
result:

\bigskip

\noindent \textbf{Theorem 1 }(\textbf{Subcontinuity Theorem}, after Solecki
[Sol2, Th. 1(ii)]).\textit{\ For }$G\ $\textit{Polish and amenable at }$%
1_{G} $ \textit{and }$\mathbf{t}$ \textit{a null sequence, there is }$\sigma
=\sigma (\mathbf{t})\in \mathcal{P}(G)$ \textit{such that for each }$K\in
\mathcal{K}_{+}(\sigma )$\textit{\ there is a subsequence }$\mathbf{s}=%
\mathbf{s}(K):=\{t_{m(n)}\}$ \textit{with}%
\[
\lim\nolimits_{n}\sigma (Kt_{m(n)})=\nu (K)\qquad \text{(}n\in \mathbb{N}%
\text{)},\text{ so}\qquad \sigma _{\_}^{\mathbf{s}}(K)>0.
\]%
\textit{\ }

\noindent \textbf{Proof.} For $\mathbf{t}=\{t_{n}\}$ null, put $\mu
_{n}:=2^{n-1}\sum\nolimits_{m\geq n}2^{-m}\delta _{t_{m}^{-1}}\in \mathcal{P}%
(G);$ then $1_{G}$ $\in $ \textrm{supp}$(\mu _{n})\supseteq
\{t_{m}^{-1}:m>n\}.$ By definition of amenability at $1_{G}$, in $\mathcal{P}%
(G)$ there are $\sigma $ and $\sigma _{n}\ll \mu _{n},$ with $\sigma
_{n}\ast \sigma (K)\rightarrow \sigma (K)$ for all $K\in \mathcal{K}(G).$
For $n\in \mathbb{N}$ choose $\alpha _{mn}\geq 0$ with $\sum\nolimits_{m\geq
n}\alpha _{mn}=1$ $(n\in \mathbb{N})$ and with $\sigma
_{n}:=\sum\nolimits_{m\geq n}\alpha _{mn}\delta _{t_{m}^{-1}}.$

Fix $K\in \mathcal{K}_{+}(\sigma )$ and $\theta $ with $0<\theta <1.$ As $K$
is compact, $\sigma _{n}\ast \sigma (K)\rightarrow \sigma (K);$ then w.l.o.g.%
\[
\sigma _{n}\ast \sigma (K)>\theta \sigma (K)\qquad (n\in \mathbb{N}).
\]%
Then for each $n$
\[
\sup \{\sigma (Kt_{m}):m\geq n\}\cdot \sum\nolimits_{m\geq n}\alpha
_{mn}\geq \sum\nolimits_{m\geq n}\alpha _{mn}\sigma (Kt_{m})>\theta \sigma
(K).
\]%
So for each $n$ there is $m=m(\theta )\geq n$ with%
\[
\sigma (Kt_{m})>\theta \sigma (K).
\]%
Now choose $m(n)\geq n$ inductively so that $\sigma
(Kt_{m(n)})>(1-2^{-n})\sigma (K);$ then, by Proposition 1, $\lim_{n}\sigma
(Kt_{m(n)})=\sigma (K):$ $\sigma $ is subcontinuous along $\mathbf{s}$ $%
:=\{t_{m(n)}\}$ on $K.$ $\square $

\bigskip

\textbf{Remark. }The selection above of the subsequence $\mathbf{s}$ mirrors
the role of `admissible directions' which we encounter later in
Cameron-Martin theory (\S 7.2 and below \S 8.2).

\bigskip

We are now able to deduce Solecki's interior-point theorem in a slightly
stronger form, which asserts that the sets $B_{\delta }^{\Delta }$
reconstruct the open sets of $G$ using the compact subsets of a
`non-negligible set', as follows.

\bigskip

\noindent \textbf{Theorem 2} (\textbf{Aggregation Theorem}). \textit{For }$%
G\ $\textit{Polish and amenable at }$1_{G},$\textit{\ if} $E\in \mathcal{U}%
(G)$\textit{\ is not left Haar null -- then, setting}%
\[
\hat{E}:=\bigcup\nolimits_{\delta ,\Delta >0,g\in G,\mathbf{t}}\{B_{\delta
}^{gK,\Delta }(\nu (\mathbf{t})):K\subseteq E,K\in \mathcal{K}_{+}(\nu (%
\mathbf{t})),\Delta <\nu (\mathbf{t})(gK)\},
\]%
\[
1_{G}\in \mathrm{int}(\hat{E})\subseteq \hat{E}\subseteq E^{-1}E.
\]%
\textit{In particular, for }$E$\textit{\ open, }$1_{G}\in \mathrm{int}(\hat{E%
}).$

\bigskip

\noindent \textbf{Proof.} Suppose otherwise; then for each $n$ there is%
\[
t_{n}\in B_{1/n}\backslash \hat{E}.
\]%
Consider $\sigma =\sigma (\mathbf{t}).$ As $E$ is not left Haar null, there
is $g$ with $\sigma (gE)>0.$ Choose compact $K\subseteq gE$ with $\sigma
(K)>0.$ Then with $h:=g^{-1}$ and $H:=$ $hK\subseteq E,$ $\sigma (K)=\sigma
(gH)=\sigma ^{\mathbf{s}}(gH)>0$ for some subsequence $\mathbf{s}$ $%
=\{t_{m(m)}\}$. So, as in Lemma 1, with $\Delta :=\sigma (gH)/4$ for some $%
\delta >0$%
\[
B_{\delta }^{gH\Delta }(\sigma (\mathbf{t}))\subseteq
(gH)^{-1}gH=H^{-1}H\subseteq E^{-1}E.
\]%
Choose $n$ with $n>1/\delta .$ Then $t_{n}\in B_{\delta }$ for all $m>n$; so
for infinitely many $k$%
\[
t_{m(k)}\in B_{\delta }^{gH\Delta }(\sigma (\mathbf{t}))\subseteq \hat{E},
\]%
a contradiction. As for the last assertion, for $E$ open, $D$ countable and
dense, $G\subseteq \bigcup\nolimits_{d\in D}dE,$ so for any $\mu \in
\mathcal{P}(G)$ (in particular for $\sigma $) $\mu (dE)>0$ for some $d\in D$%
, and so $E$ is not left Haar null. $\square $

\bigskip

The immediate consequence is

\bigskip

\noindent \textbf{Corollary 2} (\textbf{Solecki's Interior-Point Theorem }%
[Sol2, Th 1(ii)]). \textit{For }$G\ $\textit{Polish and amenable at }$1_{G},$%
\textit{\ if} $E\in \mathcal{U}(G)$\textit{\ is not left Haar null, then }$%
1_{G}\in \mathrm{int}(E^{-1}E)$.

\bigskip

\noindent \textbf{Corollary 2}$^{\prime }$\textbf{.} \textit{For }$G\ $%
\textit{a Polish group, }\textit{if }$B\in \mathcal{U}(G)$ \textit{is not
left Haar null and is in }$\mathcal{M}_{+}(\mu )$\textit{\ for some
subcontinuous\ }$\mu \in \mathcal{P}_{\text{fin}}(G)$\textit{, then }$($**$)$
\textit{holds for some }$\delta >0.$

\textit{In particular, }$($**$)$ \textit{holds in a locally compact group }$%
G $\textit{, for any Baire non-meagre set }$E.$

\bigskip

\noindent \textbf{Proof. }The first assertion is immediate from Lemma 1. As
for the second, for a non-meagre Baire set $E$, if $\tilde{E}$ is the
quasi-interior and $K\subseteq \tilde{E}$ is compact with non-empty
interior, then $\eta _{G}(K)>0.$ Since $\eta $ is subcontinuous, there is $%
\delta >0$ with\textit{\ }%
\[
Kt\cap K\neq \emptyset \qquad (||t||<\delta ).
\]%
A fortiori,
\[
\tilde{E}t\cap \tilde{E}\neq \emptyset \qquad (||t||<\delta );
\]%
then $U:=(Et)^{\widetilde{}}\cap \tilde{E}\neq \emptyset ,$ since $(Et)^{%
\widetilde{}}=\tilde{E}t$ (the Nikodym property of the \textit{usual}
topology of $G$). So since $U$ is open and non-meagre, also $Et\cap E\neq
\emptyset $, and so again (**). $\square $

\bigskip

The next result establishes the embeddability by (left-sided) translation of
an appropriate subsequence of a given \textit{null} sequence into a given
target set that (like-sidedly) is non-left-Haar null. This property of
embedding into a \textit{non-negligible} set, first studied in respect of
category and measure negligibility on $\mathbb{R}$ by Kestelman and much
later independently by Borwein and Ditor and thereafter also by other
authors, mostly for combinatorial challenges, has emerged as an important
general unifying principle, termed \textit{shift-compactness}, applicable in
a much wider context embracing metric groups $G$ under various topologies
refining $\mathcal{T}_{G}$ and so defining various notions of negligibility:
for the background here see [BinO2,3], [MilO]. Its consequences include
various uniform-boundedness theorems as well as the Effros and the Open
Mapping Theorems. Here we establish the said property, announced in [MilO],
in relation to the ideal $\mathcal{HN}$ of left Haar null sets. (It is a $%
\sigma $-ideal for Polish $G$ in the presence of amenability at 1 [Sol2, Th
1(i)].) This leaves open the `converse question' of a refinement topology
for which $\mathcal{HN}$ is the associated notion of negligibility; this
seems plausible under the continuum hypothesis, CH, if one restricts
attention only to Borel sets in $\mathcal{HN}$ and their subsets by lifting
a result concerning $\mathbb{R}$ in [CieJ, Cor. 4.2] to $G$ -- see also the
Remark following our next result.

\bigskip

\noindent \textbf{Theorem 3 }(\textbf{Shift-compactness Theorem} for $%
\mathcal{HN}$).\textit{\ For }$G\ $\textit{Polish and amenable at }$1_{G},$%
\textit{\ if} $E\in \mathcal{U}(G)$\textit{\ is not left Haar null and }$%
z_{n}$ \textit{is null, then there are }$s\in E$ \textit{and an infinite }$%
\mathbb{M}\subseteq \mathbb{N}$ \textit{with }%
\[
\{sz_{m}:m\in \mathbb{M}\}\subseteq E.
\]%
\textit{Indeed, this holds for quasi all }$s\in E,$ \textit{i.e. off a left
Haar null set.}

\bigskip

\noindent \textbf{Proof.} Put $t_{n}:=z_{n}^{-1},$ which is null. With $%
\sigma =\sigma (\mathbf{t})$ as in the Subcontinuity Theorem, since $E$ is
not left Haar null, there is $g$ with $\sigma (gE)>0.$ For this $g,$ put $%
\mu :=\!_{g}\sigma .$ Fix a compact $K_{0}\subseteq E$ with $\mu (K_{0})>0$
and then, passing to a subsequence of $\mathbf{t}$ as necessary (by Th. 1),
we may assume that $\mu _{-}^{\mathbf{t}}(K_{0})>0.$ Proceed to choose
inductively a sequence $m(n)\in \mathbb{N}$, and decreasing compact sets $%
K_{n}\subseteq K_{0}\subseteq E$ with $\mu (K_{n})>0$ such that%
\[
\mu (K_{n}\cap K_{n}t_{m(n)})>0.
\]%
To check the inductive step, suppose $K_{n}$ already defined. As $\mu
(K_{n})>0,$ by the Subcontinuity Theorem, there is a subsequence $\mathbf{%
s=s(}K_{n})$ of $\mathbf{t}$ with $\mu _{-}^{\mathbf{s}}(K_{n})>0.$ By Lemma
1$^{\prime },$ there is $k(n)>n$ such that $\mu (K_{n}\cap K_{n}s_{k(n)})>0.$
Putting $t_{m(n)}=s_{k(n)}$ and $K_{n+1}:=K_{n}\cap K_{n}t_{m(n)}\subseteq
K_{n}$ completes the inductive step.

By compactness, select $s$ with%
\[
s\in \bigcap\nolimits_{m\in \mathbb{N}}K_{m}\subseteq K_{n+1}=K_{n}\cap
K_{n}t_{m(n)}\qquad (n\in \mathbb{N});
\]%
choosing $k_{n}\in K_{n}\subseteq K$ with $s=k_{n}t_{m(n)}$ gives $s\in
K_{0}\subseteq E,$ and%
\[
sz_{m(n)}=st_{m(n)}^{-1}=k_{n}\in K_{n}\subseteq K_{0}\subseteq E.
\]%
Finally take $\mathbb{M}:=\{m(n):n\in \mathbb{N}\}.$

As for the final assertion, we follow the idea of the Generic Completeness
Principle [BinO1, Th. 3.4] (but with $\mathcal{U}(G)$ for $\mathcal{B}a$
there):\ define
\[
F(H):=\bigcap\nolimits_{n\in \mathbb{N}}\bigcup\nolimits_{m>n}H\cap
Ht_{m}\qquad (H\in \mathcal{U}(G));
\]%
then $F:\mathcal{U}(G)\rightarrow \mathcal{U}(G)$ and $F$ is monotone $%
(F(S)\subseteq F(T)\ $for $S\subseteq T);$ moreover, $s\in F(H)$ iff $s\in H$
and $sz_{m}\in H$ for infinitely many $m$. We are to show that $%
E_{0}:=E\backslash F(E)$ is left Haar null. Suppose otherwise. Then renaming
$g$ and $K_{0}$ as necessary, w.l.o.g. both $\mu (E_{0})>0$ and $%
K_{0}\subseteq E_{0}$ (and $\mu (K_{0})>0).$ But then, as above, $\emptyset
\neq F(K_{0})\cap K_{0}\subseteq F(E)\cap E_{0},$ a contradiction, since $%
F(E)\cap E_{0}=\emptyset .$ $\square $

\bigskip

\noindent \textbf{Remark. }In the setting of Th. 3 any non-empty open set $U$
is not left Haar null (as $\{dU:D\in D\}$ with $D\ $countable dense covers $%
G $), hence neither is $U\backslash H$ for $H\in \mathcal{HN}$. So the
(Hashimoto ideal) topology generated by such sets includes $\mathcal{HN}$
among its negligible sets.

\bigskip

\noindent \textbf{Corollary 3. }\textit{For }$G\ $\textit{Polish and
amenable at }$1_{G}$ \textit{and }$z_{n}$ \textit{null, there is }$\mu \in
\mathcal{P}(G)$ \textit{such that for }$K\in \mathcal{K}_{+}(\mu )$%
\[
K\cap Kz_{m}^{-1}\in \mathcal{M}_{+}(\mu )\text{ \qquad for infinitely many }%
m\in \mathbb{N},
\]%
\textit{iff for }$\mu $\textit{-quasi all }$s\in K$\textit{\ there is an
infinite }$\mathbb{M}\subseteq \mathbb{N}$ \textit{with }%
\[
\{sz_{m}:m\in \mathbb{M}\}\subseteq K.
\]

\noindent \textbf{Proof.} We will refer to the function $F$ of the preceding
proof. First proceed as in the proof of Th. 3 above, taking $%
t_{n}:=z_{n}^{-1}$ and $g=1_{G}$ (so that $\mu =\sigma ).$ Fix $K$ with $\mu
(K)>0.$ For the forward direction, continue as in the proof of Th. 3 with $%
K_{0}=K$ and observe that the proof above needs only that $s_{k(n)}\in
K_{n}^{-1}K_{n}$ occurs infinitely often whenever $\mu (K_{n})>0$. This
yields the desired conclusion that $\mu (K\backslash F(K))=0.$ For the
converse direction, suppose that $\mu (F(K))>0.$ Since for each $n\in
\mathbb{N}$%
\[
F(K)\subseteq \bigcup\nolimits_{m>n}K\cap Kt_{m},
\]%
we have $\mu (K\cap Kt_{m})>0$ for some $m>n;$ so
\[
K\cap Kt_{m}\in \mathcal{M}_{+}(\mu )\text{ \qquad for infinitely many }%
m.\qquad \square
\]

\bigskip

\noindent \textbf{Remark.} With $E$ as in the Shift-compactness Theorem, if $%
z_{n}\in B_{1/n}\backslash E^{-1}E,$ then $z_{n}$ is null; so, for some $%
s\in E,$ $sz_{m}\in E$ for infinitely many $m.$ Then, for any such $m,$%
\[
z_{m}\in E^{-1}E,
\]%
contradicting the choice of $z_{m}.$ So $1_{G}\in \mathrm{int}(E^{-1}E),$
i.e. $E$ has the Steinhaus-Weil property, as before.

\bigskip

The following sharpens a result due (for Lebesgue measure on $\mathbb{R}$)
to Mospan [Mos] by providing the converse below; it is antithetical to Lemma
1 (and so to Theorem 3).

\bigskip

\noindent \textbf{Proposition 6 (Mospan property). }\textit{For }$G\ $%
\textit{\ a metric group and compact }$K$\textit{, }\textit{if }$1_{G}\notin
\mathrm{int}(K^{-1}K),$\textit{\ then }$\mu _{-}(K)=0,$\textit{\ i.e. }$\mu $
\textit{is maximally discontinuous; equivalently, there is a `null sequence'
}$t_{n}\rightarrow 1_{G}$\textit{\ with }$\lim_{n}\mu (Kt_{n})=0.$\textit{\ }

\textit{Conversely, if }$\mu (K)>\mu _{-}(K)=0,$ \textit{then there is a
null sequence }$t_{n}\rightarrow 1_{G}$\textit{\ with }$\lim_{n}\mu
(Kt_{n})=0,$\textit{\ and there is a compact }$C\subseteq K$\textit{\ with }$%
\mu (K\backslash C)=0$\textit{\ with }$1_{G}\notin \mathrm{int}(C^{-1}C).$

\bigskip

\noindent \textbf{Proof.} The first assertion follows from Lemma 1. For the
converse, as in [Mos]: suppose that $\mu (Kt_{n})=0,$ for some sequence $%
t_{n}\rightarrow 1_{G}.$ By passing to a subsequence, we may assume that $%
\mu (Kt_{n})<2^{-n-1}.$ Put $D_{m}:=K\backslash \bigcap\nolimits_{n\geq
m}Kt_{n}\subseteq K;$ then $\mu (K\backslash D_{m})\leq \sum_{n\geq m}\mu
(Kt_{n})<2^{-m},$ so $\mu (D_{m})>0$ provided $2^{-m}<\mu (K).$ Now choose
compact $C_{m}\subseteq D_{m},$ with $\mu (D_{m}\backslash C_{m})<2^{-m}.$
So $\mu (K\backslash C_{m})<2^{1-m}.$ Also $C_{m}\cap C_{m}t_{n}=\emptyset ,$
for each $n\geq m,$ as $C_{m}\subseteq K;$ but $t_{n}\rightarrow 1_{G},$ so
the compact set $C_{m}^{-1}C_{m}$ contains no interior points. Hence, by
Baire's theorem, neither does $C^{-1}C,$ since $C=\bigcup\nolimits_{m}C_{m},$
which differs from $K$ by a null set. $\square $

\bigskip

\noindent \textbf{Proposition 7.} \textit{A (regular) Borel measure }$\mu $%
\textit{\ on a locally compact metric topological group }$G$ \textit{has the
Steinhaus-Weil property iff}\newline
\noindent (i) \textit{for each non-null compact set }$K,$\textit{\ the map }$%
m_{K}:t\rightarrow \mu (Kt)$\textit{\ is subcontinuous at }$1_{G};$\newline
\noindent (ii) \textit{for each non-null compact set }$K$\textit{\ there is
no `null' sequence }$t_{n}\rightarrow 1_{G}$\textit{\ with }$\mu
(Kt_{n})\rightarrow 0.$

\bigskip

\noindent \textbf{Remark. }This is immediate from Prop. 6 (cf. [Mos]).

\bigskip

We now prove a strengthening of the Subcontinuity Theorem obtained by
assuming a `concentration property'. That this property holds in a abelian
Polish group emerges from an inspection of Solecki's proof of his theorem
that an abelian Polish group is amenable at 1.

\bigskip

\noindent \textbf{Definitions.} Say that a null sequence $\mathbf{t}$ is
\textit{regular} if $\mathbf{t}$ is non-trivial, $||t_{k}||$ is
non-increasing, and
\[
||t_{k}||\leq r(k):=1\left/ [2^{k}(k+1)]\right. \qquad (k\in \mathbb{N}).
\]

For regular $\mathbf{t}$, put $\mu _{k}=\mu _{k}(\mathbf{t}%
):=2^{k-1}\sum\nolimits_{m\geq k}2^{-m}(\delta _{t_{m}^{-1}}+\delta
_{t_{m}})=\frac{1}{2}\delta _{t_{k}^{-1}}+\frac{1}{4}\delta
_{t_{k+1}^{-1}}+....$ Then $\mu _{k}(B_{r(k)})=1$ for $k\in \mathbb{N}.$
Merging $\mathbf{t}^{-1}$ with $\mathbf{t}$ by alternation of terms if
necessary, it is now convenient to assume that $\mathbf{t}$ contains as
successive pairs inverses of its terms. So, if $\nu _{k}\ll \mu _{k}$, then
\[
\nu _{k}:=\sum\nolimits_{m\geq k}a_{km}\delta _{t_{m}^{-1}},
\]%
for some non-negative sequence $\mathbf{a}_{k}:=%
\{a_{kk,}a_{k,k+1,}a_{k,k+2},...\}$ of unit $\ell _{1}$-norm. Say that $\{%
\mathbf{a}_{k}\}$ has the \textit{concentration property }if for some index $%
j$ and some $\alpha >0$
\[
a_{k,k+j}\geq \alpha >0\qquad \text{for all large }k;
\]%
then say that the sequence $\{\nu _{k}\}$ has the \textit{concentration
property}. (This will fail if $\mathbf{a}_{k}$ has $a_{k,k+k}=1,$ which
concentrates measure in an unbounded fashion.)

\bigskip

\noindent \textbf{Definition.} Say that a group $G$ is \textit{strongly
amenable at 1} if $G$ is amenable at 1, and for each \textit{regular} $%
\mathbf{t}$ a Solecki measure $\sigma (\mathbf{t})$ has associated measures $%
\sigma _{k}(\mathbf{t})\ll \mu _{k}(\mathbf{t})$ with the concentration
property.

\bigskip

\noindent \textbf{Theorem 4 }(\textbf{Strong amenability at 1}, after [Sol2,
Prop. 3.3(i)]). \textit{Any abelian Polish group }$G$\textit{\ is strongly
amenable at 1.}

\bigskip

\noindent \textbf{Proof. }This follows the construction in [Sol2] of the
reference measure in the case of $\mu _{k}(\mathbf{t})$ above. First define
the normalized restriction

\[
\sigma _{k}:=\mu _{k}|B_{r(k)}\left/ \mu _{k}(B_{r(k)})\right.
\]%
and then set%
\[
\sigma :=\ast _{k=1}^{\infty }\rho _{k}\text{ for }\rho _{k}:=\frac{1}{k+1}%
\sum\nolimits_{i=0}^{k}\sigma _{k}^{i}.
\]%
(Convolution powers intended here.) Then the argument in [Sol2] shows that $%
\sigma _{k}\ast \sigma (K)\rightarrow \sigma (K)$ (for $K$ compact).
However, as $\mathbf{t}$ is regular, $\mu _{k}\equiv \sigma _{k}$. But $%
a_{kk}=1/2$ (all $k),$ so the measures $\sigma _{k}$ here have the
concentration property. $\square $

\bigskip

\noindent \textbf{Definitions (Sequence and measure symmetrization):}

\noindent 1. Merging $\mathbf{t}^{-1}$ with $\mathbf{t}$ by alternation of
terms yields the regular sequence $\mathbf{s}%
=(s_{1},s_{2},...):=(t_{1},t_{1}^{-1},t_{2},t_{2}^{-1},...)$; we term this
the \textbf{symmetrized sequence} of $\mathbf{t}$. (It is `symmetric' in the
sense only that $||s_{2k-1}||=||s_{2k}||$.)

\noindent 2. For odd $k,$ as $(\mu _{k}(\mathbf{t})+\mu _{k}(\mathbf{t}%
^{-1}))/2$ is symmetric as a measure, taking $\sigma _{2k}(\mathbf{s}%
)=\sigma _{2k-1}(\mathbf{s}):=(\mu _{k}(\mathbf{t})+\mu _{k}(\mathbf{t}%
^{-1}))/2$ in lieu of $\mu _{k}(\mathbf{t})$ above yields each $\rho _{k}$
symmetric. So, in the abelian context of Theorem 4 above, the limiting
convolution $\sigma $ is a \textit{symmetric Solecki measure }$\sigma (%
\mathbf{t}).$

\bigskip

\noindent \textbf{Remark. }Performing the symmetrization of the Definition
above gives in the proof of Theorem 4 above that $%
a_{2k-1,2k-1}=a_{2k-1,2k}=a_{2k,2k}=a_{2k,2k+1}=1/4$, which presents \textit{%
simultaneous concentration} along $\mathbf{t}$ and $\mathbf{t}^{-1}.$

\bigskip

We now re-run the proof of Theorem 1 with improved estimates to yield:

\bigskip

\noindent \textbf{Theorem 1}$_{\text{\textbf{S}}}$ \textbf{(Strong
Subcontinuity Theorem). }\textit{For }$G$\textit{\ a Polish group that is
strongly amenable at 1, if }$\mathbf{t}$\textit{\ is regular and }$\sigma
=\sigma (\mathbf{t})$\textit{\ is a Solecki measure -- then for }$K\in
\mathcal{K}_{+}(\sigma )$%
\[
\sigma (K)=\lim_{n}\sigma (Kt_{n})=\sigma _{-}^{\mathbf{t}}(K).
\]%
\textit{Likewise, passing to the symmetrized sequence of }$\mathbf{t}$%
\textit{\ as above and to a symmetric Solecki measure }$\sigma (\mathbf{t})$
\textit{with the simultaneous concentration property (for }$\mathbf{t}$
\textit{and }$\mathbf{t}^{-1}$\textit{) corresponding to an abelian context:}%
\[
\sigma (K)=\lim_{n}\sigma (t_{n}K).
\]

\bigskip

\noindent \textbf{Proof. }Fix $\mathbf{t}$\textit{\ }and a corresponding
Solecki reference measure $\nu (\mathbf{t})$\textit{\ }and its associated
sequence $\nu _{k}$, which as in Th. 4 has the concentration property. Write
$\sigma _{k}:=\sum\nolimits_{m\geq k}a_{km}\delta _{t_{m}^{-1}}$; as $\sigma
_{k}$ has the concentration property, there are $n_{0},$ $j$ and $\alpha >0$
with
\[
a_{kk+j}\geq \alpha >0\qquad (k\geq n_{0}).
\]%
Now fix $K\ $compact with $\sigma (K)>0$ and $\varepsilon >0.$ Put%
\[
\delta :=\varepsilon \left/ \left( \frac{2}{\alpha }-1\right) \right. >0,
\]%
as $\alpha \leq 1.$ Then, by upper semicontinuity and by `amenability at 1'
(i.e. $\sigma _{k}\ast \sigma (K)\rightarrow \sigma (K)$), there is $%
n_{1}=n_{1}(\varepsilon ,K)>n_{0}$ with%
\[
\sigma (Kt_{k})\leq \sigma (K)+\delta \text{ and }\sigma _{k}\ast \sigma
(K)\geq \sigma (K)-\delta \qquad (k\geq n_{1}).
\]%
So (by upper semicontinuity) for $k\geq n_{1}$%
\[
\sum\nolimits_{m\geq k,m\neq k+j}a_{km}\sigma (Kt_{m})\leq
\sum\nolimits_{m\geq k,m\neq k+j}a_{km}(\sigma (K)+\delta )=(\sigma
(K)+\delta )(1-a_{kk+j}).
\]%
Also (by `amenability at 1') for $k\geq n_{1}$%
\begin{eqnarray*}
a_{kk+j}\sigma (Kt_{k+j}) &\geq &\sigma (K)-\delta -\sum\nolimits_{m\geq
k,m\neq k+j}a_{km}\sigma (Kt_{m}) \\
&\geq &\sigma (K)+\delta -2\delta -(\sigma (K)+\delta )(1-a_{kk+j}) \\
&=&a_{kk+j}\sigma (K)-\delta (2-a_{kk+j}).
\end{eqnarray*}%
So for $m=k+j>n_{1}+j$%
\[
\sigma (Kt_{m})=\sigma (Kt_{k+j})\geq \sigma (K)-\delta \left( \frac{2}{%
a_{kk+j}}-1\right) \geq \sigma (K)-\delta \left( \frac{2}{\alpha }-1\right)
=\sigma (K)-\varepsilon .\qquad
\]%
As for the final assertion concerning symmetrization, note that $\sigma
(t_{n}K)=\sigma (K^{-1}t_{n}^{-1})\rightarrow \sigma (K^{-1})=\sigma (K),$
by symmetry of $\sigma .\square $

\bigskip

We note an immediate corollary, needed in \S 4.

\bigskip

\noindent \textbf{Corollary 4}. \textit{For }$G,\mathbf{t}$\textit{\ and }$%
\sigma $ \textit{as in Th. 1}$_{\text{S}}$ \textit{above, and }$K,H\in
\mathcal{K}_{+}(\sigma ),\delta >0$\textit{: if }$0<\Delta <\sigma (K)$%
\textit{\ and }$0<D<\sigma (H),$\textit{\ then there is }$n$\textit{\ with}
\[
B_{\delta }^{K\Delta }\cap B_{\delta }^{HD}\supseteq \{t_{m}:m\geq n\}.
\]

\bigskip

\noindent \textbf{Proof.} Take $\varepsilon :=\min \{\sigma (K)-\Delta
,\sigma (H)-D\}>0.$ As $K,H\in \mathcal{K}_{+}(\sigma ),$ there is $n$ such
that $||t_{m}||<\delta $ for $m\geq n$ and
\[
\sigma (Kt_{m})\geq \sigma (K)-\varepsilon \geq \Delta ,\qquad \sigma
(Ht_{m})\geq \sigma (H)-\varepsilon \geq D\qquad (m\geq n).\qquad \square
\]

\bigskip

To accommodate varying sided-ness conventions, we close with the left-handed
version of Theorem 1, in which $\tilde{\nu}_{\_}^{\mathbf{s}}(K)$ (defined
at the start of the section) replaces $\nu _{-}^{\mathbf{s}}(K).$

\bigskip

\noindent \textbf{Theorem 1}$_{\text{\textbf{L}}}$\textbf{\ }(\textbf{Left
Subcontinuity Theorem}, after Solecki [Sol2, Th. 1(ii)]).\textit{\ For }$G\ $%
\textit{amenable at }$1_{G}$\textit{, and }$\mathbf{t}$ \textit{a null
sequence, there is }$\sigma =\sigma _{L}(\mathbf{t})\in \mathcal{P}(G)$
\textit{such that for each }$K\in \mathcal{K}_{+}(\sigma )$\textit{\ there
is a subsequence }$\mathbf{s}=\mathbf{s}(K):=\{t_{m(n)}\}$ \textit{with}%
\[
\lim\nolimits_{n}\sigma (t_{m(n)}K)=\sigma (K)\qquad \text{(}n\in \mathbb{N}%
\text{)},\text{ so}\qquad \tilde{\sigma}_{\_}^{\mathbf{s}}(K)>0.
\]

\noindent \textbf{Proof.} From Theorem 1 with $\sigma =\sigma (\mathbf{t}),$
recalling that $\tilde{\sigma}(E):=\sigma (E^{-1})$ so that $\sigma (E)=%
\tilde{\sigma}(E^{-1}),$ there is $\{t_{m(n)}\}$ with
\[
\tilde{\sigma}(t_{m(n)}^{-1}K^{-1})>(1-2^{-n})\tilde{\sigma}(K^{-1})\qquad
\text{(}n\in \mathbb{N}\text{)},\text{ so}\qquad \lim_{n}\tilde{\sigma}%
(t_{m(n)}^{-1}K^{-1})>0.
\]%
Since $K^{-1}$ is compact and $t_{m(n)}^{-1}$ is null, replace $K^{-1}$ by $%
K,$ $t_{m(n)}^{-1}$ by $t_{m(n)}$; then with $\mu $ for $\tilde{\sigma}$
\[
\mu (t_{m(n)}K)>(1-2^{-n})\mu (K)\qquad \text{(}n\in \mathbb{N}\text{)},%
\text{ so}\qquad \lim_{n}\mu (t_{m(n)}K)>0.\qquad \square
\]

\part{Around the Simmons-Mospan Theorem}

\begin{center}
\bigskip
\end{center}

We saw in Part I that measure may exhibit selective continuity under
translation in the presence of amenability at 1. Here we consider
obstacles/obstructions to continuity.

\bigskip

\section{A Lebesgue Decomposition}

We begin with definitions isolating left-handed components in Christensen's
notion of Haar null sets [Chr1], and Solecki's left Haar null sets [Sol2];
right-handed versions have analogous properties. As far as we are aware the
component notions in parts (ii)-(iv) below have not been previously studied.
Below $G$ is a Polish group.

\bigskip

\noindent \textbf{Definition. (}i)\textbf{\ Left }$\mu $\textbf{-null: }For $%
\mu \in \mathcal{M}(G)$, say that $N$ is \textit{left }$\mu $\textit{-null }(%
$N\in \mathcal{M}_{0}^{\text{L}}(\mu )$) if it is contained in a universally
measurable set $B$ such that%
\[
\mu (gB)=0\qquad (g\in G).
\]%
Thus a set $S$ is \textit{left Haar null} ([Sol3] after [Chr1]) if it is
contained in a universally measurable set $B$ that is left $\mu $-null%
\textit{\ }for some $\mu \in \mathcal{M}(G).$

\noindent (ii) \textbf{Left }$\mu $\textbf{-inversion:} For $\mu \in
\mathcal{M}(G),$ say that $N\in \mathcal{M}_{0}^{\text{L}}(\mu )$ is \textit{%
left invertibly }$\mu $\textit{-null }($N\in \mathcal{M}_{0}^{\text{L-inv}%
}(\mu )$) if
\[
N^{-1}\in \mathcal{M}_{0}^{\text{L}}(\mu ),
\]%
so that $N^{-1}$ is contained in a universally measurable set $B^{-1}$ such
that%
\[
\mu (gB^{-1})=0\qquad (g\in G).
\]%
\noindent (iii) \textbf{Left }$\mu $\textbf{-absolute continuity: }For $\mu
,\nu \in \mathcal{M}(G),$ $\nu $ is \textit{left absolutely continuous}
w.r.t. $\mu $ ($\nu <^{\text{L}}\mu $) if $\nu (N)=0$ for each $N\in
\mathcal{M}_{0}^{\text{L}}(\mu ),$ and likewise for the invertibility
version: $\nu <^{\text{L-inv}}\mu .$

\noindent (iv) \textbf{Left }$\mu $\textbf{-singularity:} For $\mu ,\nu \in
\mathcal{M}(G),\nu $ is \textit{left singular} w.r.t. $\mu $ (on $B$) ($\nu
\bot ^{\text{L}}\mu $ (on $B)$) if $B$ is a support of $\nu $ and $B$ is in $%
\mathcal{M}_{0}^{\text{L}}(\mu ),$ and likewise $\nu \bot ^{\text{L-inv}}\mu
.$

\bigskip

\noindent \textbf{Remark. }For $\mu $ symmetric, since%
\[
_{g^{-1}}\mu (B)=\mu _{g}(B^{-1})
\]%
if $B$ is left $\mu $-null we may conclude only that $B^{-1}$ is \textit{%
right} $\mu $-null. The `inversion property', property (ii) above, is thus
quite strong (though obvious in the abelian case).

\bigskip

Notice that each of $\mathcal{M}_{0}^{\text{L}}(\mu )$ and $\mathcal{M}_{0}^{%
\text{L-inv}}(\mu )$ forms a $\sigma $-algebra (since $g\bigcup%
\nolimits_{n}B=\bigcup\nolimits_{n}gB$ and $g\left(
\bigcup\nolimits_{n}B\right) ^{-1}=\bigcup\nolimits_{n}gB^{-1}$). This
implies the following left versions of the Lebesgue Decomposition Theorem
(we need the second one below). The `pedestrian' proof below demonstrates
that the Principle of Dependent Choice (DC) suffices, a further example that
`positive' results in measure theory follow from DC (as Solovay points out
in [Solo, p. 31]).

\bigskip

\noindent \textbf{Theorem LD.} \textit{For }$G$\textit{\ a Polish group, }$%
\mu ,\nu \in \mathcal{M}(G),$ \textit{there are }$\nu _{\text{a}},\nu _{%
\text{s}}$ $\in \mathcal{M}(G)$ \textit{with}%
\[
\nu =\nu _{\text{a}}+\nu _{\text{s}}\text{ with }\nu _{\text{a}}<^{\text{L}%
}\mu \text{ and }\nu _{\text{s}}\bot ^{\text{L}}\mu ,
\]%
\textit{and likewise, there are }$\nu _{\text{a}}^{\prime },\nu _{\text{s}%
}^{\prime }$ $\in \mathcal{M}(G)$\textit{\ with}%
\[
\nu =\nu _{\text{a}}^{\prime }+\nu _{\text{s}}^{\prime }\text{ with }\nu _{%
\text{a}}<^{\text{L-inv}}\mu \text{ and }\nu _{\text{s}}\bot ^{\text{L-inv}%
}\mu .
\]

\bigskip

\noindent \textbf{Proof. }As the proof depends on $\sigma $-additivity, it
will suffice to check the `L' case. Write $G=\bigcup\nolimits_{n}G_{n}$ with
the $G_{n}$ disjoint, universally measurable, and with each $\nu (G_{n})$
finite (say, with all but one term $\sigma $-compact, and their complement $%
\nu $-null). Put $s_{n}=\sup \{\nu (E):E\subseteq G_{n},E\in \mathcal{M}%
_{0}^{\text{L}}(\mu )\}.$ In $\mathcal{M}_{0}^{\text{L}}(\mu ),$ for each $n$
with $s_{n}>0,$ choose $E_{n,m}\subseteq G_{n}$ with $\nu
(E_{n,m})>s_{n}-1/m,$ and put $B_{n}:=\bigcup\nolimits_{m}E_{n,m}\subseteq
G_{n}.$ Then the sets $B_{n}$ are disjoint and lie in $\mathcal{M}_{0}^{%
\text{L}}(\mu ),$ as does also $B:=\bigcup\nolimits_{n}B_{n}$; moreover $\nu
(G_{n}\backslash B_{n})=0$ for each $n.$ Put $A:=G\backslash B.$ Then $\nu
(M)=0$ for $M\in \mathcal{M}_{0}^{\text{L}}(\mu )$ with $M\subseteq A,$
since $A=\bigcup\nolimits_{n}(G_{n}\backslash B_{n}).$ So $\nu _{\text{a}%
}:=\nu |A<^{\text{L}}\mu ,$ and $\nu _{\text{s}}:=\nu |B\bot ^{\text{L}}\mu
, $ since $B\in \mathcal{M}_{0}^{\text{L}}(\mu )$. $\square $

\bigskip

\noindent \textbf{Remark. } The above rests on DC; a simpler argument rests
on maximality: choose a maximal disjoint family $\mathcal{B}$ of universally
measurable sets $M\in \mathcal{M}_{0}^{\text{L}}(\mu )$ with finite positive
$\nu (M);$ then, their union $B\in \mathcal{M}_{0}^{\text{L}}(\mu )$ (as $%
\mathcal{B}$ would be countable, by the $\sigma $-finiteness of $\nu ).$

\section{Discontinuity: the Simmons-Mospan Theorem}

It is convenient to begin by repeating the gist of the Simmons-Mospan
argument here, as it is short, despite its `near perfect disguise', to
paraphrase a phrase from Loomis [Loo, p. 85]. The result follows from their
use of Fubini's Theorem and the Lebesgue decomposition theorem of \S 3, but
here we stress the dependence on Theorem FN and on left $\mu $-inversion. We
revert to the Weil left-sided convention and associated $KK^{-1}$ usage.

\bigskip

\noindent \textbf{Proposition 8 (Local almost nullity).} \textit{For }$G$
\textit{a Polish group, }$\mu \in \mathcal{M}(G)$\textit{, }$V$ \textit{open
and }$K$ \textit{compact with }$K\in \mathcal{M}_{0}^{\text{L-inv}}(\mu )$%
\textit{, i.e. }$K,K^{-1}\in \mathcal{M}_{0}^{\text{L}}(\mu ):$\newline
-- \textit{for any }$\nu \in \mathcal{M}(G),$\textit{\ }$\nu (tK)=0$\textit{%
\ for }$\mu $\textit{-almost all }$t\in V,$ \textit{and likewise }$\nu
(Kt)=0.$

\bigskip

\noindent \textbf{Proof. }For\ $\nu $ invertibly $\mu $-absolutely
continuous (as in \S 3 above), the conclusion is immediate; for general $\nu
$ this will follow from Theorem LD (\S 3), once we have proved the
corresponding singular version of the assertion: that is the nub of the
proof.

Thus, suppose that $\nu \bot ^{\text{L-inv}}\mu $\textit{\ }on\textit{\ }$K.$
For $t\in V$ let $t=$ $uw$ be any expression for $t$ as a group product of $%
u,v\in G$, and note that $\mu (uK^{-1})=0,$ as $K^{-1}\in \mathcal{M}_{0}^{%
\text{L}}(\mu )$. Let $H\ $be the set
\[
\bigcup\nolimits_{t\in V}(\{t\}\times tK),
\]%
here viewed as a union of vertical $t$-sections. We next express it as a
union of $u$-horizontal sections and apply the Fubini Null Theorem (Th. FN,
\S 1).

Since $u=tk=uwk$ is equivalent to $w=k^{-1},$ the $u$-horizontal sections of
$H$ may now be rewritten, eliminating $t$, as
\[
\{(t,u):uw=t\in V,u\in tK=uwK\}=\{(uw,u):uw\in V,uw\in uK^{-1}\}.
\]%
So $H$ may now be viewed as a union of $u$-horizontal sections as
\[
\bigcup\nolimits_{u\in G}(V\cap (uK^{-1}))\times \{u\}),
\]%
all of these $u$-horizontal sections being $\mu $-null. By Th. FN, $\mu $%
-almost all vertical $t$-sections of $H$ for $t\in V$ are $\nu $-null. As
the assumptions on $K$ are symmetric the right-sided version follows. $%
\square $

\bigskip

\noindent The result here brings to mind the Dodos Dichotomy Theorem [Dod1,
Th. A] for \textit{abelian} Polish groups $G$: if an analytic set $A$ is
witnessed as Haar-null under one measure $\mu \in \mathcal{P}(G),$ then
either $A$ is Haar-null for quasi all $\nu \in \mathcal{P}(G)$ or else it is
not Haar-null for quasi all such $\nu ,$ i.e. if $A\in \mathcal{M}_{0}(\mu )$
(omitting the unnecessary superscript L), then either $A\in \mathcal{M}%
_{0}(\nu )$ for quasi all such $\nu ,$ or $A\notin \mathcal{M}_{0}(\nu )$
for quasi all such $\nu $ (w.r.t. the Prokhorov-L\'{e}vy metric in $\mathcal{%
P}(G)$ [Dud, 11.3, cf. 9.2]). Indeed, [Dod2, Prop. 5] when $A$ is $\sigma $%
-compact $A$ is Haar-null for quasi all $\nu \in \mathcal{P}(G).$ The result
is also reminiscent of [Amb, Lemma 1.1].

Before stating the Simmons-Mospan specialization to the Haar context and
also to motivate one of the conditions in its subsequent generalizations, we
cite (and give a direct proof) of the following known result (equivalence of
Haar measure $\eta $ and its inverse $\tilde{\eta}$), encapsulated in the
formula\textbf{\ }
\[
\eta (K^{-1})=\int_{K}d\eta (t)/\Delta (t),
\]%
exhibiting the direct connection between $\eta $ and $\tilde{\eta}$ via the
modular function [HewR, 15.14], or [Hal, \S 60.5f]; this equivalence result
holds more generally between any two probability measures when one is left
and the other right quasi-invariant -- see [Xia, Cor. 3.1.4]; this is
related to atheorem of Mackey's [Mac] cf. \S 9.16. As will be seen from the
proof, in Lemma H\ below there is no need to assume the group is separable,
a compact metrizable subspace (being totally bounded) is separable.

\bigskip

\noindent \textbf{Lemma H }(cf. [Hal, \S 50(ff); \S 59 Th. D]). \textit{In a
locally compact metrizable group }$G,$\textit{\ for }$K$ \textit{compact, if
}$\eta (K)=0,$\textit{\ then }$\eta (K^{-1})=0,$ \textit{and, by regularity,
so also for }$K$ \textit{measurable.}

\bigskip

\noindent \textbf{Proof. }Fix a compact $K.$ As $K$ is compact, $\Delta $
(the modular function of $G)$ is bounded away from $0$ on $K,$ say by $M>0;$
furthermore, $K$ is separable, so pick $\{d_{n}:n\in \mathbb{N}\}$ dense in $%
K$. Then for any $\varepsilon >0$ there are two (finite) sequences $%
m(1),...m(n)\in \mathbb{N}$ and $\delta (1),...\delta (n)>0$ such that $%
\{B_{\delta (i)}d_{m(i)}:i\leq n\}$ covers $K$ and%
\[
M\sum\nolimits_{i\leq n}\eta (B_{\delta (i)})\leq \sum\nolimits_{i\leq
n}\eta (B_{\delta (i)})\Delta (d_{m(i)})=\sum\nolimits_{i\leq n}\eta
(B_{\delta (n)}d_{m(n)})<\varepsilon .
\]%
Then%
\[
\sum\nolimits_{i\leq n}\eta (d_{m(i)}^{-1}B_{\delta
(i)})=\sum\nolimits_{i\leq n}\eta (B_{\delta (i)})\leq \varepsilon /M.
\]%
But $\{d_{m(i)}^{-1}B_{\delta (i)}:i\leq n\}$ covers $K^{-1}$ by the
symmetry of the balls $B_{\delta }$ (by the symmetry of the norm); so, as $%
\varepsilon >0$ is arbitrary, $\eta (K^{-1})=0.$

As for the final assertion, if $\eta (E^{-1})>0$ for some measurable $E,$
then $\eta (K^{-1})>0$ for some compact $K^{-1}\subseteq E^{-1}$, by
regularity; then $\eta (K)>0,$ and so $\eta (E)>0$. $\square $

\bigskip

Proposition 8 and Lemma H immediately give:

\bigskip

\noindent \textbf{Theorem S }(cf. [Sak, III.11], [Mos], [BarFF, Th. 7])%
\textbf{. }\textit{For }$G\ $\textit{locally compact with left Haar measure }%
$\eta $\textit{\ and }$\nu $\textit{\ a Borel measure on }$G,$\textit{\ if
the set }$S$\textit{\ is }$\eta $\textit{-null, then for }$\eta $\textit{%
-almost all }$t$%
\[
\nu (tS)=0.
\]%
\textit{In particular, this is so for }$S$\textit{\ the support of a measure
}$\nu $\textit{\ singular with respect to }$\eta .$

\bigskip

This in turn allows us to prove the locally compact (separable) case of the
Simmons-Mospan Theorem, as stated in the Introduction. We then pursue a
non-locally compact variant.

\bigskip

\noindent \textbf{Proof of Theorem SM. }If $\mu $ is absolutely continuous
w.r.t. Haar measure $\eta ,$ then $\mu ,$ being invariant, is subcontinuous,
and Lemma 1 (in \S 2) gives the Steinhaus-Weil property. Otherwise,
decomposing $\mu $ into its singular and absolutely continuous parts w.r.t. $%
\eta $, choose $K$ a compact subset of the support of the singular part of $%
\mu $; then \textit{\ }$\mu (K)>\mu _{-}(K)=0,$ by Prop. SM above, and so
Prop. 6 (Converse part -- see \S 2) applies. $\square $

\bigskip

\noindent \textbf{Proposition 9 }(after Simmons, cf. [Sim, Lemma] and
[BarFF, Th. 8]). \textit{For }$G$ \textit{\ a Polish group, }$\mu ,\nu \in
\mathcal{M}(G)$\textit{\ and }$\nu \bot ^{\text{L-inv}}\mu $\textit{\
concentrated on a compact invertibly }$\mu $\textit{-null set }$K,$ \textit{%
there is }$B\subseteq K$\textit{\ such that }$K\backslash B$\textit{\ is }$%
\nu $\textit{-null and both }$BB^{-1}$\textit{\ and }$B^{-1}B$ \textit{have
empty interior.}

\bigskip

\noindent \textbf{Proof. }As we are concerned only with the subspace $%
KK^{-1}\cup K^{-1}K$, w.l.o.g. the group $G$ is separable. By Prop. 8, $%
Z:=\{x:\nu (xK)=0\}$ is dense and so also%
\[
Z_{1}:=\{x:\nu (K\cap xK)=0\},
\]%
since $\nu (K\cap xK)\leq \nu (xK)=0,$ so that $Z\subseteq Z_{1}.$ Take a
denumerable dense set $D\subseteq Z_{1}$ and put%
\[
S:=\bigcup\nolimits_{d\in D}K\cap dK.
\]%
Then $\nu (S)=0.$ Take $B=K\backslash S.$ If $\emptyset \neq V\subseteq
BB^{-1}$ and $d\in D\cap V,$ then for some $b_{1},b_{2}\in B\subseteq K$%
\[
d=b_{1}b_{2}^{-1}:\qquad b_{1}=db_{2}\in K\cap dK\subseteq S,
\]%
a contradiction, since $B\cap S=\emptyset .$ So $(K\backslash S)(K\backslash
S)^{-1}$ has empty interior. A similar argument based on
\[
T:=\bigcup\nolimits_{d\in D}Kd\cap K
\]%
ensures that also $(K\backslash S\backslash T)^{-1}(K\backslash S\backslash
T)$ has empty interior. $\square $

\bigskip

In order to generalize the Simmons Theorem from its locally compact context
we will need to cite the following result. Here $\mathbb{Q}_{+}:=\mathbb{Q}%
\cap (0,\infty )$ denotes the positive rationals, and $B_{\delta }^{K,\Delta
}(\sigma ):=\{z\in B_{\delta }:\sigma (Kz)>\Delta \}$ as in \S 2.

\bigskip

\noindent \textbf{Theorem 5 (Disaggregation Theorem, }[BinO8, Th.6, Prop. 4]%
\textbf{). }\textit{Let }$G\ $\textit{be a Polish group that is strongly
amenable at 1, and let} $\mathbf{t}$ \textit{be a regular null sequence. For}
$\sigma =\sigma (\mathbf{t})$ \textit{there are a countable family }$%
\mathcal{H}$ \textit{with }$\mathcal{H\subseteq K}_{+}(\sigma ),$ \textit{a
countable set} $D=D(\mathcal{H})\subseteq G$ \textit{dense in }$G,$ \textit{%
and a dense subset }$G(\sigma )$\textit{\ on which the sets below are the
sub-basic sets of a metrizable topology:}%
\[
B_{\delta }^{K,\Delta }(\sigma )\qquad (K\in \mathcal{H},\delta ,\Delta \in
\mathbb{Q}_{+},\mathit{\ }\Delta <\sigma (K)).
\]%
\textit{In particular, the space }$G(\sigma )$ \textit{is continuously and
compactly embedded in }$G.$ \textit{Moreover, each such open set contains a
cofinal subsequence of} $\mathbf{t}$.

\bigskip

For a proof we refer the reader to the companion paper [BinO8]; the result
relies on Corollary 4 above. The subspace $G(\sigma )$ here is a topological
analogue of the Cameron-Martin subspace $H(\gamma )$ of a locally convex
topological vector space equipped with a Radon Gaussian measure $\gamma $ --
see \S 8.2-3.

We are now ready for the promised generalization. This requires equivalence
of the Solecki measure and its inverse -- valid at least in Polish abelian
groups (see Th. 4 of \S 2 on strong amenability at 1).

\bigskip

\noindent \textbf{Theorem 6 }(\textbf{Generalized Simmons Theorem}, cf.
[Sim, Th. 2])\textbf{. }\textit{Let }$G\ $\textit{be a Polish group that is
strongly amenable at 1 (e.g. if }$G$\textit{\ is abelian), let }$\sigma
=\sigma (\mathbf{t})$ \textit{be a Solecki measure corresponding to a
regular null sequence} $\mathbf{t},$ \textit{which we assume is equivalent
to its inverse }$\tilde{\sigma}$ \textit{(e.g. if }$G$\textit{\ is abelian),}
\textit{and let }$G(\sigma )$ \textit{be the dense subgroup of the preceding
theorem. Then:}\newline
$\nu \in \mathcal{M}(G)$ \textit{is left invertibly-singular w.r.t. }$\sigma
$\textit{\ iff }$\nu $\textit{\ has a support that is a }$\sigma $\textit{%
-compact union of compact sets }$K_{n}$\textit{\ with each of the compact
sets }$K_{n}K_{n}^{-1}$\textit{\ and }$K_{n}^{-1}K_{n}$\textit{\ nowhere
dense (equivalently: having empty interior) in the topology of the subgroup }%
$G(\sigma )$\textit{, as above.}

\bigskip

\noindent \textbf{Proof. }Suppose that $\mu \in \mathcal{M}(G)$ is
subcontinuous. If $\nu \in \mathcal{M}(G)$, by Theorem LD write%
\[
\nu =\nu _{\text{a}}+\nu _{\text{s}}\text{ with }\nu _{\text{a}}<^{\text{%
L-inv}}\mu \text{ and }\nu _{\text{s}}\bot ^{\text{L-inv}}\mu .
\]%
If $\nu $ is concentrated as in the statement of the theorem on a $\sigma $%
-compact set $B$ with $B^{-1}B$ having empty interior in $G(\sigma )$, then $%
\nu _{\text{a}}=0,$ and so $\nu $ is left invertibly-singular w.r.t. $\mu .$
Indeed, as $\nu $ is concentrated on $B,$ so is $\nu _{\text{a}}.$ We claim
that $\nu _{\text{a}}(B)=0.$ Otherwise, $\nu _{\text{a}}(K_{n})>0$ for some
compact $K_{n}\subseteq B.$ So $K=K_{n}\notin \mathcal{M}_{0}^{\text{L-inv}%
}(\sigma ),$ as $\nu _{\text{a}}<^{\text{L-inv}}\sigma .$ The argument now
splits into two cases, according as $K\notin \mathcal{M}_{0}^{\text{L}%
}(\sigma )$ or $K^{-1}\notin \mathcal{M}_{0}^{\text{L}}(\sigma )$

First, suppose that $\sigma (gK)>0$ for some $g\in G;$ then, by Lemma 1,
there are $\delta >0$ and $0<\Delta <\sigma _{-}^{\mathbf{t}}(gK)$ with
\[
B_{\delta }^{gK\Delta }(\sigma )\subseteq (gK)^{-1}gK=K^{-1}K\subseteq
B^{-1}B,
\]%
contradicting the above property of $B.$

Next, suppose that $\sigma (Kg)=\sigma ^{-1}(g^{-1}K^{-1})>0$ for some $g\in
G;$ so $\sigma (g^{-1}K^{-1})>0,$ as $\tilde{\sigma}$ is equivalent to $%
\sigma .$ Then, again by Lemma 1, there are $\delta >0$ and $0<\Delta
<\sigma _{-}^{\mathbf{t}}(g^{-1}K^{-1})$ with
\[
B_{\delta }^{g^{-1}K^{-1}\Delta }(\sigma )\subseteq
(g^{-1}K^{-1})^{-1}g^{-1}K^{-1}=KK^{-1}\subseteq BB^{-1},
\]%
again contradicting the above property of $B.$ So $\nu =\nu _{\text{a}}$ is
invertibly singular w.r.t. $\mu .$

The rest of the proof is as in Simmons [Sim, Th. 2], using Prop. 6: the
Baire-category argument still holds, since compactness implies closure under
the $G(\sigma (\mathbf{t}))$-topology, the latter being a finer topology;
avoidance of interior points requires second countability, assured by Th. 5.
$\square $

\bigskip

\noindent \textbf{Corollary 5 (Simmons Theorem}: [Sim, Th. 2])\textbf{. }%
\textit{For }$G$\textit{\ separable and locally compact and }$\eta $ \textit{%
left Haar measure:}\newline
$\nu \in \mathcal{M}(G)$ \textit{is singular w.r.t. }$\eta $\textit{\ iff }$%
\nu $\textit{\ has a support that is a }$\sigma $\textit{-compact union of
compact sets }$K_{n}$\textit{\ with each of the compact sets }$%
K_{n}K_{n}^{-1}$\textit{\ nowhere dense (equivalently: having empty
interior).}

\bigskip

For the non-separable version of the above, see \S 8.1.

\part{Weil topologies}

\section{Weil-like topologies: preliminaries}

We turn now to relatives of the \textit{Weil topology}. For background, we
refer to Weil's book [Wei, Ch. VII] and Halmos's book [Hal, Ch. XII] (see
also \S 8.4). Weil regarded his result as a \textit{Converse Haar Theorem, }%
in retrieving the topological-group structure from the measure-algebra
structure [Fre] as encoded by the Haar-measurable subsets -- cf. [Kod].
(Here one may work either, following Weil, to within a dense embedding in a
locally compact group, as in the Remark to Theorem 7M below, or, following
Mackey, uniquely up to homeomorphism, granted the further assumption of an
analytic Borel structure [Mac, Th. 7.1]; for further information see \S %
8.16.) The alternative view below throws light on this result in that the
measure structure is already encoded by the density topology $\mathcal{D}$
via the Haar density theorem, for which see [Mue], [Hal, \S 61(5), p. 268],
cf. [BinO2, \S 7; Th. 6.10], [BinO5]; this view is partially implicit in
[Amb], where refinement of one invariant measure $\mu _{1}$ by another $\mu
_{2}$ holds when sets in $\mathcal{M}_{+}(\mu _{2})$ contain sets in $%
\mathcal{M}_{+}(\mu _{1})$ (as in the refinement of one topology by
another). This falls within the broader aim of retrieving a \textit{%
topological} group structure from a given (one-sidedly) invariant topology $%
\tau $ on a group $G,$ when $\tau $ arises from refinement of a topological
group structure (i.e. starting from a\textit{\ semitopological} group
structure $(G,\tau )).$ Also relevant here are \textit{Converse
Steinhaus-Weil} results, as in Prop. 7 of \S 3\ above (see also \S 8.4). For
background on group-norms see the textbook treatment in [ArhT, \S 3.3] (who
trace this notion back to Markov) or [BinO2], but note their use of
`pre-norm' for what we call (following Pettis [Pet]) a \textit{pseudo-norm};
for quasi-interiors and regular open sets see \S 8.6. Thus a norm $||\cdot
||:G\rightarrow \lbrack 0,\infty )$ satisfies all the three conditions 1-3
below and generates a right-invariant metric $d(x,y)=||xy^{-1}||$ and so a
topology $\mathcal{T=T}_{d},$ just as a right-invariant metric $d$ derives
from a\ separable topology $\mathcal{T}_{G}$ and generates, via the
Birkhoff-Kakutani Theorem ([HewR, Th. 8.3], [Gao, Th. 2.1.1]), the norm $%
||x||=d(x,1_{G}).$ A pseudo-norm differs in possibly lacking condition 1.i.
(so generates a \textit{pseudo}-metric).

1.i (positivity): $||g||>0$ for $g\neq 1_{G}$, and 1.ii: $||1_{G}||=0;$

2 (subadditivity): $||gh||\leq ||g||+||h||,$

3 (symmetry): $||g^{-1}||=||g||.$

\bigskip

Recall that a set function $\lambda $ defined on $\mathcal{U}(G)$ is a
\textit{submeasure }if it is monotone and subadditive with $\lambda
(\emptyset )=0$ (Introduction, [Fre, Ch. 39, \S 392], [Tal]); by analogy
with the term \textit{finitely additive measure }(for background see [Bin],
[TomW, Ch. 12]; cf. [Pat]), this is a \textit{finitely subadditive outer
measure}, similarly as in Maharam [Mah], albeit in the context of Boolean
algebras, but without her positivity condition. Recall from Halmos [Hal, Ch.
II \S 10] that a submeasure is an \textit{outer measure} if in addition it
is \textit{countably subadditive}. The set function $\lambda $ is \textit{%
left invariant} if $\lambda (gE)=\lambda (E)$ for all $g\in G$ and $E\in
\mathcal{U}(G).$

Propositions 10 and 11 below are motivated by [Hal, Ch. XII \S 62, cf. Ch.
II \S 9 (2-4)], where $G$ is a locally compact group with $\lambda $ its
left Haar measure, but here the context is broader, allowing in \textit{%
amenable} groups $G$ (cf. [TomW, Ch. 12], [Pat]). The two results enable the
introduction in \S 6 of Weil-like topologies generated from families of
left-invariant \textit{pseudo-metrics} derived from invariant submeasures.
The latter rely on the natural \textit{measure-metric}, also known as the
\textit{Fr\'{e}chet-Nikodym metric} ([Fre, \S 323Ad], [Hal, \S 40 Th. A],
[Bog2, p. 53, 102-3, 418]); see [Drew1,2] (cf. [Web]) for the related
literature of Fr\'{e}chet-Nikodym topologies and their relation to the
Vitali-Hahn-Saks Theorem. Maharam [Mah] studies sequential continuity of the
order relation (of inclusion, here in the measure algebra), and requires
positivity to obtain a (measure-) \textit{metric}; see Talagrand [Tal] (cf.
[Fre, \S 394] and the literature cited there) for a discussion of
pathological submeasures (the only measures they dominates under $\ll $
being trivial), and [ChrH] for corresponding exotic abelian Polish groups.

In the setting of a locally compact group $G,$ these pseudo-metrics are
implicit in work of Struble: initially, in 1953 [Str1], he used a
(`sampler') family of pre-compact open sets $\{E_{t}:t>0\}$ to construct a
mean on $G,$ thereby refering to a one-parameter family of pseudo-metrics
corresponding to the sets $E_{t}$; some twenty years later in 1974 [Str2]
(cf. [DieS, Ch. 8]) identifies a left-invariant (proper) metric on $G$ by
taking the supremum of pseudo-metrics, each generated from some open set in
a countable open base at $1_{G}.$ The pseudo-metric makes a very brief
appearance in Yamasaki's textbook treatment [Yam, Ch. 1] of Weil's theorem.

\bigskip

\noindent \textbf{Proposition 10} (\textbf{Weil pseudo-norm}, cf. [Fre, \S %
392H], [Yam, Ch. 1, Proof of Th. 4.1]). \textit{For }$G$\textit{\ a Polish
group, }$\lambda \in \mathcal{M}_{\text{sub}}(G),$\textit{\ a left-invariant
submeasure on }$\mathcal{U}(G)$\textit{, and }$E\in \mathcal{U}(G)$\textit{\
with }$\lambda (E)>0,$\textit{\ put}%
\[
||g||_{E}^{\lambda }:=\lambda (gE\triangle E)\qquad (g\in G).
\]%
\textit{Then }$||.||_{E}$\textit{\ defines a group pseudo-norm with
associated right-invariant pseudo-metric}%
\[
d_{E}^{\lambda }(g,h)=||gh^{-1}||_{E}^{\lambda }\qquad (g,h\in G).
\]%
\textit{Likewise, for }$\lambda $ \textit{right-invariant, a pseudo-norm is
defined by}%
\[
||g||_{E}^{\lambda }:=\lambda (E\triangle Eg)\qquad (g\in G).
\]

\bigskip

\noindent \textbf{Proof.} Since $\lambda (\emptyset )=0,$ $%
||1_{G}||_{E}^{\lambda }=0.$ By left invariance under $a,$%
\[
||a^{-1}||_{E}^{\lambda }=\lambda (a^{-1}E\triangle E)=\lambda
(a(a^{-1}E\triangle E))=\lambda (E\triangle aE)=||a||_{E}^{\lambda }.
\]%
Also,
\[
||ab||_{E}^{\lambda }\leq ||a||_{E}^{\lambda }+||b||_{E}^{\lambda }
\]%
follows from monotonicity, subadditivity and $\lambda (abE\triangle
aE)=\lambda (bE\triangle E):$%
\begin{eqnarray*}
\lambda (abE\backslash E\cup E\backslash abE) &\leq &\lambda (abE\backslash
aE)\cup (aE\backslash E)\cup (E\backslash aE)\cup (aE\backslash abE)) \\
&=&\lambda (abE\backslash aE)\cup (aE\backslash abE)\cup (aE\backslash
E)\cup (E\backslash aE)) \\
&\leq &\lambda (abE\triangle aE)+\lambda (E\triangle aE)=\lambda
(bE\triangle E)+\lambda (E\triangle aE).\qquad \square
\end{eqnarray*}

\bigskip

\noindent \textbf{Corollary 6 }(Kneser for Haar measure, [Kne, Hilfs. 4]).
\textit{For }$G$\textit{\ a Polish group, }$\lambda \in \mathcal{M}_{\text{%
sub}}(G),$\textit{\ a left-invariant submeasure on }$\mathcal{U}(G)$\textit{%
, and }$E\in \mathcal{U}(G)$\textit{\ with }$\lambda (E)>0,$ \textit{the set}%
\[
H:=\{g\in G:\lambda (gE\triangle E)=0\}
\]%
\textit{is a subgroup of }$G$ \textit{closed under the norm }$%
||g||_{E}^{\lambda }$\textit{.}

\bigskip

\noindent \textbf{Proof.} Indeed $H=\{g\in G:||g||_{E}^{\lambda }=0\},$ and
so $H$ is a subgroup, since for $g,h\in H,$ $||gh^{-}||_{E}^{\lambda }\leq
||g||_{E}^{\lambda }+||h||_{E}^{\lambda }=0.$ $\square $

\bigskip

Recall now (from the Introduction) that a subset of a Polish group $G$ is
\textit{left Haar null} if it is contained in a universally measurable set $%
B $ such that for some $\mu \in \mathcal{P}(G)$%
\[
\mu (gB)=0\qquad (g\in G).
\]%
It is \textit{Haar null} [Sol1] (cf. [HofT, p. 374]) if it is contained in a
universally measurable set $B$ such that for some $\mu \in \mathcal{P}(G)$%
\[
\mu (gBh)=0\qquad (g,h\in G).
\]%
This motivates the following application of Proposition 10 beyond Haar
measure. Extending the notation of \S 3, below $\mathcal{M}_{0}^{L}(G)$
(resp. $\mathcal{M}_{0}(G)$) denotes the family of left-Haar-null (resp.
Haar-null) sets of $G$, and we write%
\[
\mathcal{U}_{+}^{L}(G):=\mathcal{U}(G)\backslash \mathcal{M}%
_{0}^{L}(G),\qquad \mathcal{U}_{+}(G):=\mathcal{U}(G)\backslash \mathcal{M}%
_{0}(G).
\]

Prop. 10 may be applied to the following measures; those constructed from $%
\mu $ a normalized counting measure (of finite support) are studied in
[Sol1].

\bigskip

\noindent \textbf{Proposition 11.} \textit{In a Polish group }$G,$\textit{\
for }$\mu \in \mathcal{P}(G)$\textit{\ put}%
\begin{eqnarray*}
\mu _{L}^{\ast }(E) &:&=\sup \{\mu (gE):g\in G\}\qquad (E\in \mathcal{U}(G)),
\\
\hat{\mu}(E) &:&=\sup \{\mu (gEh):g,h\in G\}\qquad (E\in \mathcal{U}(G)).
\end{eqnarray*}%
\textit{Then }$\mu _{L}^{\ast }$ \textit{(resp. }$\hat{\mu}$\textit{) is a
left invariant (resp. bi-invariant) submeasure, which is positive for }$E\in
\mathcal{U}_{+}^{L}(G)$ \textit{(resp. for }$E\in \mathcal{U}_{+}(G)$\textit{%
), i.e. for universally measurable, non-left-Haar null (resp. non-Haar null)
sets.}

\bigskip

\noindent \textbf{Proof}. We consider only $\hat{\mu}$, as the case $\mu
_{L}^{\ast }$ is similar and simpler (through the omission of $h$ and $b$
below). The set function $\hat{\mu}$ is well defined, with%
\[
\mu (E)\leq \hat{\mu}(E)\leq 1\qquad (E\in \mathcal{U}(G)),
\]%
since $\mu $ is a probability measure; it is bi-invariant, since%
\[
\hat{\mu}(aEb):=\sup \{\mu (gaEbh):g,h\in G\}=\sup \{\mu (gEh):g,h\in G\},
\]%
and $G$ is a group. Furthermore, for $B\in \mathcal{U}(G)$%
\[
\mu (gBh)\leq \hat{\mu}(B)\leq 1,\qquad (g,h\in G).
\]%
So, for $\mu \in \mathcal{P}(G)$
\[
0<\hat{\mu}(B)\leq 1\qquad (B\in \mathcal{U}_{+}(G)),
\]%
since there are $g,h\in G\ $with $\mu (gBh)>0.$ Countable subadditivity
follows (on taking suprema of the leftmost term over $g,h$) from%
\[
\mu (g(\bigcup\nolimits_{n}A_{n})h)\leq \sum\nolimits_{n}\mu (gA_{n}h)\leq
\sum\nolimits_{n}\hat{\mu}(gA_{n}h)=\sum\nolimits_{n}\hat{\mu}(A_{n}),
\]%
for any sequence of sets $A_{n}\in \mathcal{U}(G)$. $\square $

\bigskip

\noindent \textbf{Definition. }For $\mu \in \mathcal{P}(G),E\in \mathcal{U}%
(G),$ put
\[
B_{\varepsilon }^{E}(\mu ):=\{x\in G:||x||_{E}^{\mu }<\varepsilon \}.
\]%
Our next step uses Prop. 11 to inscribe these balls into $EE^{-1}$ for all
small enough $\varepsilon >0$.

\bigskip

\noindent \textbf{Lemma 2 }(\textbf{Self-intersection Lemma}). \textit{In a
Polish group }$G$ \textit{for }$E\in \mathcal{U}_{+}(G),$ \textit{and
respectively for }$E\in \mathcal{U}_{+}^{L}(G),$ \textit{and }$\mu \in
\mathcal{P}(G),$%
\[
\left.
\begin{array}{cc}
1_{G}\in B_{\varepsilon }^{E}(\hat{\mu})\subseteq EE^{-1} & \qquad
(0<\varepsilon <\hat{\mu}(E)), \\
1_{G}\in B_{\varepsilon }^{E}(\mu _{L}^{\ast })\subseteq EE^{-1} & \qquad
(0<\varepsilon <\mu _{L}^{\ast }(E)).%
\end{array}%
\right\}
\]%
\textit{Equivalently, for }$0<\varepsilon <\hat{\mu}(E),$\textit{\ and
respectively for }$0<\varepsilon <\mu _{L}^{\ast }(E)$\textit{, }%
\[
E\cap xE\neq \emptyset \qquad (x\in B_{\varepsilon }^{E}(\hat{\mu}));\qquad
E\cap xE\neq \emptyset \qquad (x\in B_{\varepsilon }^{E}(\mu _{L}^{\ast })).
\]

\noindent \textbf{Proof. }We check only the $\hat{\mu}$ case; the other is
similar and simpler (through the omission of $h$ below). For $E\in \mathcal{U%
}_{+}(G),$ since $\hat{\mu}(E)>0$ by Prop. 11, we may pick $g,h\in G$ such
that $\varepsilon _{E}:=\mu (gEh)>0.$ Consider $x$ and $\varepsilon >0$ with
$||x||_{E}^{\hat{\mu}}<\varepsilon \leq \varepsilon _{E}.$ If $E$ and $xE$
are disjoint, then
\begin{eqnarray*}
\varepsilon _{E} &=&\mu (gEh)\leq \mu (g(E\cup xE)h)\leq \hat{\mu}(g(E\cup
xE)h)=\hat{\mu}(E\cup xE) \\
&=&\hat{\mu}(xE\triangle E)=||x||_{E}^{\hat{\mu}}<\varepsilon \leq
\varepsilon _{E},
\end{eqnarray*}%
a contradiction. So $E$ and $xE$ do meet. Now first pick $t\in xE\cap E\ $%
and next $s\in E$ so that $t=xs;$ then $x=ts^{-1}\in EE^{-1}.$ The argument
is valid when $\varepsilon _{E}=\mu (gEh)$ assumes any value in $(0,\hat{\mu}%
(E)]$, as in Prop. 1M. For the converse, if $B_{\varepsilon }^{E}(\hat{\mu}%
)\subseteq EE^{-1},$ proceed as in Prop. 1M. $\square $

\bigskip

We need a simple analogue of a result due to Weil ([Wei, Ch. VII, \S 31],
cf. [Hal, Ch. XII \S 62]). Below $\tau _{1}$ denotes the $\tau $-open\textit{%
\ neighbourhoods of} $1_{G}$. For $G$ locally compact with $\lambda =\eta
_{G}$, the identity%
\begin{equation}
2\eta (E)-2\eta (E\cap xE)=\eta (E\triangle xE)=1-2\int
1_{E}(t)1_{E^{-1}}(t^{-1}x)d\eta (t)  \tag{$\dag $}
\end{equation}%
connects the continuity of the (pseudo-) norm to $\mathcal{T}_{d}$%
-continuity of translation in the topological group structure $(G,\mathcal{T}%
_{d})$ of the locally compact group, and to continuity of the convolution
function here (for $E$ of finite $\eta $-measure) -- see [HewR, Th. 20.16];
see also [HewR, Th. 20.17] for the well-known connection between the
Steinhaus-Weil Theorem and convolution. Such a continuity condition result
guarantees that $B_{\varepsilon }^{E}(\eta )$ contains points other than $%
1_{G}.$

\bigskip

\noindent \textbf{Lemma 3 }(\textbf{Fragmentation Lemma}; cf. [Hal, Ch. XII
\S 62 Th. A]). \textit{For }$\lambda \in \mathcal{M}_{\text{sub}}(G)$\textit{%
\ a left-invariant submeasure on a Polish group }$G$\textit{\ equipped with
a finer right-invariant topology }$\tau $\textit{\ with }$1_{G}$-\textit{%
open-nhd family }$\tau _{1}\subseteq \mathcal{U}_{+}^{L}(G)$\textit{:}

\textit{if the map }%
\[
x\mapsto ||x||_{E}^{\lambda }
\]%
\textit{is continuous under }$\tau $ \textit{at }$x=1_{G}$ \textit{for each }%
$E\in \mathcal{U}_{+}^{L}(G)$

\noindent \textit{-- then, for each }$\emptyset \neq E,F\in \tau $\textit{\
and }$\varepsilon >0$ \textit{with }$\varepsilon <\lambda (E)$\textit{,
there exists }$H\in \tau _{1}$\textit{\ with }$HH^{-1}\subseteq FF^{-1}$%
\textit{\ and}%
\[
||h^{\prime }h^{-1}||_{E}^{\lambda }<\varepsilon \text{\qquad }(h,h^{\prime
}\in H):\text{\qquad }HH^{-1}\subseteq B_{\varepsilon }^{E},
\]%
\textit{so that }\textrm{diam}$_{E}^{\lambda }(H)\leq \varepsilon .$

\bigskip

\noindent \textbf{Proof.} Pick any $f\in F,$ and $D\in \tau _{1}$ satisfying
$||x||_{E}^{\lambda }<\varepsilon /2$ for all $x\in D$. As $\tau $ is
right-invariant and $1_{G}\in D\cap Ff^{-1}\in \tau $, pick $H\in \tau _{1}$
with $H\subseteq D\cap Ff^{-1}$; then%
\[
HH^{-1}=Hff^{-1}H^{-1}\subseteq FF^{-1}.
\]%
For $h,h^{\prime }\in H,$ as $h,h^{\prime }\in D$%
\[
||h^{\prime }f(hf)^{-1}||_{E}^{\lambda }=||h^{\prime }h^{-1}||_{E}^{\lambda
}\leq ||h^{\prime }||_{E}^{\lambda }+||h^{-1}||_{E}^{\lambda }=||h^{\prime
}||_{E}^{\lambda }+||h||_{E}^{\lambda }<\varepsilon .\qquad \square
\]

\bigskip

In the presence of a refinement topology $\tau $ on the group $G,$ the lemma
motivates further notation: write $\mathcal{P}_{\text{cont}}(G,\tau ),$ or
just%
\[
\mathcal{P}_{\text{cont}}(\tau ):=\{\mu \in \mathcal{P}(G,\mathcal{T}%
_{d}):g\mapsto ||g||_{E}^{\hat{\mu}}:=\hat{\mu}(gE\triangle E)\text{ is }%
\tau \text{-continuous at }1_{G}\}.
\]

Of necessity attention here focuses on continuity. The characterization
question as to which topologies $\tau $ yield a non-empty $\mathcal{P}(\tau
) $ is in part answered by Theorem 7M below. Indeed, for Haar measure $\eta $
in the locally compact case,
\[
\mu \in \mathcal{P}_{\text{cont}}(\tau )\qquad (\mu \ll \eta ,\tau \supseteq
\mathcal{T}_{d}),
\]%
by ($\dag $) in the presence of $d\mu /d\eta $ as a kernel:%
\begin{equation}
||x||_{E}^{\mu }=1-2\int 1_{E}(t)1_{E^{-1}}(t^{-1}x)\frac{d\mu }{d\eta }%
d\eta (t).  \tag{$\dag \dag $}
\end{equation}%
However, $\mathcal{P}(G)$ will contain measures $\mu $ singular with respect
to $\eta :$ for such $\mu $, by Th. SM there will be Borel subsets $B$ of
positive $\mu $-measure such that $BB^{-1}$ has void $\mathcal{T}_{d}$%
-interior.

\section{Weil-like topologies: theorems}

Prop. 11 now yields the following result, which embraces known Hashimoto
topologies [BinO5] in both the Polish abelian setting, where the left Haar
null sets form a $\sigma $-ideal (Christensen [Chr1]), and likewise in (the
not necessarily abelian) Polish groups that are \textit{amenable at }$1$
(Solecki [Sol1]); this includes, as additive groups, $F$- (hence also
Banach) spaces -- cf. [BinO5,6], where use is made of Hashimoto topologies.

\bigskip

\noindent \textbf{Theorem 7}.\textit{\ Let }$G$\textit{\ be a Polish group\
and }$\tau $ \textit{both a left- and a right-invariant refinement topology
with }$1_{G}$\textit{-open-nhd} \textit{family\ }$\tau _{1}\subseteq
\mathcal{U}_{+}(G)$\textit{.}

\noindent \textit{Then both the families }$\{AA^{-1}:$\textit{\ }$A\in \tau
_{1}\}$ \textit{and }$\{B_{\varepsilon }^{E}(\hat{\mu}):\emptyset \neq E\in
\tau ,\mu \in \mathcal{P}(\tau )$\textit{\ and }$0<\varepsilon \leq \hat{\mu}%
(E)\}$ \textit{generate neighbourhoods of the identity under which }$G$
\textit{is a topological group. Moreover, the pseudo-norms}%
\[
\{||.||_{E}^{\hat{\mu}}:\emptyset \neq E\in \tau ,\mu \in \mathcal{P}_{\text{%
cont}}(\tau )\}
\]%
\textit{are downward directed by refinement as follows: for }$\emptyset \neq
E,F\in \tau _{1}$\textit{, }$\lambda ,\mu \in \mathcal{P}(\tau )$ \textit{and%
} $\varepsilon <\min \{\hat{\lambda}(E),\hat{\mu}(F)\}\},$\textit{\ there is
}$H\in \tau _{1}$ \textit{such that for }$0<\delta <\min \{\tilde{\lambda}%
(H),\hat{\mu}(H)\}$\textit{\ }%
\[
B_{\delta }^{H}(\lambda )\cap B_{\delta }^{H}(\mu )\subseteq B_{\varepsilon
}^{E}(\lambda )\cap B_{\varepsilon }^{F}(\mu ).
\]

\bigskip

\noindent \textbf{Proof.} The proof is similar to but simpler than that of
[Hal, Ch. XII \S 62 Th. A]. Given two (non-left-Haar-null) sets $E,F\in \tau
_{1}$ and $\varepsilon <\min \{\hat{\lambda}(E),\hat{\mu}(F)\},$ by the
Fragmentation Lemma (Lemma 3 of \S 5) applied separately to $\hat{\lambda}$
and to $\hat{\mu},$ there are $A,B\in \tau _{1}$ with
\[
AA^{-1}\subseteq B_{\varepsilon }^{E}(\hat{\lambda}),\quad BB^{-1}\subseteq
B_{\varepsilon }^{F}(\hat{\mu}).
\]%
Take any $H\in \tau _{1}$ with $H\subseteq A\cap B;$ then
\[
HH^{-1}\subseteq AA^{-1}\cap BB^{-1}.
\]%
Since $H\in \mathcal{U}_{+}(G)$ (as $\tau _{1}\subseteq \mathcal{U}_{+}(G)),$
take $\delta $ with $0<\delta <\min \{\hat{\lambda}(H),\hat{\mu}(H)\};$ then
by ($\ast $)%
\[
B_{\delta }^{H}(\hat{\lambda})\cap B_{\delta }^{H}(\hat{\mu})\subseteq
HH^{-1}\subseteq AA^{-1}\cap BB^{-1}\subseteq B_{\varepsilon }^{E}(\hat{%
\lambda})\cap B_{\varepsilon }^{F}(\hat{\mu}).
\]%
(So `mutual refinement' holds between the sets of the form $AA^{-1}$ and
those of the form $B_{\varepsilon }^{E}$.) As $||\cdot ||_{E}^{\hat{\mu}}$
is a pre-norm,%
\[
B_{\varepsilon /2}^{E}(\hat{\mu})B_{\varepsilon /2}^{E}(\hat{\mu}%
)^{-1}=B_{\varepsilon /2}^{E}(\hat{\mu})B_{\varepsilon /2}^{E}(\hat{\mu}%
)\subseteq B_{\varepsilon }^{E}(\hat{\mu}).
\]%
By the Fragmentation Lemma again, given any $x\in G$ and $\varepsilon >0,$
choose $H\in \tau _{1}$ with $HH^{-1}\subseteq B_{\varepsilon }^{E}(\tilde{%
\mu}).$ Then with $F:=xH\in \tau ,$%
\[
B_{\varepsilon }^{F}(\hat{\mu})=\{z:||z||_{F}^{\hat{\mu}}<\varepsilon
\}\subseteq (xH)(xH)^{-1}=xHH^{-1}x^{-1}\subseteq xB_{\varepsilon }^{E}(\hat{%
\mu})x^{-1}.
\]%
Finally, for any $x_{0}$ with $||x_{0}||_{E}^{\hat{\mu}}<\varepsilon ,$ put $%
\delta :=\varepsilon -||x_{0}||_{E}^{\hat{\mu}}.$ Then for $||y||_{E}^{\hat{%
\mu}}<\delta ,$
\[
||x_{0}\cdot y||_{E}^{\hat{\mu}}\leq ||x_{0}||_{E}^{\hat{\mu}}+||y||_{E}^{%
\hat{\mu}}<||x_{0}||_{E}^{\hat{\mu}}+\varepsilon -||x_{0}||_{E}^{\hat{\mu}%
}<\varepsilon ,
\]%
i.e.
\[
x_{0}B_{\delta }^{E}(\hat{\mu})\subseteq B_{\varepsilon }^{E}(\hat{\mu}%
).\qquad \square
\]

Specializing to locally compact groups yields as a corollary, on writing $%
B_{\varepsilon }^{E}:=B_{\varepsilon }^{E}(\eta ):$

\bigskip

\noindent \textbf{Theorem 7M}.\textit{\ For }$G$\textit{\ a locally compact
group with left Haar measure }$\eta ,$\textit{\ if:}

\noindent (i)\textit{\ }$\tau $ \textit{is both a left- and a
right-invariant refinement topology with }$\tau _{1}\subseteq \mathcal{M}%
_{+} $\textit{,}

\noindent (ii)\textit{\ for every non-empty }$E\in \tau ,$\textit{\ the
pseudo-norm}%
\[
g\mapsto ||g||_{E}:=\eta (gE\triangle E)\qquad (g\in G)
\]%
\textit{is continuous under }$\tau $\textit{\ at }$g=1_{G}$

\noindent -- \textit{then both the families }$\{AA^{-1}:$\textit{\ }$A\in
\tau _{1}$\} \textit{and }$\{B_{\varepsilon }^{E}:$\textit{\ }$\emptyset
\neq E\in \tau $ \textit{and }$0<\varepsilon \leq 2\eta (E)\}$ \textit{%
generate neighbourhoods of the identity under which }$G$ \textit{is a
topological group. Moreover, the pseudo-norms}%
\[
\{||.||_{E}:\emptyset \neq E\in \tau \}
\]%
\textit{are downward directed by refinement; indeed, for }$\emptyset \neq
E,F\in \tau $\textit{\ and} $\varepsilon <2\min \{\eta (E),\eta (F)\},$%
\textit{\ there is }$H\in \tau _{1}$ \textit{such that for }$0<\delta <\eta
(H)$\textit{\ }%
\[
B_{\delta }^{H}\subseteq B_{\varepsilon }^{E}\cap B_{\varepsilon }^{F}.
\]

\bigskip

\noindent \textbf{Proof.} It is enough to replace $\mathcal{P}(G)$ by $%
\{\eta \}$ (so that $\lambda $ and $\mu $ both refer to $\eta $), and to
note that if $xE$ and $E$ are disjoint, then $\eta (xE\triangle E)=2\eta
(E), $ so that in Lemma 2 the bound $\eta ^{\ast }(E)$ in the restriction
governing inclusion may be replaced by $2\eta (E)$. $\square $

\bigskip

\noindent \textbf{Remark. }As in [Hal, Ch. XII \S 62 Th. F], but by the
Fragmentation Lemma (and by the countable additivity of $\eta $), the
Weil-like topology on a locally compact $G$ in Theorem 7M is locally bounded
(norm-totally-bounded in some ball). Then $G$ with the Weil-like topology
may be densely embedded in its completion $\hat{G},$ which is in turn
locally compact, being locally complete and (totally) bounded. However, the
corresponding argument in the case of the preceeding more general Theorem 7
fails, since $\hat{\mu}$ there is not necessarily countably additive.

\bigskip

Finally, we give a category version of Theorem 7M, as an easy corollary;
indeed, our main task is merely to define what is meant by `mutatis
mutandis' in the present context. Given the assumption $\tau _{1}\subseteq
\mathcal{B}_{+},$ we are entitled to refer to the usual quasi-interior of
any $E\in \mathcal{B}_{+},$ denoted below by $\tilde{E},$ as in Cor. 2$%
^{\prime };$ we also write $\tilde{B}_{\varepsilon }^{E}$ for $%
B_{\varepsilon }^{\tilde{E}}(\eta ).$

\bigskip

\noindent \textbf{Theorem 7B}.\textit{\ For }$G$\textit{\ a locally compact
group with left Haar measure }$\eta ,$\textit{\ if:}

\noindent (i)\textit{\ }$\tau $ \textit{is both a left- and a
right-invariant refinement topology with }$\tau _{1}\subseteq \mathcal{B}%
_{+} $\textit{\ and with the left Nikodym property (preservation of category
under left shifts),}

\noindent (ii)\textit{\ for every non-empty }$E\in \tau $\textit{\ the
pseudo-norm}%
\[
g\mapsto ||g||_{\tilde{E}}:=\eta (g\tilde{E}\triangle \tilde{E})\qquad (g\in
G)
\]%
\textit{is continuous under }$\tau $\textit{\ at }$g=1_{G}$

\noindent -- \textit{then both the families }$\{AA^{-1}:$\textit{\ }$A\in
\tau _{1}$\} \textit{and }$\{\tilde{B}_{\varepsilon }^{E}:$\textit{\ }$%
\emptyset \neq E\in \tau $ \textit{and }$0<\varepsilon \leq 2\eta (\tilde{E}%
)\}$ \textit{generate neighbourhoods of the identity under which }$G$
\textit{is a topological group. Moreover, the pseudo-norms}%
\[
\{||.||_{\tilde{E}}:\emptyset \neq E\in \tau \}
\]%
\textit{are downward directed by refinement; indeed, for }$\emptyset \neq
E,F\in \tau $\textit{\ and} $\varepsilon <2\min \{\eta (\tilde{E}),\eta (%
\tilde{F})\},$\textit{\ there is }$H\in \tau _{1}$ \textit{such that for }$%
0<\delta <2\eta (\tilde{H})$\textit{\ }%
\[
\tilde{B}_{\delta }^{H}\subseteq \tilde{B}_{\varepsilon }^{E}\cap \tilde{B}%
_{\varepsilon }^{F}.
\]

\noindent \textbf{Proof.} In place of the inclusion of Lemma 2 we note a
result stronger than that valid for $\tilde{E}$ (i.e. inclusion only in $%
\tilde{E}\tilde{E}^{-1}$): since meagreness is translation-invariant (the
`Nikodym property' of [BinO5]), $(xE)^{\widetilde{}}=x\tilde{E}$ for
non-meagre Baire $E,$ so $x\tilde{E}\cap \tilde{E}\neq \emptyset $ implies $%
xE\cap E\neq \emptyset $, and so again%
\[
\tilde{B}_{\varepsilon }^{E}=B_{\varepsilon }^{\tilde{E}}\subseteq EE^{-1};
\]%
here again in Lemma 2 the bound $\eta ^{\ast }(E)$ in the restriction
governing inclusion may be replaced by $2\eta (E)$. The proof of Theorem 7
may now be followed verbatim, but for the replacement of $\mathcal{P}(G)$ by
$\{\eta \}$, using the stronger inclusion just observed, and $B_{\varepsilon
}^{\cdot }(\eta )$ by $\tilde{B}_{\varepsilon }^{\cdot }.$ $\square $

\bigskip

\noindent \textbf{Remark. }The last result follows more directly from Th. 7M
in a context where there exists on $G$ a \textit{Marczewski measure} (see
[TomW, Ch. 13, cf. Ch. 11]), i.e. a finitely additive invariant measure on $%
\mathcal{B}$ vanishing on bounded members of $\mathcal{B}_{0}$; this
includes $\mathbb{R}$, $\mathbb{R}^{2},\mathbb{S}^{1}$, albeit under AC
[TomW, Cor. 13.3]; cf. [Myc], but not $\mathbb{R}^{d}$ for $d\geq 3$ [DohF].

\bigskip

With the groundwork of Part I on translation-continuity for \textit{compacts}
completed, we close by establishing the promised dichotomy associated with
the map%
\[
x\mapsto ||x||_{E}^{\mu }=\mu (xE\triangle E),
\]%
for measurable $E:$ Theorem FN (\S 1) creates a duality between the
vanishing of the $F$-based pseudo-norm and a \textit{dichotomy} for $x$%
-translates of $E^{-1}$ in relation to $F$ according as $x\in E$ or $x\notin
E,$ which are thus unable in each case to distinguish between the points of $%
F.$ Below we write $\forall ^{\mu }$ for the generalized quantifier
\textquotedblleft for $\mu $-a.a.' (cf. [Kec, 8.J]).

\bigskip

\noindent \textbf{Theorem 8 (Almost Inclusion-Exclusion). }\textit{For }$G$%
\textit{\ a Polish group } $\mu\in \mathcal{P}(G)$\textit{\ and non-null }$%
\mu $\textit{-measurable }$E,F$\textit{, the vanishing }$\mu $\textit{-a.e.
on }$F$ \textit{of the }$E$\textit{-norm under }$\mu :$%
\[
||x||_{F}^{\mu }=\mu (xE\triangle E)=0\qquad (x\in F),
\]%
\textit{is equivalent to the following Almost Inclusion-Exclusion for
translates of }$E^{-1}$\textit{:}

\noindent (i) Inclusion: $F$\textit{\ is }$\mu $\textit{-almost covered by }$%
\mu $\textit{-almost every translate }$xE^{-1}$\textit{for }$x\in E:$%
\[
\mu (F\backslash xE^{-1})=0\qquad (\forall ^{\mu }x\in E),
\]%
\noindent (ii) Exclusion:\textit{\ }$F$\textit{\ is }$\mu $\textit{-almost
disjoint from }$\mu $\textit{-almost every} \textit{translate }$xE^{-1}$%
\textit{for }$x\notin E:$%
\[
\mu (F\cap xE^{-1})=0\qquad (\forall ^{\mu }x\notin E).
\]%
\noindent \textbf{Proof. }By Theorem FN (\S 1), applied to the set $H$ of
Prop. 3, i.e.%
\[
H:=\bigcup\nolimits_{x\in F}\{x\}\times (xE\triangle E),
\]%
$H$ has vertical sections $H_{x}$ almost all $\mu $-null iff $\mu $-almost
all of its horizontal sections $H^{y}$ are $\mu $-null. But, since $y\in xE$
iff $x\in yE^{-1}$, $H^{y}=F\backslash yE^{-1}$ for $y\in E$ and $%
H^{y}:=F\cap yE^{-1}$ for $y\in G\backslash E.$ $\square $

\bigskip

\noindent \textbf{Remark. }If the inclusion side of the dichotomy of Th. 8
holds for all $x\in E,$ then $F\subseteq EE^{-1}.$ The converse direction
may fail: consider $E=(1,2)\subseteq \mathbb{R}$ and $F=(-1,1),$ so that $%
E-E=F,$ but no translate of $-E$ may cover $F.$

\section{Other interior-point properties}

\subsection{\textbf{The Steinhaus property }$AA^{-1}$\textbf{\ versus the
Steinhaus property }$AB^{-1}$}

We clarify the relation between two versions of the Steinhaus interior
points property: the simple (sometimes called `classical') version
concerning sets $AA^{-1}$ and the composite, more embracing one, concerning
sets $AB^{-1}$, for sets from a given family $\mathcal{H}$. The latter is
connected to a strong form of metric transitivity: Kominek [Kom] shows, for
a general separable Baire topological group $G$ equipped with an inner
regular measure $\mu $ defined on some $\sigma $-algebra $\mathcal{M}$, that
$AB^{-1}$ has non-empty interior for all $A,B\in \mathcal{M}_{+}$ iff for
each countable dense set $D$ and each $E\in \mathcal{M}_{+}$ the set $%
X\backslash DE\in \mathcal{M}_{0}$; this is recalled in Theorem K\ below.
The composite property is thus related to the Smital property, for which see
[BarFN]. Care is required when moving to the alternative property for $AB,$
since the family $\mathcal{H}$ need not be preserved under inversion.

In general the simple property does not imply the composite:\ Mato\u{u}skov%
\'{a} and Zelen\'{y} [MatZ] show that in any non-locally compact abelian
Polish group there are closed non-(left) Haar null sets $A,B$ such that $A+B$
has empty interior. Recently, Jab\l o\'{n}ska [Jab2] has shown that likewise
in any non-locally compact abelian Polish group there are closed non-Haar
meager sets $A,B$ such that $A+B$ has empty interior. Bartoszewicz and M.
and T. Filipczak [BarFF, Ths. 1, 4] analyze the Bernoulli product measure on
$\{0,1\}^{\mathbb{N}}$ with $p$ the probability of the digit $1;$ see \S %
9.15. The product space may be regarded as comprising canonical binary digit
expansions of the additive reals modulo 1 (in which case the measure is not
invariant). Here the (Borel) set $A$ of binary expansions with asymptotic
frequency $p$ of the digit `1' has $[0,1)$ as its difference set iff $\frac{1%
}{4}\leq p\leq \frac{3}{4};$ however $A+A$ has empty interior unless $p=%
\frac{1}{2}$ (the base $2$ simple-normal-numbers case).

Below we identify some conditions on a family of sets $A$ with the simple $%
AA^{-1}$ property which do imply the $AB^{-1}$ property. What follows is a
generalization to a group context of relevant observations from [BinO9] from
the classical context of $\mathbb{R}$.

The motivation for the definition below is that its subject, the space $H,\ $%
is a subgroup of a topological group $G$ from which it inherits a
(necessarily) translation-invariant (i.e. either-sidedly) topology $\tau .$
Various notions of `density at a point' give rise to `density topologies'
[BinO5], which are translation-invariant since they may be obtained via
translation to a fixed reference point: early examples, which originate in
spirit with Denjoy as interpreted by Haupt and Pauc [HauP], were studied
intensively in [GofW], [GofNN], soon followed by [Mar] and [Mue]; more
recent examples include [FilW] and others investigated by the Wilczy\'{n}ski
school, cf. [Wil].

Proposition 12 below embraces as an immediate corollary the case $H=G$ with $%
G\ $locally compact and $\sigma $ the Haar density topology (see [BinO8]).
Proposition 13 proves that Proposition 12 applies also to the ideal topology
(in the sense of [LukMZ]) generated from the ideal of Haar null sets of an
abelian Polish group.

We recall that a group $H$ carries a \textit{left semi-topological}
structure $\tau $ if the topology $\tau $ is left invariant [ArhT] ($hU\in
\tau $ iff $U\in \tau $); the structure is \textit{semi-topological} if it
is also right invariant, i.e. briefly: $\tau $ is translation invariant. $H$
is a \textit{quasi-topological group} under $\tau $ if $\tau $ is both left
and right invariant and inversion is $\tau $-continuous.

\bigskip

\noindent \textbf{Definition.} For $H$ a group with a translation-invariant
topology $\tau ,$ call a topology $\sigma \supseteq \tau $ a \textit{%
Steinhaus refinement }if:

\noindent i) \textrm{int}$_{\tau }$($AA^{-1})\neq \emptyset $ for each
non-empty $A\in \sigma ,$ and

\noindent ii) $\sigma $ is involutive-translation invariant: $hA^{-1}\in
\sigma $ for all $A\in \sigma $ and all $h\in H.$

\bigskip

Property (ii) above (called simply `invariance' in [BarFN]) calls apparently
for only left invariance, but in fact, via double inversion, delivers
translation invariance, since $Uh=(h^{-1}U^{-1})^{-1}$; then $H$ under $%
\sigma $ is a semi-topological group with a continuous inverse, so a \textit{%
quasi-topological group}.

\bigskip

\noindent \textbf{Proposition 12.} \textit{If }$\tau $\textit{\ is
translation-invariant, and }$\sigma \supseteq \tau $\textit{\ is a Steinhaus
refinement topology, then }\textrm{int}$_{\tau }$($AB^{-1})\neq \emptyset $%
\textit{\ for non-empty }$A,B\in \sigma .$ \textit{In particular, as }$%
\sigma $\textit{\ is preserved under inversion, also }\textrm{int}$_{\tau }$(%
$AB)\neq \emptyset $\textit{\ for }$A,B\in \sigma $.\textit{\ }

\bigskip

\noindent \textbf{Proof.} Suppose $A,B\in \sigma $ are non-empty; as $%
B^{-1}\in \sigma ,$ choose $a\in A$ and $b\in B;$ then by (ii)%
\[
1_{H}\in C:=a^{-1}A\cap b^{-1}B^{-1}\in \sigma .
\]%
By (i), for some non-empty $W\in \tau ,$%
\[
W\subseteq CC^{-1}=(a^{-1}A\cap b^{-1}B)\cdot (A^{-1}a\cap B^{-1}b)\subseteq
(a^{-1}A)\cdot (B^{-1}b).
\]%
As $\tau $ is translation invariant, $aWb^{-1}\in \tau $ and
\[
aWb^{-1}\subseteq AB^{-1},
\]%
the latter since for each $w\in W$ there are $x\in A,y\in B^{-1}$ with
\[
w=a^{-1}x.yb:\qquad awb^{-1}=xy\in AB^{-1}.
\]%
So \textrm{int}$_{\tau }$($AB^{-1})\neq \emptyset .$ $\square $

\bigskip

\noindent \textbf{Corollary 7.} \textit{In a locally compact group the Haar
density topology is a Steinhaus refinement.}

\bigskip

\noindent \textbf{Proof. }Property (i) follows from Weil's theorem since
density-open sets are non-null measurable; left translation invariance in
(ii) follows from left invariance of Haar measure, while involutive
invariance holds, as any measurable set of positive Haar measure has
non-null inverse ([HewR, 15.14], cf. \S 4 Lemma H). $\square $

\bigskip

A weaker version, inspired by metric transitivity, comes from applying the
following concept.

\bigskip

\noindent \textbf{Definition.} Say $H$ \textit{acts transitively} on $%
\mathcal{H}$ for each $A,B\in \mathcal{H}$ if there is $h\in H$ with $A\cap
hB\in \mathcal{H}$.

\bigskip

Thus a locally compact topological group acts transitively on the non-null
Haar measurable sets (in fact, either-sidedly); this follows from Fubini's
Theorem [Hal, 36C], via the average theorem [Hal, 59.F]:%
\[
\int_{G}|g^{-1}A\cap B|dg=|A|\cdot |B^{-1}|\qquad (A,B\in \mathcal{M}),
\]%
($g=ab^{-1}$ iff $g^{-1}a=b)$ -- cf. [TomW, \S 11.3 after Th. 11.17].

[MatZ] show that in any non-locally compact abelian Polish group $G$ there
exist two non-Haar null sets, $A,B\notin \mathcal{HN}$, such that $A\cap
hB\in \mathcal{HN}$ for all $h;$ that is, $G$ here does \textit{not} act
transitively on the non-Haar null sets.

\bigskip

\noindent \textbf{Definition.} In a quasi-topological group $(H,\tau )$ say
that a proper $\sigma $-ideal $\mathcal{H}$ has the \textit{Simple Steinhaus
Property }$AA^{-1}$ if $AA^{-1}$ has interior points for universally
measurable subsets $A\notin \mathcal{H}.$ This follows [BarFN].

\bigskip

\noindent \textbf{Proposition 12}$^{\prime }$ (cf. [Kha, Th. 1])\textbf{.}
\textit{In a group }$(H,\tau )$ \textit{with }$\tau $\textit{\
translation-invariant, if }$H$ \textit{acts transitively on a family of
subsets }$\mathcal{H}$\textit{\ with the simple Steinhaus property, then }$%
\mathcal{H}$ \textit{has the composite Steinhaus property: }\textrm{int}$%
_{\tau }$($AB^{-1})\neq \emptyset $\textit{\ for }$A,B\in \mathcal{H}$.%
\textit{\ Furthermore, if }$\mathcal{H}$\textit{\ is preserved under
inversion, then also }\textrm{int}$_{\tau }$($AB)\neq \emptyset $\textit{\
for }$A,B\in \mathcal{H}$.

\bigskip

\noindent \textbf{Proof. }For $A,B\in \mathcal{H}$ choose $h$ with $C:=A\cap
hB\in \mathcal{H}$; then%
\[
CC^{-1}h=(A\cap hB)(A^{-1}\cap B^{-1}h^{-1})\subseteq AB^{-1}.\qquad \square
\]

\bigskip

\noindent \textbf{Proposition 13.} \textit{If }$(H,\tau )$\textit{\ is a
quasi-topological group (i.e. }$\tau $\textit{\ is invariant with continuous
inversion) carrying a left invariant }$\sigma $\textit{-ideal }$\mathcal{H}$
\textit{with the Steinhaus property and }$\tau \cap \mathcal{H}=\{\emptyset
\}$\textit{, then the ideal-topology }$\sigma $\textit{\ with basis }%
\[
\mathcal{B}:=\{U\backslash N:U\in \tau ,N\in \mathcal{H}\}
\]%
\textit{is a Steinhaus refinement of }$\tau $\textit{.}

\textit{In particular, for }$(H,\tau )$\textit{\ an abelian Polish group,
the ideal topology generated by the }$\sigma $\textit{-ideal of Haar null
subsets is a Steinhaus refinement.}

\bigskip

\noindent \textbf{Proof. }If $U,V\in \mathcal{B}$ and $w\in U\cap V,$ choose
$M,N\in \mathcal{H}$ and $W_{M},W_{N}\in \tau $ such that $x\in
(W_{M}\backslash M)\subseteq U$ and $x\in (W_{N}\backslash N)\subseteq V.$
Then as $M\cup N\in \mathcal{H}$,%
\[
x\in (W_{M}\cap W_{N})\backslash (M\cup N)\in \mathcal{B}.
\]%
So $\mathcal{B}$ generates a topology $\sigma $ refining $\tau .$ With the
same notation, $hU=hW_{M}\backslash hM\in \sigma ,$ as $hM\in H,$ and $%
U^{-1}=W_{M}^{-1}\backslash M^{-1}$. Finally, $UU^{-1}$ has non-empty $\tau $%
-interior, as $U\notin \mathcal{H}$ and is non-empty.

As for the final assertion concerned with an abelian Polish group context,
note that if $N$ is Haar null ($N\in \mathcal{HN}$), then $\mu (hN)=0$ for
some $\mu \in \mathcal{P}(G)$ and all $h\in H,$ so $hN\in \mathcal{HN}$ for
all $h\in H.$ Furthermore, if $A\notin \mathcal{HN}$ then $A^{-1}\notin
\mathcal{HN}$, otherwise $\mu (hA^{-1})=0$ for some $\mu \in \mathcal{P}(G)$
and all $h\in H;$ then, taking $\tilde{\mu}(B)=\mu (B^{-1})$ for Borel $B,$
we have $\tilde{\mu}(A)=0$ and $\tilde{\mu}(hA)=\mu (A^{-1}h^{-1})=0$ for
all $h\in H,$ a contradiction. $\square $

\bigskip

\noindent \textbf{Remark.} A left Haar null set need not be right Haar null:
for one example see [ShiT], and for more general non-coincidence see Solecki
[Sol1, Cor. 6]. So the argument in Prop. B does not extend to the family of
left Haar null sets $\mathcal{LHN}$ of a \textit{non-commutative} Polish
group. Indeed, Solecki [Sol2, Th. 1.4] shows in the context of a countable
product of countable groups that the simpler Steinhaus property holds for $%
\mathcal{HN}$ iff $\mathcal{HN}=\mathcal{LHN}$.

\bigskip

Next, a result from [Kom]. Recall that $\mu $ is \textit{quasi-invariant} if
$\mu $-nullity is translation invariant. The transitivity assumption (of
co-nullity) is motivated by \textit{Sm\'{\i}tal's lemma}, which refers to a
countable dense set -- see [KucS].

\bigskip

\noindent \textbf{Theorem K. (}[Kom, Th. 5]). \textit{If }$\mu \in P(G)$%
\textit{\ is quasi-invariant and there exists a countable subset }$%
H\subseteq G\ $\textit{with }$HM$\textit{\ co-null for all }$M\in M_{+}(\mu
),$\textit{\ then int}$(AB^{-1})\neq \emptyset $\textit{\ for all }$A,B\in
\mathcal{M}_{+}(\mu ).$

\bigskip

\noindent \textbf{Proof. }By regularity, we may assume $A,B\in \mathcal{M}%
_{+}(\mu )$ are compact, so $AB^{-1}$ is compact. Fix $g\in G$; then by
quasi-invariance $\mu (gB)>0.$ So by the transitivity assumption, both $%
G\backslash HgB$ and $G\backslash HA$ are null, and so $HA\cap HgB\neq
\emptyset .$ Say $h_{1}a=h_{2}gb,$ for some $a\in A,b\in B,h_{1},h_{2}\in H;$
then $g=h_{2}^{-1}h_{1}ab^{-1}.$ As $g$ was arbitrary,%
\[
G=\dbigcup\nolimits_{h\in H}h_{2}^{-1}h_{1}AB^{-1}.
\]%
By Baire's Theorem, as $H$ is countable, \textrm{int}$(AB^{-1})\neq
\emptyset .$ $\square $

\subsection{Borell's interior-point property}

For completeness of this overview of the Steinhaus-Weil interior-point
property, we offer in brief here the context and statement of a (by now)
classical Steinhaus-like result in probability theory; this differs in that
the Polish group now specializes to an infinite-dimensional topological
vector space and the reference measure is Gaussian, so no longer invariant.
We refer to the companion paper [BinO8] for further details and background
literature, and to our generalizations to Polish groups and to other
reference measures.

For $X$ a locally convex topological vector space, $\gamma $ a probability
measure on the $\sigma $-algebra of the cylinder sets generated by the dual
space $X^{\ast }$ (equivalently, for $X$ separable Fr\'{e}chet, e.g.
separable Banach, the Borel sets), with $X^{\ast }\subseteq L^{2}(\gamma ):$
then $\gamma $ is called \textit{Gaussian} on $X$ (gamma for Gaussian,
following [Bog1]) iff $\gamma \circ \ell ^{-1}$ defined by
\[
\gamma \circ \ell ^{-1}(B)=\gamma (\ell ^{-1}(B))\qquad (\text{Borel }%
B\subseteq \mathbb{R})
\]%
is Gaussian (normal) on $\mathbb{R}$ for every $\ell \in X^{\ast }\subseteq
L^{2}(\gamma )$. For a monograph treatment of Gaussianity in a Hilbert-space
setting, see Janson [Jan]. Write $\gamma _{h}(K):=\gamma (K+h)$ for the
translate by $h.$ \textit{Relative} \textit{quasi-invariance} of $\gamma
_{h} $ and $\gamma ,$ that for all compact $K$
\[
\gamma _{h}(K)>0\text{ iff }\gamma (K)>0,
\]%
holds relative to a set of vectors $h\in X$ (the\textit{\ admissible
directions}) forming a vector subspace known as the\textit{\ Cameron-Martin
space, }$H(\gamma )$. Then $\gamma _{h}$ and $\gamma $ are equivalent, $%
\gamma \sim \gamma _{h},$ iff $h\in H(\gamma ).$ Indeed, if $\gamma \sim
\gamma _{h}$ fails, then the two measures are mutually singular, $\gamma
_{h}\bot \gamma $ (the Hajek-Feldman Theorem -- cf. [Bog1, Th. 2.4.5,
2.7.2]).

Continuing with the assumption above on $X^{\ast },$ as $X\subseteq X^{\ast
\ast }\subseteq L^{2}(\gamma ),$ one can equip $H=H(\gamma )$ with a norm
derived from that on $L^{2}(\gamma ).$ In brief, this is done with reference
to a natural covariance under $\gamma $ obtained by regarding $f\in X^{\ast
} $ as a random variable and working with its zero-mean version $f-\gamma
(f);$ then, for $h\in H,$ $\delta _{h}^{\gamma },$ the (shifted) evaluation
map defined by $\delta _{h}^{\gamma }(f):=f(h)-\gamma (f)$ for $f\in X^{\ast
},$ is represented as $\langle f-\gamma (f),\hat{h}\rangle _{L^{2}(\gamma )}$
for some $\hat{h}\in L^{2}(\gamma ).$ (Here for $\gamma $ symmetric $\gamma
(f)=0,$ so $\delta _{h}^{\gamma }=\delta _{h}$ is the Dirac measure at $h$.)
This is followed by identifying $h$ with $\hat{h}$ (for $h\in H$), and $%
|h|_{H}:=||\hat{h}||_{L^{2}(\gamma )}$ is a norm on $H$ arising from the
inner product%
\[
(h,k)_{H}:=\int_{X}\hat{h}(x)\hat{k}(x)d\gamma (x).
\]%
Formally, the construction first requires an extension of the domain of $%
\delta _{h}^{\gamma }$ to $X_{\gamma }^{\ast },$ the closed span of $%
\{x^{\ast }-\gamma (x^{\ast }):x^{\ast }\in X^{\ast }\}$ in $L^{2}(\gamma ),$
a Hilbert subspace in which to apply the Riesz Representation Theorem.

We may now state the Steinhaus-like property due, essentially in this form,
to Christer Borell. ([LeP, Prop. 1] offers a weaker,
`one-dimensional-section' form, which would now be called, as in e.g. [Brz],
a `radial' form.)

\bigskip

\noindent \textbf{Theorem B (Borell's Interior-point Theorem}, [Bor, Cor.
4.1] -- see [Bog1, p. 64]). \textit{For }$\gamma $\textit{\ a Gaussian
measure on a locally convex topological space }$X$\textit{\ with }$X^{\ast
}\subseteq L^{2}(\gamma )$\textit{, and }$A$\textit{\ any non-null }$\gamma $%
\textit{-measurable subset }$A$\textit{\ of }$X,$\textit{\ the difference
set }$A-A$\textit{\ contains a }$|.|_{H}$\textit{-open nhd (neighbourhood)
of }$0$\textit{\ in the Cameron Martin space }$H=H(\gamma ),$\textit{\ i.e. }%
$(A-A)\cap H$\textit{\ contains a }$H$\textit{-open nhd of }$0.$

\bigskip

This flows from the continuity in $h$ of the density of $\gamma _{h}$ wrt $%
\gamma $ ([Bog1, Cor. 2.4.3]), as given in the \textit{%
Cameron-Martin-Girsanov formula}:%
\begin{equation}
\exp \left( \hat{h}(x)-\frac{1}{2}||\hat{h}||_{L^{2}(\gamma )}^{2}\right)
\tag{$CM$}
\end{equation}%
(where $\hat{h}$ `Riesz-represents' $h,$ i.e. $x^{\ast }(h)=\langle x^{\ast
},\hat{h}\rangle ,$ for $x^{\ast }\in X^{\ast },$ as above). Thus here a
modified Steinhaus Theorem holds: the \textit{relative-interior-point theorem%
}.

\section{Complements}

\noindent \textbf{1}.\textbf{\ }\textit{Non-separability. }Simmons
establishes his theorems without separability by using the Kakutani-Kodaira
`separable quotient' approximation theorem for metric groups [MonZ, I.2.6],
cf. [HewR, Th. 8.7] and the recent generalization [Hu]. (Starting with a
fixed compact support permits a reduction to the case of a compactly
generated group, which is followed by working in a \textit{separable}
quotient group.)

The links between the Effros theorem (\S 8.7 below), the Baire theorem and
the Steinhaus-Weil theorem are pursued at length in [Ost2]. There, any
separability assumption is avoided. Instead \textit{sequential} methods are
used, for example shift-compactness arguments (\S 8.4 below).

\noindent \textbf{2}. \textit{Cameron-Martin theory. }In the context of
Gaussian measures in probability theory, Bogachev [Bog1] gives a thorough
treatment of Cameron-Martin theory (on translation of Wiener and other
Gaussian measures; cf. abstract Wiener spaces, reproducing-kernel Hilbert
spaces, etc.). Here results of Steinhaus-Weil type occur, in the context of
topological vector spaces [Bog1, p.64], as we have seen in \S 7.2, albeit
only relative to a specific \textit{subspace}. In the companion paper to
this, [BinO8], we develop an analogous theory in the context of topological
groups.

\noindent \textbf{3.}\textit{\ Smallness. }When the (classical)
Cameron-Martin subspace of a topological vector space $X$, alluded to in \S %
8.1 (and mentioned below Th. 5 in \S 6), is a separable Hilbert space (as
when the Gaussian measure $\gamma $ is a Radon measure, see [Bog1, 3.2.7],
cf. [Bog1, p. 62]), so complete under its norm and hence analytic [Rog],
being $\gamma $-null it has empty interior. So by a result of Dodos [Dod3,
Cor. 9], it is (generically) left-Haar null -- left-Haar null for quasi all $%
\mu \in \mathcal{P}(X),$ in the sense of the L\'{e}vy metric on $\mathcal{P}%
(X).$

\noindent \textbf{4}. \textit{Inclusion-Exclusion dichotomy. }Our Weil-like
analysis in Part III focuses on inclusions amongst sets of the form $%
EE^{-1}, $ the exception being the Inclusion-Exclusion of a set $F$ by an $E$%
- or non-$E$-translate of $E^{-1}$ in Theorem 8 (a dichotomy as between $E$
and its complement). This places most of our study on one side of a related
inclusion-exclusion dichotomy, for subsets $H,B$ of a group $G$ one has
either inclusion , or `near-disjointness':%
\[
HH^{-1}\subseteq BB^{-1},\qquad \text{or}\qquad HH^{-1}\cap
BB^{-1}=\{1_{G}\}.
\]%
Inclusion may be equivalently re-phrased to the meeting of distinct pairs of
$H^{-1}$-translates of $B:$%
\begin{equation}
kB\cap k^{\prime }B\neq \emptyset \qquad (k,k^{\prime }\in H^{-1}),  \tag{I}
\end{equation}%
whereas exclusion to their disjointness:%
\begin{equation}
kB\cap k^{\prime }B=\emptyset \qquad (\text{distinct }k,k^{\prime }\in
H^{-1}).  \tag{E}
\end{equation}

The duality of the relation of (E) to the results in Th. 8 is clarified by
observing that $\mu (F\cap xE^{-1})=0,$ for a.a. $x\in C,$ is equivalent to $%
\mu (C\cap yE)=0,$ for a.a. $y\in F$. Indeed,
\[
0=\int \int 1_{C}(x)1_{F}(y)1_{xE^{-1}}(y)d(\mu \times \mu )=\int \int
1_{F}(y)1_{C}(x)1_{yE}(x)d(\mu \times \mu ).
\]

The condition (E) gives rise to $\mathcal{I}_{0}$, the $\sigma $-ideal
introduced in Balcerzak et al. [BalRS], generated by Borel sets $B$ having
perfectly many disjoint translates, as in (E) above with $H^{-1}$ a perfect
compact set (i.e. compact and dense-in-itself); continuum-many disjoint
translates of a compactum also emerge in a theorem of Ulam concerning a
non-locally compact Polish group: see [Oxt1, Th. 1]. Such \textit{perfect
exclusions} offer a combinatorial tool, akin to \textit{shift-compactness}
(as in Th. 3, subsequence embedding under translation of a null sequence
into a non-negligible set -- cf. [BinO2,3] [MilO]), and play a key role in
the context of groups with \textit{ample generics}; see for instance the
small-index property of [HodHLS].

Solecki [Sol3] proves a `Fubini for negligibles'-type theorem (cf. Theorem
FN\ in \S 1 above): the non-negligible vertical sections (relative to a
uniformly Steinhaus ideal) of a planar $\mathcal{I}_{0}$-negligible set form
a horizontal $\mathcal{I}_{0}$-negligible set. The ideal $\mathcal{I}_{0}$
is of interest, as it violates the countable (anti)-chain condition, [BalRS].

\noindent \textbf{5}. \textit{Steinhaus-Weil property of a Borel measure. }%
In a locally compact group $G,$ the family $\mathcal{M}_{+}(\mu )$ of finite
non-null measurable sets of a Borel measure $\mu $ on $G$ fails to have the
Steinhaus-Weil property iff there are a null sequence $z_{n}\rightarrow
1_{G} $ and a non-null compact set $K$ with $\lim_{n}\mu (t_{n}K)=0.$
Equivalently, this is so iff the measure $\mu $ is not absolutely continuous
with respect to Haar measure: these observations motivated Th. SM.

Alternative characterizations of absolute continuity and singularity of
Radon measures on a locally compact group include the following.

First for absolute continuity: \newline
(i) factorization under convolution: $\mu =f\ast \nu $ with $f\in
L^{1}(G,\eta )$ [LiuR];\newline
(ii) \textit{continuity} of the map $g\mapsto \mu _{g}$, using the total
variation norm on the space of measures [LiuR];\newline
(iii) \textit{separability of the orbit} $\{\mu _{g}:g\in G\}$ in the space
of measures, under the same norm as in (ii) [Lar], [LiuRW], [Tam];\newline
(iv) `Glicksberg separability': for a fixed sequence $\{g_{n}\},$ each range
set $\mathcal{R}(K):=\{\mu _{g}(K):g\in K\},$ for $\mu $-null compacts sets $%
K,$ contains $\{\mu (g_{n}K):n\in \mathbb{N}\}$ as a dense subset [Gli];%
\newline
(v) \textit{outer-equiregularity} of sequences of translates of the $\mu $%
-null compact sets $K$: for every $\varepsilon >0$, compact $\mu $-null $K,$
and sequence $\{g_{n}\}$ there is an open $W\supseteq K$ with $\mu
(g_{n}(W\backslash K))<\varepsilon $ for all $n$ [Gli]

Then for singularity: `self-singularity', e.g. $\mu \perp \mu _{g}$ for $%
\eta $-almost all $g,$ pursued in [LiuR], [LiuRW], [Pro], cf. [BinO8].

\noindent \textbf{6}. \textit{Regular open sets. }Recall that $U\ $is
regular open if $U=$\textrm{int}$($\textrm{cl}$U),$ and that \textrm{int}$($%
\textrm{cl}$U)$ is itself regular open; for background see e.g. [GivH, Ch.
10]. For $\mathcal{D=D}_{\mathcal{B}}$ the Baire-density topology of a
normed topological group, let $\mathcal{D}_{\mathcal{B}}^{RO}$ denote the
regular open sets. For $D\in \mathcal{D}_{\mathcal{B}}^{RO}$, put%
\[
N_{D}:=\{t\in G:tD\cap D\neq \emptyset \}=DD^{-1},\qquad \mathcal{N}%
_{1}:=\{N_{D}:1_{G}\in D\in \mathcal{D}_{RO}\};
\]%
then $\mathcal{N}_{1}$ is a base at $1_{G}$ (since $1_{G}\in C\in \mathcal{D}%
_{RO}$ and $1_{G}\in D\in \mathcal{D}_{RO}$ yield $1_{G}\in C\cap D\in
\mathcal{D}_{RO})$ comprising $\mathcal{T}$-neighbourhoods that are $%
\mathcal{D}_{\mathcal{B}}$-open (since $DD^{-1}=\bigcup \{Dd^{-1}:d\in D\})$%
. We raise the (metrizability)\ question, by analogy with the Weil topology
of a measurable group (see \S 5 and \S 8.4 above): with $\mathcal{D}_{%
\mathcal{B}}$ above replaced by a general density topology $\mathcal{D}$ on
a group $G,$ when is the topology generated by $\mathcal{N}_{1}$ on $G$ a
norm topology? Some indications of an answer may be found in [ArhT, \S 3.3].
We note the following answer in the context of Theorem 7B; compare Struble's
Theorem [Str2], or [DieS, Ch. 8]. If there exists a separating sequence $%
D_{n},$ i.e. such that for each $g\neq 1_{G}$ there is $n$ with $%
||g||_{D_{n}}=1,$ then%
\[
||g||:=\sum\nolimits_{n}2^{-n}||g||_{D_{n}}
\]%
is a norm, since it is separating and, by the Nikodym property, $(D\cap
g^{-1}D)=g^{-1}(gD\cap D)\in \mathcal{B}_{0}$.

\noindent \textbf{7}. \textit{The Effros Theorem }asserts that a transitive
continuous action of a Polish group $G$ on a space $X$ of second category in
itself is necessarily `open', or more accurately is microtransitive (the
(continuous) evaluation map $e_{x}:g\mapsto g(x)$ takes open neighbourhoods $%
E$ of $1_{G}$ to open neighbourhoods that are the orbit sets $E(x)$ of $x$).
It emerges that this assertion is very close to the shift-compactness
property: see [Ost2]. The Effros Theorem reduces to the Open Mapping Theorem
when $G,X\ $are Banach spaces regarded as additive groups, and $G$ acts on $%
X\ $by a linear surjection $L:G\rightarrow X$ via $g(x)=L(g)+x.$ Indeed,
here $e_{0}(E)=L(E).$ For a neat proof, choose an open neighbourhood $U$ of $%
0$ in $G$ with $E\supseteq U-U$; then $L(U)$ is Baire (being analytic) and
non-meagre (since $\{L(nU):n\in N\}$ covers $X),$ and so $L(U)-L(U)\subseteq
L(E)$ is an open neighbourhood of $0$ in $X.$

\noindent \textbf{8}. \textit{Haar null and left Haar null. }The two
families, which are both left and right translation-invariant (cf. [Sol2, p.
696] -- if $\mu \in \mathcal{P}(G)$ witnesses that $E$ is left Haar null,
then $\mu _{g^{-1}}$ witnesses that $Eg$ is left Haar null), coincide in
Polish abelian groups, and in locally compact second countable groups (where
they also coincide with the sets of Haar measure zero -- by an application
of the Fubini theorem). The former family, however, is in general smaller;
indeed, non-Haar null sets need not have the Steinhaus-Weil property -- see
[Sol2].

\noindent \textbf{9}. \textit{Beyond local compactness:\ Haar
category-measure duality. }In the absence of Haar measure, the definition of
left Haar null subsets of a topological group $G$ requires $\mathcal{U}(G),$
the universally measurable sets -- by dint of the role of the totality of
(probability) measures on $G$. The natural dual of $\mathcal{U}(G)$ is the
class $\mathcal{U}_{\mathcal{B}}(G)$ of \textit{universally Baire sets},
defined for $G$ with a Baire topology as those sets $B$ whose preimages $%
f^{-1}(B)$ are Baire in any compact Hausdorff space $K$ for any continuous $%
f:K\rightarrow G$. Initially considered in [FenMW] for $G=\mathbb{R}$, these
have attracted continued attention for their role in the investigation of
axioms of determinacy and large cardinals -- see especially [Woo]; cf.
[MarS].

Analogously to the left Haar null sets, define a \textit{left Haar meagre}
set as any set $M$ coverable by a universally Baire set $B$ for which there
are a compact Hausdorff space $K$ and a continuous $f:K\rightarrow G$ with $%
f^{-1}(gB)$ meagre in $K$ for all $g\in G.$ These were introduced, in the
abelian Polish group setting with $K$ metrizable, by Darji [Dar], cf.
[Jab1], and shown there to form a $\sigma $-ideal of meagre sets
(co-extensive with the meagre sets for $G\ $locally compact).

\noindent \textbf{10}. \textit{Metrizability and Christensen's Theorem.} An
analytic topological group is metrizable; so if also it is a Baire space,
then it is a Polish group -- [HofT, Th. 2.3.6].

\noindent \textbf{11}. \textit{Strong Kemperman property: qualitative versus
quantitative measure theory. }We note that property (*) of the Introduction
and Lemma 1 corresponds to the following quantitative property on the line,
stated in terms of Haar (i.e. Lebesgue) measure $\eta (\cdot )=|\cdot |$ for
sets open in the Lebesgue density topology $\mathcal{D}_{\mathcal{L}}$:

\noindent (iv)* strong Kemperman property (see [Kem], [Kuc, Lemma 3.7.2]):
for $0\in U\in \mathcal{D}_{\mathcal{L}}$ there is $\delta >0$ so that for
all $|t|<\delta $%
\[
|U\cap (t+U)|\geq \varepsilon .
\]%
This is connected with continuity of the Haar norm. Indeed, since%
\[
|U\cap (t+U)|=|U|-|U\triangle (t+U)|/2,
\]%
the inequality above is equivalent to
\[
||t||_{U}^{\eta }:=|U\triangle (t+U)|\leq 2(|U|-\varepsilon ).
\]%
The latter holds for any $0<\varepsilon <|U|\ $and for sufficiently small $%
t, $ by the continuity of the norm $||t||_{U}^{\eta }$.

\noindent \textbf{12}. \textit{Proof of Theorem FN of \S 1.}\textbf{\ }For $%
\mu $-null\textbf{\ }$N\subseteq G$ the set $N\times G$ is $\mu \times \nu $%
-null, so (by passing to the complement of the null exceptional set of the
theorem) we may assume w.l.o.g. that the exceptional set of $A$ is empty. By
inner regularity, it suffices to show that $(\mu \times \nu )(K)=0$ for all
compact $K\subseteq A.$

For $K$ compact, denote by $F\ $the (compact) projection of $K$ on the first
axis. Let $\varepsilon >0.$ By Prop. 5, for any $x\in F$ there is an open
neighbourhood $U_{x}$ of $x$ and open $V_{x}$ with $\nu (V_{x})<\varepsilon $
and
\[
K\cap (U_{x}\times G)\subseteq R_{x}:=U_{x}\times V_{x}.
\]%
By compactness of $F,$ there are $U^{j}\times V^{j}$ for $i=1,...,n,$ with $%
U^{j},V^{j}$ open and $\nu (V^{j})<\varepsilon $ such that
\[
F\subseteq \bigcup\nolimits_{j}U^{j}:\qquad K\subseteq
\bigcup\nolimits_{j}U^{j}\times V^{j}.
\]%
To disjoin the sets $U^{j},$ put%
\[
S^{j}:=U^{j}\backslash \bigcup\nolimits_{j<i}U^{j}:\qquad
\bigcup\nolimits_{j}U^{j}=\bigcup\nolimits_{j}S^{j}.
\]%
Then%
\[
F=\bigcup\nolimits_{j}F\cap S^{j}:\qquad K\subseteq
\bigcup\nolimits_{j}S^{j}\times V^{j}.
\]%
So%
\[
\mu (K)\leq \sum\nolimits_{j}(\mu \times \nu )(S^{j}\times
V^{j})=\sum\nolimits_{j}\mu (S^{j})\nu (V^{j})\leq \sum\nolimits_{j}\mu
(S^{j})\cdot \varepsilon =\varepsilon \mu (F).
\]%
As $\varepsilon >0$ was arbitrary, $\mu (K)=0.\qquad \square $

\noindent \textbf{13}. \textit{Proof of Proposition 5 (Almost everywhere
upper semicontinuity)}\textbf{. }As this is much as in Prop. 2, suffice it
to indicate the necessary adjustments. For measurable $E\subseteq G$
measurable and compact $F$ with $\mu (F)>0,$ refer to
\[
H:=\Phi _{E}(F)=\bigcup\nolimits_{x\in F}\{x\}\times xE.
\]%
Instead of the finite union of rectangles $K,$ choose a compact set $%
K=K_{\varepsilon }\subseteq H$ with%
\[
(\mu \times \mu )(H\backslash K)<\varepsilon ^{2}.
\]%
Then obtain a compact set $C=C_{\varepsilon }\subseteq F\backslash
S_{\varepsilon }$ with $\mu (F\backslash C)<\varepsilon $ and%
\[
\mu (H_{x}\backslash K_{x})\leq \varepsilon \qquad (x\in C_{\varepsilon }).
\]%
For $x\in C_{\varepsilon },$ by upper semicontinuity of $t\mapsto \mu
(K_{t}),$ there is a neighbourhood $U_{x}^{\varepsilon }$ of $x$ with%
\[
\mu (K_{y})<\mu (K_{x})+\varepsilon \qquad (y\in U_{x}^{\varepsilon }\cap
C).
\]%
So%
\[
\mu (H_{y})<\mu (K_{y})+\varepsilon <\mu (K_{x})+2\varepsilon <\mu
(H_{x})+3\varepsilon ,
\]%
i.e.%
\[
\mu (H_{y})<\mu (H_{x})+3\varepsilon \qquad (x\in C,y\in U_{x}^{\varepsilon
}\cap C).
\]%
Thereafter the inductive argument of Prop. 4 follows almost verbatim. $%
\square $

\noindent \textbf{14}. \textit{Carlson--Simpson Theorem.} The
shift-compactness property in Th. 3 \S 3 concerned with embeddings appears
to be related to the following van der Waerden-like theorem in Ramsey Theory
(cf. [BinO1, \S 6]). This is due to Carlson and Simpson [CarS] and is a
generalization of the famous Hales-Jewett Theorem; for a proof and related
literature see [DodKT]:

\textit{For every integer }$k\geq 2$\textit{\ and every finite coloring of
the set of all (finite) words over }$[k]:=\{1,...,k\}$\textit{\ there exist
a word }$c$\textit{\ over }$[k]$\textit{\ and a sequence of left-variable
words }$w_{n}(v)$\textit{\ over }$[k]\cup \{v\}$\textit{\ such that}%
\[
\{c\}\cup \{c\cdot w_{1}(a_{1})\cdot ...\cdot
w_{n}(a_{n}):a_{1},...,a_{n}\in \lbrack k]\text{ and }n\in \mathbb{N}\}
\]%
\textit{is monochromatic. }

Here $\cdot $ denotes concatenation and $w(a)$ is obtained by substituting $%
a $ for $v$. A \textit{word} over $S$ is a string of symbols from the
(finite) collection $S$ of symbols (the vocabulary) and a \textit{%
left-variable word} $w(v)$ over $[k]\cup \{v\}$ means a word over $[k]\cup
\{v\}$ beginning with the symbol $v,$ where $v\notin \lbrack k].$

\noindent \textbf{15}. \textit{Infinite Bernoulli convolutions}. The
measures above in [BarFF] are related to infinite Bernoulli convolutions
with parameter (probability) $p\in (0,1)$. These have been studied since the
work of Erd\H{o}s in 1939 and 1940; for background and references, see e.g.
Cooper [Coo], Solomyak [Solom]. By the law of pure types from probability
theory, they are either absolutely continuous (the generic case), or
continuous singular. This second case occurs if $p$ is a PV-number -- an
algebraic number of Pisot-Vijayaraghavan type.

\noindent \textbf{16. }\textit{Quasi-invariance and the Mackey topology of
analytic Borel groups. }We stop to comment on the force of full
quasi-invariance of a measure in connection with a Steinhaus triple $%
(H,G,\mu )$ with $H$ and $G$ completely metrizable. Both groups, being
absolutely Borel, are analytic spaces. So both carry a `standard' Borel
structures with $H$ a Borel substructure of $G.$ Mackey [Mac] investigates
such Borel groups, defining also a (Borel) measure $\mu $ to be `standard'
if it has a Borel support It emerges that every $\sigma $-finite Borel
measure in an analytic Borel space is standard [Mac, Th. 6.1]. Of interest
to us is Mackey's notion of a `measure class' $C_{\mu },$ comprising all
Borel measures $\nu $ with the same null sets as $\mu :$ $\mathcal{M}%
_{0}(\nu )=\mathcal{M}_{0}(\mu ).$ Such a measure class may be closed under
translation, and may be right or left invariant; then the common null sets
are themselves invariant, and so may be viewed as witnessing
quasi-invariance of the measure $\mu .$ Mackey shows that a Borel group with
a one-sided invariant measure class has a both-sidedly invariant measure
class [Mac, Lemma 7.2]; furthermore, if the class is countably generated,
then the class contains a left-invariant and a right-invariant measure [Mac,
Lemma 7.3]. This enables Mackey to improve on Weil's theorem in showing that
an analytic Borel group $G$ with a one-sidedly invariant measure class, in
particular one generated by a quasi-invariant measure, has a unique locally
compact topology both making $G$ a topological group and generating the
given Borel structure.

\section{References}

\noindent \lbrack Amb] {\ W. Ambrose, {Measures on locally compact
topological groups.} \textsl{Trans. Amer. Math. Soc.}\ \textbf{61} (1947),
106-121.}\newline
[ArhT] {A. Arhangelskii, M. Tkachenko, \textsl{Topological groups and
related structures. }World Scientific, 2008.}\newline
[BeeV] {G. Beer, L. Villar, {Borel measures and Hausdorff distance.} \textsl{%
Trans. Amer. Math. Soc. }\textbf{307} (1988), 763--772.}\newline
[BalRS] {M. Balcerzak, A. Ros\l anowski, S. Shelah, {Ideals without ccc.}%
\textsl{\ J. Symbolic Logic} \textbf{63} (1998), 128-148.}\newline
[Bana] W. Banaszczyk, \textsl{Additive subgroups of topological vector
spaces.} Lecture Notes in Mathematics \textbf{1466}. Springer, 1991.\newline
[BarFF] A. Bartoszewicz, M. Filipczak, T. Filipczak, On supports of
probability Bernoulli-like measures. \textsl{J. Math. Anal. Appl.} \textbf{%
462} (2018), 26--35.\newline
[BarFN] A. Bartoszewicz, M. Filipczak, T. Natkaniec, On Sm\'{\i}tal
properties. \textsl{Topology Appl.} \textbf{158} (2011), 2066--2075.\newline
[Bin] {N. H. Bingham, {Finite additivity versus countable additivity.}\
\textsl{Electronic J. History of Probability and Statistics}, \textbf{6}
(2010), 35p.}\newline
[BinGT] {N. H. Bingham, C. M. Goldie and J. L. Teugels, \textsl{Regular
variation}{,} 2nd ed., Cambridge University Press, 1989 (1st ed. 1987).}%
\newline
[BinO1]{\ N. H. Bingham and A. J. Ostaszewski, {Kingman, category and
combinatorics.} \textsl{Probability and mathematical genetics} (Sir John
Kingman Festschrift, ed. N. H. Bingham and C. M. Goldie), 135-168, \textsl{%
London Math. Soc. Lecture Notes in Mathematics} \textbf{378}, CUP, 2010.}%
\newline
[BinO2]{\ N. H. Bingham and A. J. Ostaszewski, {Normed groups: Dichotomy and
duality.} \textsl{Dissert. Math. }\textbf{472} (2010), 138p.}\newline
[BinO3] {N. H. Bingham and A. J. Ostaszewski, {Dichotomy and infinite
combinatorics: the theorems of Steinhaus and Ostrowski. }\textsl{Math. Proc.
Camb. Phil. Soc. }\textbf{150} (2011), 1-22.}\newline
[BinO4]{\ N. H. Bingham and A. J. Ostaszewski, {Category-measure duality:
convexity, mid-point convexity and Berz sublinearity}, \textsl{Aequationes
Mathematicae}, \textbf{91} (2017), 801-836 (arXiv:1607.05750).}\newline
[BinO5]{\ N. H. Bingham and A. J. Ostaszewski, {Beyond Lebesgue and Baire
IV: Density topologies and a converse Steinhaus-Weil theor{em}.} \textsl{%
Topology and its Applications} \textbf{239} (2018), 274-292
(arXiv:1607.00031).}\newline
[BinO6]{\ N. H. Bingham and A. J. Ostaszewski, Additivity, subadditivity and
linearity: Automatic continuity and quantifier weakening. \textsl{Indag.
Math.} (N.S.) \textbf{29} (2018), 687--713. (arXiv 1405.3948v3).}\newline
[BinO7] N. H. Bingham and A. J. Ostaszewski, Set theory and the analyst,
\textsl{European J. Math.}, to appear, arXiv:1801.09149v2.\newline
[BinO8] {N. H. Bingham and A. J. Ostaszewski, Beyond Haar and
Cameron-Martin: the Steinhaus support. arXiv}1805.02325v2.\newline
[Bog1] V. I. Bogachev, \textsl{Gaussian Measures}, Math. Surveys \&
Monographs \textbf{62}, Amer Math Soc., 1998.\newline
[Bog2]{\ V. I. Bogachev, \textsl{Measure theory. }Vol. I, II.
Springer-Verlag, Berlin, 2007.}\newline
[Bog3] {V. I. Bogachev, \textsl{Differentiable measures and the Malliavin
calculus.} Math. Surv. Mon. \textbf{164}, American Math. Soc, 2010.}\newline
[Bor] {K. C. Border,\textsl{\ Fixed point theorems with applications to
economics and game theory.} Cambridge University Press, 2$^{\text{nd}}$ ed.
1989 (1$^{\text{st}}$ ed.1985).}\newline
[Bor] C. Borell, Gaussian Radon measures on locally convex spaces, \textsl{%
Math. Scand.} \textbf{36} (1976), 265-284.\newline
[Bou] {A. Bouziad, {Continuity of separately continuous group actions in }$p$%
{-spaces}.\textsl{\ {Topology Appl.} }\textbf{71} (1996), 119-124.}\newline
[Brz] J. Brzd\k{e}k, Subgroups of the group $\mathbb{Z}_{n}$ and a
generalization of the Go\l \k{a}b-Schinzel functional equation. \textsl{%
Aequat. Math.} \textbf{43} (1992), 59--71.\newline
[CarS] T. J. Carlson and S. G. Simpson, A dual form of Ramsey's theorem,
\textsl{Adv. Math.} \textbf{53} (1984),265-290{.}\newline
[Chr1] {J. P. R. Christensen, {On sets of Haar measure zero in abelian
Polish groups.} Proceedings of the International Symposium on Partial
Differential Equations and the Geometry of Normed Linear Spaces (Jerusalem,
1972).\textsl{\ Israel J. Math.} \textbf{13} (1972), 255--260 (1973).}%
\newline
[Chr2] {J. P. R. Christensen, \textsl{Topology and Borel structure.
Descriptive topology and set theory with applications to functional analysis
and measure theory.} North-Holland Mathematics Studies \textbf{10}, 1974.}%
\newline
[ChrH] {J. P. R. Christensen, W. Herer,{\ On the existence of pathological
submeasures and the construction of exotic topological groups. }\textsl{%
Math. Ann. }\textbf{213} (1975), 203--210.}\newline
[CieJ] {\ K. Ciesielski, J. Jasi\'{n}ski, {Topologies making a given ideal
nowhere dense or meager}. \textsl{Topology Appl. }\textbf{63} (1995),
277--298.}\newline
[Coo] M. J. P. Cooper, Dimension, measure and infinite Bernoulli
convolutions. \textsl{Math. Proc. Cambridge Phil. Soc.} \textbf{124} (1998),
135-149.\newline
[DieS] {\ J. Diestel, A. Spalsbury, \textsl{The joys of Haar measure.} Grad.
Studies in Math. \textbf{150}. Amer. Math. Soc., 2014.}\newline
[Dar] {U. B. Darji, {On Haar meager sets.} \textsl{Topology Appl.}\textbf{160%
} (2013), 2396--2400.}\newline
[Dod1] P. Dodos, Dichotomies of the set of test measures of a Haar-null set.
\textsl{Israel J. Math.} \textbf{144} (2004), 15-28.\newline
[Dod2] P. Dodos, On certain regularity properties of Haar-null sets. \textsl{%
Fund. Math.} \textbf{191} (2004), 97-109.\newline
[Dod3] P. Dodos, The Steinhaus property and Haar-null sets. \textsl{Bull.
Lond. Math. Soc.} \textbf{41} (2009), 377--384.\newline
[DodKT] P. Dodos, V. Kanellopoulos, K. Tyros, A density version of the
Carlson and Simpson theorem, \textsl{J. Eur. Math. Soc.} \textbf{16} (2014),
2097--2164.\newline
[Drew1] {L. Drewnowski, {Topological rings of sets, continuous set
functions, integration.} I, II, III. \textsl{Bull. Acad. Polon. Sci. S\'{e}%
r. Sci. Math. Astronom. Phys.} \textbf{20} (1972), 269--276; \textbf{20}
(1972), 277--286; \textbf{20} (1972), 439--445.}\newline
[Drew2] {\ L. Drewnowski, {On control submeasures and measures.} \textsl{%
Studia Math.} \textbf{50} (1974), 203--224.}\newline
[Dud] R. M. Dudley. \textsl{Real analysis and probability.} Revised ed.,
Cambridge Studies in Advanced Mathematics \textbf{74}. Cambridge University
Press, 2002 ({1$^{\text{st}}$ ed. }1989).\newline
[Eng] {R. Engelking, \textsl{General topology}{.} 2$^{\text{nd}}$ ed.,
Heldermann, 1989 (1$^{\text{st}}$ ed. PWN, 1977).}\newline
[Erd1] P. Erd\H{o}s, On a family of symmetri Bernoulli convolutions. \textsl{%
Amer. J. Math.} \textbf{61} (1939), 974-975.\newline
[Erd2] P.Erd\H{o}s, On the smoothness properties of Bernoulli convolutions.
\textsl{Amer. J. Math.} \textbf{62} (1940), 180-186.\newline
[FenMW] {\ Q. Feng, M. Magidor, H. Woodin, {Universally Baire sets of reals,}
in H. Judah, W. Just, H. Woodin (eds.), \textsl{Set theory of the continuum}%
, 203--242, Math. Sci. Res. Inst. Publ. \textbf{26}, Springer, 1992.}\newline
[FenN] {\ J. E. Fenstad, D. Normann, {On absolutely measurable sets.}
\textsl{Fund. Math.} \textbf{81} (1973/74), 91--98.}\newline
[For] {M. K. Fort, Jr., {A unified theory of semi-continuity.} \textsl{Duke
Math. J. }\textbf{16} (1949), 237--246.}\newline
[Fre] {D. Fremlin, \textsl{Measure theory}{\ Vol. 3: \textsl{Measure algebras%
}.} Corrected 2$^{\text{nd}}$ printing of the 2002 original. Torres Fremlin,
Colchester, 2004.}\newline
[FreNR] {\ D. Fremlin, T. Natkaniec, I. Rec\l aw, {Universally
Kuratowski-Ulam spaces.}\ \textsl{Fund. Math.} \textbf{165} (2000), 239--247.%
}\newline
[Ful] {R. V. Fuller, {Relations among continuous and various non-continuous
functions.} \textsl{Pacific J. Math.} \textbf{25} (1968), 495--509.}\newline
[Gao] {\ Su Gao, \textsl{Invariant descriptive set theory}{. }Pure and
Applied Mathematics \textbf{293}. CRC Press, 2009.}\newline
[GivH] {\ S. Givant and P. Halmos, \textsl{Introduction to Boolean algebras}%
, Springer 2009.}\newline
[Gli] I. Glicksberg, Some remarks on absolute continuity on groups. \textsl{%
Proc. Amer. Math. Soc.} \textbf{40} (1973), 135--139.\newline
[Gow1] C. Gowrisankaran, Radon measures on groups. \textsl{Proc. Amer. Math.
Soc.} \textbf{25} (1970), 381--384.\newline
[Gow2] C. Gowrisankaran, Quasi-invariant Radon measures on groups. \textsl{%
Proc. Amer. Math. Soc. }\textbf{35} (1972), 503--506.\newline
[GroE] {\ K.-G. Grosse-Erdmann, {An extension of the Steinhaus-Weil theorem.}
\textsl{Colloq. Math.} \textbf{57} (1989), 307--317.}\newline
[Hal] {P. R. Halmos, \textsl{Measure theory}, Grad. Texts in Math. \textbf{18%
}, Springer 1974 (1$^{\text{st}}$ ed. Van Nostrand, 1950).}\newline
[HewR] {\ E. Hewitt, K. A. Ross, \textsl{Abstract harmonic analysis}, Vol.
I, Grundl. math. Wiss. \textbf{115}, Springer 1963 [Vol. II, Grundl. \textbf{%
152}, 1970].}\newline
[Hey] H. Heyer, \textsl{Probability measures on locally compact groups.}
Ergebnisse Math. Grenzgebiete \textbf{94}, Springer, 1977.\newline
[HodHLS] {W. Hodges, I. Hodkinson, D. Lascar, S. Shelah, {The small index
property for }$\omega ${-stable }$\omega ${-categorical structures and for
the random graph.} \textsl{J. London Math. Soc.} \textbf{48} (1993),
204--218.}\newline
[HofT] {\ J. Hoffmann-J\o rgensen, F. Tops\o e, {Analytic spaces and their
application}, in [Rog, Part 3].}\newline
[HolN] {\ L. Hol\'{a}, B. Novotn\'{y}, {Subcontinuity}.\textsl{\ Math.
Slovaca} \textbf{62} (2012), 345--362.}\newline
[Hu] Z. Hu, A generalized Kakutani-Kodaira theorem.\textsl{\ Proc. Amer.
Math. Soc.} \textbf{133} (2005), 3437--3440{.}\newline
[Jab1] {\ E. Jab\l o\'{n}ska, {Some analogies between Haar meager sets and
Haar null sets in abelian Polish groups.}\textsl{\ J. Math. Anal. Appl. }%
\textbf{421} (2015), 1479--1486.}\newline
[Jab2] E. Jab\l o\'{n}ska, A theorem of Piccard's type in abelian Polish
groups.\textsl{\ Anal. Math.} \textbf{42} (2016), 159--164.\newline
[Jan] S. Janson, \textsl{Gaussian Hilbert spaces.} Cambridge Tracts in
Mathematics \textbf{129}. Cambridge University Press, 1997.\newline
[Kal] {\ O. Kallenberg, \textsl{Foundations of modern probability}{.} 2$^{%
\text{nd}}$ ed., Springer, 2002 (1$^{\text{st}}$ ed. 1997).}\newline
[Kec] {\ A. S. Kechris, \textsl{Classical descriptive set theory}{,}
Graduate Texts in Mathematics \textbf{156}, Springer, 1995.}\newline
[Kem] {J. H. B. Kemperman, }A general functional equation{. \textsl{Trans.
Amer. Math. Soc.} \textbf{86} (1957), 28--56.}\newline
[Kha] A. Kharazishvili, Some remarks on the Steinhaus property for invariant
extensions of the Lebesgue measure, \textsl{Eur. J. Math}., Online, 2018.%
\newline
[Kne] M. Kneser, Summenmengen in lokalkompakten abelschen Gruppen, \textsl{%
Math. Z.} \textbf{66} (1956), 88-110.\newline
[Kod] {K. Kodaira, {Uber die Beziehung zwischen den Massen und den
Topologien in einer Gruppe.} \textsl{Proc. Phys.-Math. Soc. Japan} (3)
\textbf{23} (1941), 67-119.}\newline
[Kom] Z. Kominek, On an equivalent form of a Steinhaus theorem, \textsl{%
Math. (Cluj)} \textbf{30} (53)(1988), 25-27.\newline
[Kuc] {\ M. Kuczma,\textsl{\ An introduction to the theory of functional
equations and inequalities. Cauchy's equation and Jensen's inequality}{.}
2nd ed., Birkh\"{a}user, 2009 (1st ed. PWN, Warszawa, 1985).}\newline
[KucS] M. Kuczma, J. Sm\'{\i}tal, On measures connected with the Cauchy
equation. \textsl{Aequationes Math.} \textbf{14} (1976), no. 3, 421--428.%
\newline
[Lar] Larsen, R. Measures with separable orbits. \textsl{Proc. Amer. Math.
Soc. }\textbf{19 (}1968), 569--572. \newline
[LeP] R. D. LePage, Subgroups of paths and reproducing kernels, \textsl{Ann.
Prob.} \textbf{1} (1973), 345-347.\newline
[Loo] {\ L. H. Loomis, \textsl{An introduction to abstract harmonic analysis.%
} Van Nostrand, 1953.}\newline
[LukMZ] {J. Luke\v{s}, Jaroslav, J. Mal\'{y}, L. Zaj\'{\i}\v{c}ek, \textsl{%
Fine topology methods in real analysis and potential theory.} Lecture Notes
in Math. \textbf{1189}, Springer, 1986.}\newline
[LiuR] T. S. Liu, A. van Rooij, Transformation groups and absolutely
continuous measures, \textsl{lndag. Math.} \textbf{71} (1968), 225-231.%
\newline
[LiuRW] T. S. Liu, A. van Rooij, J-K Wang, Transformation groups and
absolutely continuous measures II, \textsl{Indag. Math.} \textbf{73 }(1970),
57--61.\newline
[Mac] {G. W. Mackey, {Borel structure in groups and their duals.} \textsl{%
Trans. Amer. Math. Soc.} \textbf{85} (1957), 134--165.}\newline
[Mah] {D. Maharam, {An algebraic characterization of measure algebras}.
\textsl{Ann. of Math. }(2) \textbf{48} (1947). 154--167.}\newline
[MarS] {D. A. Martin, J. R. Steel, {Projective determinacy.} \textsl{Proc.
Nat. Acad. Sci. U.S.A}. \textbf{85} (1988), 6582--6586.}\newline
[Mic] {E. Michael, {Topologies on the spaces of subsets.} \textsl{Trans.
Amer. Math. Soc.} \textbf{71} (1951), 152-182.}\newline
[MilO] {H. I. Miller and A. J. Ostaszewski, {Group action and
shift-compactness.} \textsl{J. Math. Anal. App. }\textbf{392} (2012), 23-39.}%
\newline
[MonZ] {D. Montgomery, L. Zippin,\textsl{\ Topological transformation groups.%
} Krieger, 1974 (Interscience, 1955).}\newline
[Mos] {Y. V. Mospan, {A converse to a theorem of Steinhaus-Weil}. \textsl{%
Real An. Exch. }\textbf{31} (2005), 291-294.}\newline
[Mue] {B. J. Mueller, {Three results for locally compact groups connected
with the Haar measure density theorem.} \textsl{Proc. Amer. Math. Soc.}
\textbf{16} (6) (1965), 1414-1416.}\newline
[Mun] {J. R. Munkres, \textsl{Topology, a first course}, Prentice-Hall, 1975.%
}\newline
[Myc] {J. Mycielski, {Finitely additive measures.} \textsl{Coll. Math.,}
\textbf{42} (1979), 309-318.}\newline
[Ost1] {A. J. Ostaszewski, {Families of compact sets and their universals.}
\textsl{Mathematika} \textbf{21} (1974), 116--127.}\newline
[Ost2] {A. J. Ostaszewski, {Effros, Baire, Steinhaus-Weil and
non-separability.} \textsl{Topology and its Applications }(Mary Ellen Rudin
Memorial Volume) \textbf{195} (2015), 265-274.}\newline
[Oxt1] {\ Oxtoby, John C. {Invariant measures in groups which are not
locally compact.} \textsl{Trans. Amer. Math. Soc.} \textbf{60} (1946),
215--237.}\newline
[Oxt2] {J. C. Oxtoby, \textsl{Measure and category}, 2nd ed. Graduate Texts
in Math. 2, Springer, 1980 (1$^{\text{st}}$ ed. 1972).}\newline
[PacS] J. Pachl, J. Stepr\={a}ns, Continuity of convolution and SIN groups. $%
\QTR{sl}{Canad.Math.Bull.}$ \textbf{60} (2017), 845--854.\newline
[Par] {K. R. Parthasarathy,\textsl{\ Probability measures on metric spaces}.
Academic Press, 1967.}\newline
[Pat] {A. L. T. Paterson, \textsl{Amenability.} Math. Surveys and Mon.
\textbf{29}, Amer. Math. Soc., 1988}\newline
[Pet] {B. J. Pettis, {On continuity and openness of homomorphisms in
topological groups.} \textsl{Ann. of Math. }(2) \textbf{52} (1950), 293--308.%
}\newline
[Pic] {\ S. Piccard, {Sur les ensembles de distances des ensembles de points
d'un espace Euclidien.\ }\textsl{M\'{e}m. Univ. Neuch\^{a}tel} \textbf{13},
212 pp. 1939.}\newline
[Pro] V. Prokaj, A characterization of singular measures, \textsl{Real Anal.
Exchange} \textbf{29} (2003/2004), 805--812.\newline
[Rog] {\ C. A. Rogers, J. Jayne, C. Dellacherie, F. Tops\o e, J. Hoffmann-J%
\o rgensen, D. A. Martin, A. S. Kechris, A. H. Stone, \textsl{Analytic sets,}
Academic Press, 1980.}\newline
[Rud] W. Rudin, \textsl{Functional analysis}, 2$^{\text{nd}}$ ed.
McGraw-Hill, 1991 (2$^{\text{nd}}$ ed.).\newline
[Sak] S. Saks,\textsl{\ Theory of the integral. }2$^{\text{nd}}$ ed. Dover,
1964 (1$^{\text{st}}$ ed. 1937, Monografie Mat. \textbf{VII}).\newline
[Sch] L. Schwartz, \textsl{Radon measures on arbitrary topological spaces
and cylindrical measures.} Tata Institute of Fundamental Research Studies in
Mathematics \textbf{6}, Oxford University Press, 1973.\newline
[Sho] {\ R.M. Shortt, {Universally measurable spaces: an invariance theorem
and diverse characterizations}, \textsl{Fundamenta Math. }\textbf{121}
(1984), 169-176. }\newline
[ShiT] H. Shi, B. S. Thomson, Haar null sets in the space of automorphisms
on [0,1]. \textsl{Real Anal. Exchange} \textbf{24} (1998/99), 337--350.%
\newline
[Sim] {\ S. M. Simmons, {A converse Steinhaus-Weil theorem for locally
compact groups}. \textsl{Proc. Amer. Math. Soc.} \textbf{49} (1975), 383-386.%
}\newline
[Sol1] {\ S. Solecki, {Size of subsets of groups and Haar null sets}.
\textsl{Geom. Funct. Anal.} \textbf{15} (2005), 246--273.}\newline
[Sol2] {\ S. Solecki, {Amenability, free subgroups, and Haar null sets in
non-locally compact groups}. \textsl{Proc. London Math. Soc.} (3) \textbf{93}
(2006), 693--722.}\newline
[Sol3] {\ S. Solecki, {A Fubini theorem}. \textsl{Topology Appl. }\textbf{154%
} (2007), 2462--2464.}\newline
[Solo] {\ R. M. Solovay, {A model of set-theory in which every set of reals
is Lebesgue measurable.}}\textsl{\ Ann. of Math. }{(2) \textbf{92} (1970),
1--56.}\newline
[Solom] B. Solomyak, Notes on Bernoulli convolutions, in: Fractal geometry
and applications: a jubilee of Beno\^{\i}t Mandelbrot. \textsl{Proc. Sympos.
Pure Math.} \textbf{72}, Part 1, 207--230, Amer. Math. Soc., 2004.\newline
[Ste] {\ H. Steinhaus, {Sur les distances des points de mesure positive}.
\textsl{Fund. Math.} \textbf{1} (1920), 83-104.}\newline
[Str1] {\ R. Struble, {Almost periodic functions on locally compact groups}.%
\textsl{\ Proc. Nat. Acad. Sci. U. S. A.} \textbf{39} (1953). 122--126.}%
\newline
[Str2] {\ R. Struble, {Metrics in locally compact groups}. \textsl{%
Compositio Math.} \textbf{28} (1974), 217--222. }\newline
[Tal] {\ M. Talagrand, {Maharam's problem}. \textsl{Ann. of Math.} (2)
\textbf{168} (2008), 981--1009.}\newline
[Tam] K.W. Tam, On measures with separable orbit. \textsl{Proc. Amer. Math.
Soc.} \textbf{23} (1969) 409--411.\newline
[TomW] G. Tomkowicz, S. Wagon, \textsl{The Banach-Tarski paradox.} {%
Cambridge University Press, 2016 (1}$^{\text{st}}${\ \ ed. 1985).}\newline
[vDo] {\ E. van Dowen, {Fubini's theorem for null sets}. \textsl{Amer. Math.
Monthly} \textbf{96} (1989), 718-721.}\newline
[Web] {\ H. Weber, {FN-topologies and group-valued measures}. in E. Pap (ed)
\textsl{Handbook of measure theory}, Vol. I, 703--743, North-Holland, 2002.}%
\newline
[Wei] {A. Weil, \textsl{L'int\'{e}gration dans les groupes topologiques},
Actualit\'{e}s Scientifiques et Industrielles \textbf{1145}, Hermann, 1965 (1%
$^{\text{st }}$ ed. 1940).}\newline
[Woo] {\ W. H. Woodin,\ \textsl{The axiom of determinacy, forcing axioms,
and the nonstationary ideal.}{\ }2$^{\text{nd}}$ ed., De Gruyter Series in
Logic and its Applications\ \textbf{1}. De Gruyter, 2010.}\newline
[Xia] D. X. Xia,\textsl{\ Measure and integration theory on
infinite-dimensional spaces. Abstract harmonic analysis.} Pure and App.
Math. \textbf{48}. Academic Press, 1972.\newline
[Yam] Y. Yamasaki, \textsl{Measures on infinite-dimensional spaces.} World
Scientific, 1985.\newline
[Zak] {E. Zakon, {A remark on the theorems of Lusin and Egoroff}. \textsl{%
Canad. Math. Bull.}\textbf{\ 7} (1964), 291--295.}

\bigskip

\noindent Mathematics Department, Imperial College, London SW7 2AZ;
n.bingham@ic.ac.uk \newline
Mathematics Department, London School of Economics, Houghton Street, London
WC2A 2AE; A.J.Ostaszewski@lse.ac.uk

\end{document}